\documentclass{siamart250211}
\usepackage{amsmath,amsfonts,amsopn,amssymb}
\usepackage{graphicx}
\usepackage[caption=false]{subfig}
\usepackage{rotating}
\usepackage{multirow,relsize,makecell}
\usepackage{url}
\usepackage{cite}
\usepackage{booktabs}
\usepackage{float}
\usepackage{bm}
\usepackage{esint}
\usepackage{makecell}
\usepackage[most]{tcolorbox}
\ifpdf
\DeclareGraphicsExtensions{.eps,.pdf,.png,.jpg}
\else
\DeclareGraphicsExtensions{.eps}
\fi
\usepackage{algorithm, algorithmic}
\usepackage[shortlabels]{enumitem}
\usepackage{accents,stmaryrd}

\numberwithin{theorem}{section}
\numberwithin{figure}{section}
\numberwithin{table}{section}

\newsiamremark{assumption}{Assumption}
\newsiamremark{remark}{Remark}
\newsiamremark{question}{Question}
\newsiamremark{example}{Example}
\crefname{assumption}{Assumption}{Assumptions}
\crefname{remark}{Remark}{Remarks}
\crefname{example}{Example}{Examples}

\title{Long-time asymptotics of advection-diffusion of finite element differential forms}
\author{Kaibo Hu\footnotemark[1], Jindong Wang\footnotemark[1], and Umberto Zerbinati%
\thanks{Mathematical Institute, University of Oxford,
Oxford OX2 6GG, United Kingdom
(\email{kaibo.hu@maths.ox.ac.uk},
\email{jindong.wang@maths.ox.ac.uk},
\email{umberto.zerbinati@maths.ox.ac.uk}).}}

\newcommand{\HA}{\mathcal{H}_\theta^k}

\usepackage[normalem]{ulem}
\usepackage{color}

\begin{document}

\maketitle
\begin{abstract}
This paper establishes global numerical analysis for the long-time behavior of Lie advection--diffusion of differential forms discretized by finite element exterior calculus (FEEC). For potential flows we prove a weighted Hodge decomposition implying convergence to harmonic forms of dimension equal to the Betti number, yielding a differential-form antidynamo theorem. For general advection, a mixed FEEC scheme is introduced and shown to preserve key structural properties: the dimension of closed stationary discrete $k$-forms is at least $b_k$. Under verified spectral-preservation assumptions (satisfied for small $h$ by spectral approximation and for arbitrary $h$ by a lowest-order exponential-fitting $M$-matrix discretization), discrete antidynamo theorems are obtained. In particular, on two-dimensional manifolds the numerical $1$-forms converge exponentially to a uniquely determined closed stationary form for every diffusion coefficient $\eta>0$, rigorously excluding spurious dynamo growth. The results guarantee that FEEC faithfully inherits the topological long-time asymptotics of the continuous system.
\end{abstract}


\section{Introduction}\label{sec:intro}

Advection--diffusion of differential forms arises in a wide range of important phenomena, including scalar advection--diffusion, magnetohydrodynamics (MHD), and Fokker--Planck equations \cite{arnold1998topological,heumann2012fully,heumann2010eulerian}. These models admit the unified form
$$
\frac{\partial \omega}{\partial t}+L_{\beta}\omega+\eta \Delta_{\mathrm{HL}}\omega=0,
$$
where $\omega$ is a differential $k$-form, $L_{\beta}\omega$ denotes the Lie derivative of $\omega$ along a vector field $\beta$ modeling transport, $\Delta_{\mathrm {HL}}\omega:=(d^{k-1}\delta^{k}+\delta^{k+1} d^{k})\omega$ is the Hodge Laplacian modeling diffusion, and $\eta$ is a constant coefficient measuring the relative strength between advection and diffusion.

Understanding the long-time behavior of advection--diffusion of differential forms (referred to as {\it Lie advection-diffusion} below) is of fundamental importance. For densities, this is closely related to topological and Morse-theoretic aspects of Fokker--Planck and other nonlinear equations \cite[Section V.3]{arnold1998topological}. For MHD problems, asymptotics plays a central role in understanding force-free fields in solar physics and in the design of equilibrium fields for fusion devices \cite{pontin2020parker,boozer2019curl}, while long-time behavior is crucial for understanding dynamo action. Numerical computation has become indispensable for investigating these phenomena. Finite element methods have been developed for advection--diffusion of differential forms and full MHD problems \cite{heumann2012fully,heumann2010eulerian,hu2017stable,hu2021helicity,hu2019structure,heumann2015stabilized,hiptmair2018splitting,wang2024discontinuous,luo2025local,zhang2025structure,zhang2024decoupled}. Rigorous numerical analysis is therefore required to understand the reliability and interpretation of computational results.

Most existing works (such as the error analysis in \cite{beirao2025error}) focus on {\it local numerical analysis}, namely the convergence of numerical solutions as the mesh size $h\to0$ and/or the polynomial degree $p\to\infty$ for $p$/$hp$ methods, either at a fixed time or over a finite time interval. However, many problems arising in the applications above concern long-time asymptotics rather than finite-time behavior. An important example is dynamo theory, which studies the mechanism by which magnetic fields are generated and sustained by the motion of conducting fluids, with fundamental applications to geophysical and astrophysical systems such as planets and stars. Mathematically, the dynamo problem asks whether the magnetic field $\bm B$, viewed as a differential 2-form in 3D, can exhibit exponential growth under the evolution
\begin{equation}\label{eq:intro_induction_vector}
  \partial_t\bm B-\operatorname{curl}(\bm \beta\times \bm B)
  +R_m^{-1}\operatorname{curl}\operatorname{curl}\bm B=0,
  \qquad
  \operatorname{div}\bm B=0 .
\end{equation}
Here $\bm\beta$ is the prescribed velocity field and $R_m$ is the magnetic Reynolds number. A central question in dynamo theory is whether the transport term can overcome magnetic diffusion and generate exponentially growing magnetic fields, particularly in the large-$R_m$ regime \cite{arnold1998topological}. Dynamo theory remains an active area of research, and many fundamental questions are still open. Nevertheless, several antidynamo theorems identify situations in which dynamo action cannot occur \cite{arnold1998topological}. Antidynamo theory addresses the complementary question: which geometric, topological, or structural conditions rule out such growth?

Numerical methods are widely used to investigate dynamo and antidynamo phenomena. However, local numerical analysis is insufficient for addressing many questions of practical and theoretical interest in this setting. For example, convergence as $h\to0$ over a finite time interval does not determine whether numerical dynamo action (exponential growth of numerical solutions) will occur for a particular implementation with fixed discretization parameters. This calls for numerical analysis of global and long-time behavior. Similar global questions have already attracted attention in geometric numerical integration, for example in the study of KAM-type behavior and long-time preservation properties of symplectic integrators \cite{shang1999kam}.

Lie advection--diffusion is closely related to global, topological, and dynamical-system properties. In the spirit of structure-preserving numerical methods, we aim to formulate and analyze discrete analogues of such questions arising in topological hydrodynamics. One example of this paradigm is the long-time relaxation of MHD systems under topological constraints \cite{he2026helicity,farrell2026global} (see \cite[Chapter III]{arnold1998topological}). In the context of dynamo theory, several aspects interact. Spectral properties of the operators provide information about the evolution, while pseudospectra reveal additional dynamical features \cite{trefethen1997pseudospectra,trefethen2020spectra}. The convergence of eigenvalue problems and pseudospectra for Lie advection--diffusion as $h\to0$ has been studied in another work \cite{spectra-paper}. However, spectra and pseudospectra alone do not fully characterize long-time dynamics, and that analysis remains local in the sense that it studies convergence of (pseudo)spectra as $h\to0$. 
That analysis does not address the practically relevant situation of fixed discretization parameters ($h$ or $p$), where one wishes to determine whether a numerical scheme faithfully predicts dynamo or antidynamo behavior. For example, a computation at fixed discretization parameters may produce positive eigenvalues that later disappear under mesh refinement, leading to potentially misleading conclusions about long-time growth. This motivates the development of {\it global numerical analysis} for advection--diffusion of differential forms, which is the primary focus of this paper.

In this paper, we study properties of finite element methods for Lie advection--diffusion, particularly properties that do not require asymptotic restrictions on discretization parameters such as sufficiently small mesh size. In other words, in addition to asymptotic results as $h\rightarrow 0$, we focus on the behavior of numerical schemes at finite discretization parameters, which is the regime relevant in practical computation. We show that certain classes of numerical methods faithfully inherit structural properties of the underlying continuous dynamical systems arising in topological hydrodynamics. In particular, we establish antidynamo theorems for schemes based on finite element exterior calculus (FEEC) \cite{wu2020simplex,hu2017stable,arnold2006finite,arnold2018finite}. Consequently, these results help exclude undesirable situations in which the continuous system does not exhibit dynamo action while the numerical discretization spuriously predicts exponential growth.

More broadly, the discussion extends beyond antidynamo theorems and addresses advection--diffusion from a Hodge-theoretic viewpoint. Classical Hodge theory establishes that harmonic forms, namely the kernel of the Hodge Laplacian, have dimension equal to the Betti numbers of the underlying manifold. FEEC preserves this property and therefore provides a faithful discretization of Hodge theory. When advection is introduced, a natural question is how this correspondence changes. Topological hydrodynamics provides a partial answer \cite{arnold2013some,arnold1998topological}: when diffusion dominates, the dimension of the generalized harmonic forms, namely the kernel of the Lie advection--diffusion operator, still coincides with the Betti numbers, while in general the dimension is bounded below by the Betti numbers. Numerical evidence for this phenomenon was presented in \cite{spectra-paper}. In this paper, we show that FEEC discretizations faithfully preserve these structural properties.

Although our theory does not cover all cases of Lie advection--diffusion, paralleling the limitations of the continuous theory itself, a particularly transparent situation occurs when the drift is a potential flow. If there exists a smooth function $\theta$ such that $\bm\beta^\flat=\eta d\theta$, then the generalized Hodge Laplacian can be written as a weighted Hodge Laplacian, which is self-adjoint and positive semidefinite. This is closely related to Witten's deformation of the de Rham complex \cite{witten1982supersymmetry}. The weighted Hodge decomposition then implies that the space of weighted harmonic $k$-forms is naturally isomorphic to the $k$-th de~Rham cohomology group. Since the governing operator has non-negative real spectrum, the evolution preserves the weighted harmonic component contained in the kernel, while all orthogonal components decay to zero. Consequently, the long-time limit is completely determined by the weighted harmonic component of the initial data. This yields a differential-form version of an antidynamo mechanism for gradient-type drifts.

The paper is organized as follows. \Cref{sec:prelim} reviews notation from Hodge theory and FEEC and introduces the generalized Hodge Laplacian. Section~\ref{sec:weight} develops the weighted Hodge theory for potential flows, while \Cref{sec:discretization} presents the FEEC discretization. The finite element antidynamo theorems are established in \Cref{sec:anti}. \Cref{sec:verify} verifies the structural assumptions for sufficiently fine meshes and for lowest-order exponential-fitting discretizations on arbitrary mesh size. Numerical experiments are presented in \Cref{sec:experiments} to show the long-time evolution of the proposed discretizations. 
Finally, \Cref{sec:conclude} concludes the paper.

\section{Preliminaries}\label{sec:prelim}
This section introduces the necessary mathematical tools and notation. We begin with an overview of Hodge theory and its generalization through incorporating Lie advection. Following this, we analyze the spectral properties of the generalized operator to characterize its \emph{a priori} steady-state behavior. To conclude, we present finite element exterior calculus (FEEC), which serves as the discrete framework for our subsequent numerical analysis.

\subsection{Hodge theory and its generalization}
Unless otherwise stated, all metric dependent constructions and results are
understood with respect to a fixed Riemannian metric. Euclidean settings are specified explicitly when required. We first recall the basic notation from Hodge theory and the de~Rham complex. Given a compact, connected, oriented manifold $M$, possibly with boundary, let \(L^2\Lambda^k(M)\) denote the Hilbert space of
square-integrable differential \(k\)-forms. The \(L^2\) inner product is
\begin{equation}\label{eq:L2_inner_product}
  (u,v)
  :=
  \int_M u\wedge\star v,
  \qquad
  u,v\in L^2\Lambda^k(M),
\end{equation}
where \(\star:\Lambda^k(M)\to\Lambda^{n-k}(M)\) is the Hodge star. The
associated norm is denoted by
\[
  \|\omega\|^2
  =
  \|\omega\|_{L^2\Lambda^k}^2
  :=
  (\omega,\omega).
\]

For \(s\in\mathbb R\), we denote by \(H^s\Lambda^k(M)\) the Sobolev space of
\(k\)-forms whose coefficients belong to \(H^s(M)\), and we write its
norm as \(\|\cdot\|_{H^s}\). In particular, given exterior derivatives $d^k$ satisfying $d^k\circ d^{k-1}=0$, we set
\begin{equation}\label{eq:H_Lambda_space}
  V^k
  :=
  H\Lambda^k(M)
  :=
  \left\{
    \omega\in L^2\Lambda^k(M):
    d^k\omega\in L^2\Lambda^{k+1}(M)
  \right\},
\end{equation}
equipped with the norm
\begin{equation}\label{eq:H_Lambda_norm}
  \|\omega\|_{V^k}^2
  =
  \|\omega\|_{H\Lambda^k}^2
  :=
  \|\omega\|^2+\|d^k\omega\|^2 .
\end{equation}
We then have the $L^2$-based de~Rham domain complex
\begin{equation}
    0\longrightarrow H\Lambda^0(M)\overset{d}{\longrightarrow}H\Lambda^1(M)\overset{d}{\longrightarrow}\cdots \overset{d}{\longrightarrow}H\Lambda^n(M)\longrightarrow0.
\end{equation}
The coderivative \(\delta^k:\Lambda^k(M)\to\Lambda^{k-1}(M)\) is given by
\begin{equation}\label{eq:codifferential_definition}
\star  \delta^k
  =
  (-1)^{k} d^{n-k}\star .
\end{equation} 
In analogy with $H\Lambda^k(M)$, we define the space
\begin{equation}\label{eq:Hstar_space}
  H^*\Lambda^k(M)
  :=
  \left\{
    \omega\in L^2\Lambda^k(M):
    \delta^k\omega\in L^2\Lambda^{k-1}(M)
  \right\}.
\end{equation}
For smooth forms \(u\in\Lambda^{k-1}(M)\) and
\(v\in\Lambda^k(M)\), Stokes' theorem gives
\begin{equation}\label{eq:stokes_codifferential}
  (d^{k-1}u,v)
  =
  (u,\delta^k v)
  +
  \int_{\partial M}\operatorname{tr}(u)\wedge\operatorname{tr}(\star v).
\end{equation}
The spaces $H\Lambda^k$ and $H^*\Lambda^k$ with homogeneous boundary conditions are defined as
\begin{align}
  \mathring H\Lambda^k(M)
  &:=
  \{\omega\in H\Lambda^k(M):\operatorname{tr}\omega=0\},
  \label{eq:Hring_space}\\
  \mathring H^*\Lambda^k(M)
  &:=
  \{\omega\in H^*\Lambda^k(M):\operatorname{tr}(\star\omega)=0\}.
  \label{eq:Hstar_ring_space}
\end{align}
Throughout this paper, we work in settings in which the exterior
derivatives \(d^k\) in the underlying de~Rham Hilbert complex have
closed ranges. The exterior derivative is regarded as an unbounded operator with domain $ \mathcal D(d^k)=H\Lambda^k(M)$. The domain of the coderivative depends on the existence of the boundary. In the presence of boundary, we take $ \mathcal D(\delta^k)=\mathring H^*\Lambda^k(M)$, whereas in the boundaryless case, $  \mathcal D(\delta^k)=H^*\Lambda^k(M)$.
The space of harmonic \(k\)-forms is
defined by
\begin{equation}\label{eq:harmonic_space_domain}
\begin{aligned}
     \mathfrak H^k
  &=\left\{
    \omega\in\mathcal D(d^k):
    d^k\omega=0,\ (\omega,d^{k-1}\zeta)=0,\forall \zeta\in \mathcal D(d^{k-1})
  \right\}\\
  &=\left\{
    \omega\in\mathcal D(d^k)\cap\mathcal D(\delta^k):
    d^k\omega=0,\ \delta^k\omega=0
  \right\}.
\end{aligned}
\end{equation}
The Hodge Laplacian on \(k\)-forms is
\begin{equation}\label{eq:hodge_laplacian}
  \Delta_{\mathrm{HL}}^k
  :=
  \delta^{k+1}d^k+d^{k-1}\delta^k,
\end{equation}
with domain
\[
  \mathcal D(\Delta_{\mathrm{HL}}^k)
  =
  \left\{
    \omega\in\mathcal D(d^k)\cap\mathcal D(\delta^k):
    d^k\omega\in\mathcal D(\delta^{k+1}),\
    \delta^k\omega\in\mathcal D(d^{k-1})
  \right\}.
\]
The \(L^2\)-Hodge decomposition is the
orthogonal direct sum
\begin{equation}\label{eq:L2_hodge_decomposition}
  L^2\Lambda^k(M)
  =
  \mathcal R(d^{k-1})
  \oplus
  \mathcal R(\delta^{k+1})
  \oplus
  \mathfrak H^k ,
\end{equation}
where $ \mathfrak H^k=\mathcal N(\Delta_{\mathrm{HL}}^k)$. Here,  \(\mathcal N(T)\) and \(\mathcal R(T)\) denote the kernel and
range of an operator \(T\), respectively. The de Rham cohomology group associated with the complex is
\begin{equation}\label{eq:de_rham_cohomology}
  H_{\mathrm{dR}}^k(M)
  :=
  \mathcal N(d^k)/\mathcal R(d^{k-1}).
\end{equation}
Hodge theory identifies each de Rham cohomology class with a unique harmonic
representative. Equivalently, 
$  \mathfrak H^k
  \cong
  H_{\mathrm{dR}}^k(M)$.
Consequently, the \(k\)-th Betti number satisfies
\begin{equation}\label{eq:betti_harmonic_dimension}
  b_k(M)
  =
  \dim H_{\mathrm{dR}}^k(M)
  =
  \dim \mathfrak H^k .
\end{equation}

We now consider the transport of differential forms by a given velocity field \(\bm\beta\). The advection-diffusion evolution of \(\omega\in V^k\) with velocity field \(\bm\beta\) and diffusion coefficient \(\eta>0\) is given by
\begin{equation}\label{eq:ad}
  \frac{\partial \omega}{\partial t} + L_{\bm\beta} \omega +\eta \Delta_{\rm HL}\omega=0,
\end{equation}
where \(L_{\bm\beta}\) denotes the Lie derivative along \(\bm\beta\). By Cartan's formula,
\begin{equation}
  L_{\bm\beta}\omega = di_{\bm\beta}\omega+i_{\bm\beta}d\omega,
\end{equation}
where \(i_{\bm\beta}\) is the contraction with \(\bm\beta\).

Since the topological properties of the Hodge Laplacian
\(\Delta_{\mathrm{HL}}\) are well understood, it is natural to view
\begin{equation}\label{eq:generalized_hodge_laplacian}
  \Delta_{\eta,\bm\beta}
  :=
  \eta\Delta_{\mathrm{HL}}+L_{\bm\beta},
  \qquad \eta>0,
\end{equation}
as a generalized Hodge Laplacian. The operator
\(\Delta_{\eta,\bm\beta}\) is a second-order elliptic operator. On a compact manifold, elliptic regularity for a sufficiently large positive shift $\Delta_{\eta,\bm\beta}^k+cI$, together with the compact Sobolev embedding, implies that the shifted resolvent is compact on $L^2\Lambda^k(M)$. Hence $\Delta_{\eta,\bm\beta}^k$ has a compact resolvent. In particular, its spectrum consists of isolated eigenvalues of finite algebraic multiplicity, and therefore
$\mathcal N(\Delta_{\eta,\bm\beta}^k)$ is finite-dimensional.

The kernel $\mathcal N(\Delta_{\eta,\bm\beta}^k)$ therefore plays the role of a  space of generalized harmonic $k$-forms. Its elements are  the stationary solutions of the advection–diffusion evolution, i.e., the time-invariant part of the dynamics. However, unlike the classical Hodge Laplacian, the operator $\Delta_{\eta,\bm\beta}^k$ is, in general, non-self-adjoint and indefinite. As a result, the long-time behavior of the system may not converge to a steady state given by the kernel of $\Delta_{\eta,\bm\beta}$. This raises several natural questions: how does the generalized harmonic space $\mathcal N(\Delta_{\eta,\bm\beta}^k)$ reflect the topology of the underlying domain? Does the advection–diffusion evolution converge to this generalized harmonic space as $t\to\infty$? At the discrete
level, the first question is answered in \Cref{sec:weight} for potential
flows and in \cref{thm:betti} for general advection (a lower bound by the
Betti numbers), while the second is the subject of the finite element
antidynamo theorems of \Cref{sec:anti}.

\subsection{Spectral properties of $\Delta_{\eta,\beta}$}
Beyond characterizing the kernel of the generalized Hodge Laplacian operator
\(\Delta_{\eta,\bm\beta}^k\), we also require a priori information on its
spectral structure in order to understand the long-time dynamics of the
evolution. The spectrum is well-understood for the scalar cases ($k=0$ and $k=n$). For $0$-forms, the operator is a second-order elliptic operator, and obeys the maximum principle. For $n$-forms, the operator corresponds to the adjoint of the $0$-form operator via Hodge duality, subject to a reversed velocity field. The spectral behavior of $\Delta_{\eta,\bm\beta}$ on $0$-forms serves as the cornerstone for analyzing the global dissipation of the system. Rooted in the maximum and comparison principles, the operator exhibits the following fundamental spectral properties.

\begin{theorem}[spectral properties of $\Delta_{\eta,\bm\beta}$ for $k=0$] \label{thm:spec}
For any diffusion coefficient $\eta > 0$ and velocity field $\bm\beta \in L^{\infty}(M; TM)$, let $\sigma(\Delta_{\eta,\bm\beta}^0)$ denote the spectrum of the advection-diffusion operator of $0$-forms. $\Delta_{\eta,\bm\beta}^0$ satisfies the following properties:

\begin{enumerate}
  \item {principal eigenvalue:} There exists a real, algebraically simple principal eigenvalue $\lambda_0 \in \mathbb{R}$. In particular, $\lambda_0=0$ in the closed, periodic, and Neumann cases, while $\lambda_0>0$ under Dirichlet boundary conditions.
  \item {sectorial distribution:} The operator is sectorial. Its spectrum is confined within a sectorial region in the right half of the complex plane:
  \[
  \sigma(\Delta_{\eta,\bm\beta}^0) \subset \Sigma_{\vartheta, \lambda_0}
  := \left\{ z \in \mathbb{C} \;\big|\; |\arg(z - \lambda_0)| \le \vartheta \right\}
  \]
  for some angle $\vartheta \in [0, \pi/2)$.
  \item {spectral gap:} \(\Delta_{\eta,\bm\beta}^0\) has compact resolvent and there exists a constant $\mu$ such that all $\lambda \in \sigma(\Delta_{\eta,\bm\beta}^0) \setminus \{\lambda_0\}$ satisfy ${\rm Re}(\lambda) \ge\mu> \lambda_0$.
\end{enumerate}
\end{theorem}

The spectral properties stated in \cref{thm:spec} are well established (see, e.g., \cite{gilbarg1977elliptic,pazy2012semigroups,evans2022partial}). and to clarify the boundary-condition conventions used here, we include a proof for its simplified version in Appendix~\ref{apx:spectrum} for $\beta\in C^{0,1}(M;TM)$ in the closed, periodic, and Neumann cases.

The spectral structure described in \cref{thm:spec} is precisely the property that a
structure-preserving discretization should inherit. It will be
formalized at the discrete level as Assumption~\ref{asm:spectrum} in
\Cref{sec:anti} and verified for concrete schemes in \Cref{sec:verify}.

\subsection{Finite element exterior calculus}
We briefly introduce the basic setup for finite element exterior calculus (FEEC). Let \(\mathcal{T}_h\) be a shape-regular triangulation of the manifold \(M\), with mesh size \(h:=\max_{T\in\mathcal{T}_h}\operatorname{diam}(T)\), we have finite element subcomplexes $V_h^k \subset V^k$, which form the discrete de Rham complex
\begin{equation}
 0\longrightarrow V_h^0\overset{d}{\longrightarrow}V_h^1\overset{d}{\longrightarrow}\cdots \overset{d}{\longrightarrow}V_h^n\longrightarrow0.
\end{equation}
Standard choices for these discrete spaces include polynomial differential forms, $V_h^k = \mathcal{P}_r^-\Lambda^k(\mathcal{T}_h)$ and $\mathcal P_r\Lambda^k(\mathcal{T}_h)$.

Similarly, the discrete cohomology groups are defined as
\begin{equation}
  H_h^k:= (\mathcal N(d^k)\cap V_h^k)/(d^{k-1}V_h^{k-1}).
\end{equation}
A fundamental result in finite element exterior calculus (FEEC) is that the cohomology of the discrete complex is isomorphic to the cohomology of the continuous complex \cite{arnold2006finite,christiansen2010finite}. Consequently, $b_k = \dim H^k = \dim H_h^k$.

We introduce the discrete $L^2$ adjoint operator $\delta_h^k:V_h^k\rightarrow V_h^{k-1}$ defined by the following relation
\begin{equation}
  (\delta_h u_h, v_h) = (u_h,dv_h), \qquad \forall (u_h,v_h)\in V_h^{k}\times V_h^{k-1}.
\end{equation}
Additionally, we define the $L^2$ projection operator $Q_h^k: L^2\Lambda^k \rightarrow V_h^k$ through
\begin{equation}
  (Q_h^k u,v_h) = (u,v_h) , \qquad \forall v_h\in V_h^{k}.
\end{equation}

At the discrete level, one has
\begin{equation}
   \mathfrak B_h^k = d{V_h^{k-1}} ,\quad \mathfrak Z_h^k=\mathcal N(d^k)\cap V_h^k.
\end{equation}
The discrete harmonic form is defined by analogy as
\begin{equation}
   \mathfrak H_h^k = \{v\in\mathfrak Z_h^{k}: v\perp \mathfrak B^k_h\},
\end{equation}
which is isomorphic to $ H_h^k$. The discrete Hodge decomposition then takes the form
\begin{equation}
  V_h^k = \mathfrak Z^{k,\perp}_h \oplus {\mathfrak B}_h^k\oplus \mathfrak H_h^k = {\mathfrak B}_h^{k,*}\oplus \mathfrak B_h^k\oplus\mathfrak H_h^k,
\end{equation}
where ${\mathfrak B}_h^{k,*}=\delta_h^{k+1}{V_h^{k+1}}$.

\section{Weighted Hodge theory and antidynamo for potential flows}\label{sec:weight}
Building on the general framework of \Cref{sec:prelim}, we now examine a
case in which the generalized Hodge Laplacian recovers a self-adjoint structure. The advection field is assumed to be a potential field at the scale determined
by the diffusion coefficient: there exists a scalar function \(\theta\)
such that
\begin{equation}\label{eq:potential_flow}
  \bm\beta^\flat=\eta d\theta.
\end{equation}
For a vector field \(\bm\beta\), let $j_{\bm\beta}:=\bm\beta^\flat\wedge : \Lambda^k(M)\to \Lambda^{k+1}(M)$. Under the standard \(L^2\) inner product, \(j_{\bm\beta}\) is the formal adjoint
of the contraction \(i_{\bm\beta}\). Define the deformed exterior derivative
\begin{equation}\label{eq:dtheta_def}
  d_\theta\omega:=e^{-\theta}d(e^\theta\omega)
  =d\omega+d\theta\wedge\omega.
\end{equation}
Then \eqref{eq:potential_flow} gives
\begin{equation}\label{eq:eta_dtheta_flux}
  \eta d_\theta=\eta d+j_{\bm\beta},
\end{equation}
which identifies the advection--diffusion flux \(\eta d^k+j_{\bm\beta}^k\) and motivates the exponential-fitting flux used in the finite element discretization introduced later.

Let \(\theta\in W^{1,\infty}(M)\). Then there exists constants $0<E_0\le E^0<\infty$ such that
\begin{equation}\label{eq:weight_bounds}
  0<E_0\le e^{-\theta(x)}\le E^0<\infty,
  \qquad \text{ for a.e. } x\in M.
\end{equation}
Indeed, one may take $E_0=e^{-\|\theta\|_{L^{\infty}(M)}}$ and $E^0=e^{\|\theta\|_{L^{\infty}(M)}}$. We then define the corresponding weighted $L^2$ inner product.

\begin{definition}[weighted inner product]
The weighted \(L^2\) inner product on \(k\)-forms is
\begin{equation}\label{eq:weighted_inner_product}
  (\alpha,\gamma)_\theta
  :=\int_M e^{-\theta}\alpha\wedge\star\gamma.
\end{equation}
The induced norm is denoted by \(\|\omega\|_\theta^2=(\omega,\omega)_\theta\),
and the completion of \(C^\infty\Lambda^k(M)\) is denoted by
\(L^2_\theta\Lambda^k(M)\). By \eqref{eq:weight_bounds}, the weighted and
unweighted \(L^2\)-norms are equivalent.
\end{definition}

This weighted inner product determines the corresponding adjoint of the exterior derivative.
\begin{lemma}[weighted adjoint]\label{lem:adjoint}
The adjoint of \(d^k\) with respect to \((\cdot,\cdot)_\theta\) is
\begin{equation}\label{eq:weighted_codifferential}
  \delta_\theta^{k+1}:=e^\theta\delta^{k+1} e^{-\theta}, \quad \mathcal D(\delta_\theta^{k+1})=e^\theta\mathcal D(\delta^{k+1}).
\end{equation}
Equivalently, $\delta_\theta=\delta+i_{\nabla\theta}$.
\end{lemma}

\begin{proof}
Let \((d^k)^{*_\theta}\) denote the Hilbert adjoint of \(d^k\) with
respect to the weighted inner product. By definition,
\(\gamma\in\mathcal D((d^k)^{*_\theta})\) if and only if there exists
\(z\in L^2_\theta\Lambda^k(M)\) such that
\begin{equation}
  (d^k\alpha,\gamma)_\theta
  =
  (\alpha,z)_\theta,
  \qquad
  \forall\alpha\in\mathcal D(d^k).
\end{equation}
This condition is equivalent to
\begin{equation}
  (d^k\alpha,e^{-\theta}\gamma)
  =
  (\alpha,e^{-\theta}z),
  \qquad
  \forall\alpha\in\mathcal D(d^k).
\end{equation}
Since \(\delta^{k+1}=(d^k)^*\) in the standard \(L^2\) inner product,
the latter is equivalent to
\begin{equation}
  e^{-\theta}\gamma\in\mathcal D(\delta^{k+1}),
  \qquad
  \delta^{k+1}(e^{-\theta}\gamma)=e^{-\theta}z.
\end{equation}
Therefore,
\begin{equation}
  \gamma\in e^{\theta}\mathcal D(\delta^{k+1}),\qquad
  z
  =
  e^\theta\delta^{k+1}(e^{-\theta}\gamma).
\end{equation}
This proves both the action and the exact domain of the weighted
Hilbert adjoint.

Using the identity $i_X\omega=(-1)^{k(n-k)}\star (\star\omega\wedge X^\flat)=(-1)^{n(k+1)}\star (X^\flat\wedge\star\omega)$ (see \cite[Lemma 8.2.1]{hirani2003discrete}), we have, for
\(\gamma\in\Lambda^k\):
\begin{align*}
\delta(e^{-\theta}\gamma)
&=
(-1)^{n(k+1)+1}\star d(e^{-\theta}\star\gamma) \\
&=
e^{-\theta}\delta\gamma
+
(-1)^{n(k+1)+1}\star(de^{-\theta}\wedge\star\gamma) \\
&=
e^{-\theta}\delta\gamma
+
e^{-\theta}(-1)^{n(k+1)}
\star(d\theta\wedge\star\gamma) \\
&=
e^{-\theta}\bigl(\delta\gamma+i_{\nabla\theta}\gamma\bigr).
\end{align*}
Therefore,
\[
    \delta_\theta:=e^\theta\delta e^{-\theta}
    =
    \delta+i_{\nabla\theta}.
\]
\end{proof}

\begin{remark}[domain of the codifferential]
Regard the exterior derivative $d^{k-1}:
  \mathcal D(d^{k-1})=H\Lambda^{k-1}(M)
  \subset L^2\Lambda^{k-1}(M)
  \longrightarrow L^2\Lambda^k(M)$, its Hilbert adjoint $\delta^k:=(d^{k-1})^*$ has domain \cite{arnold2006finite}:
\[
  \mathcal D(\delta^k)
  =
  \begin{cases}
    H^*\Lambda^k(M),
      & \partial M=\emptyset,\\[1mm]
    \mathring H^*\Lambda^k(M),
      & \partial M\neq\emptyset,
  \end{cases}
\]
For the weighted codifferential $\delta_\theta^k$, given
$\theta\in W^{1,\infty}(M)$,
\[
\mathcal D(\delta_\theta^k)
=
e^\theta \mathcal D(\delta^k)
=
\mathcal D(\delta^k).
\]
Indeed, multiplication by $e^{\pm\theta}$ defines an isomorphism on
$H^*\Lambda^k(M)$ and, when $\partial M\neq\emptyset$, preserves the
homogeneous boundary condition defining
$\mathring H^*\Lambda^k(M)$.
\end{remark}

We can then define the associated weighted Hodge Laplacian and weighted harmonic forms.
\begin{definition}[weighted Laplacian and harmonic forms]
The weighted Hodge Laplacian on \(k\)-forms is
\begin{equation}\label{eq:weighted_laplacian}
  \Delta_\theta^k
  :=d^{k-1}\delta_\theta^k+\delta_\theta^{k+1}d^k.
\end{equation}
The space of weighted harmonic \(k\)-forms is
\begin{equation}\label{eq:weighted_harmonic_space}
  \HA:=\mathcal N(\Delta_\theta^k).
\end{equation}
\end{definition}

Because \(\bm\beta^\flat=\eta d\theta\), \cref{lem:adjoint} implies
\begin{equation}\label{eq:potential_operator_identity}
  \eta\Delta_\theta^k
  =\eta\Delta_{\mathrm{HL}}^k+d^{k-1}i_{\bm\beta}^k+i_{\bm\beta}^{k+1}d^k
  =\eta\Delta_{\mathrm{HL}}^k+L_{\bm\beta}^k.
\end{equation}
Thus the advection--diffusion operator with a potential-flow is exactly a weighted
Hodge Laplacian multiplied by \(\eta\).

\subsection{Weighted Hodge decomposition}
With these preparations, we now establish the weighted Hodge decomposition.
\begin{lemma}[equivalence and orthogonality]\label{lem:ortho}
A form \(\omega\in\HA\) if and only if \(d\omega=0\) and
\(\delta_\theta\omega=0\). Moreover,
\(\mathcal R(d^{k-1})\), \(\mathcal R(\delta_\theta^{k+1})\), and \(\HA\) are
mutually orthogonal with respect to \((\cdot,\cdot)_\theta\).
\end{lemma}

\begin{proof}
If \(\omega\in\HA\), then
\begin{align*}
  0=(\Delta_\theta^k\omega,\omega)_\theta
  &=\|\delta_\theta\omega\|_\theta^2+
   \|d\omega\|_\theta^2.
\end{align*}
Hence \(d\omega=0\) and \(\delta_\theta\omega=0\). The converse is immediate
from the definition of \(\Delta_\theta^k\). Orthogonality follows from the
weighted adjoint relation and from \(d\circ d=0\): for example,
\((d\alpha,\delta_\theta u)_\theta=(d\circ d\alpha,u)_\theta=0\). The orthogonality
with \(\HA\) follows similarly.
\end{proof}

\begin{remark}[implications of Witten deformation]\label{rmk:witten}
  Via the unitary transformation $U\omega = e^{-\theta/2}\omega$, the weighted Hodge Laplacian $\Delta_\theta$ satisfies
  \begin{equation}
  \begin{aligned}
       U \Delta_\theta U^{-1} &= e^{-\theta/2}(\delta_\theta d +d\delta_\theta) e^{\theta/2}\\
       &= (e^{\theta/2}\delta e^{-\theta/2})(e^{-\theta/2}de^{\theta/2})+(e^{-\theta/2}de^{\theta/2})(e^{\theta/2}\delta e^{-\theta/2})\\
       &=\delta_{\theta/2}d_{\theta/2}+d_{\theta/2}\delta_{\theta/2}\\
       &=(\delta_{\theta/2}+d_{\theta/2})^2
  \end{aligned}
  \end{equation}
  Hence, \(\Delta_\theta\) is unitarily equivalent to the Witten Laplacian $(d_{\theta/2}+\delta_{\theta/2})^2$, which is self-adjoint and non-negative.
\end{remark}

\begin{remark}[closedness of the ranges]
By the standard closed-range property of the de Rham complex, $\mathcal R(d^k)$ and $\mathcal R(\delta^{k+1})$ are closed. Because the weighted codifferential satisfies $\mathcal{D}(\delta_\theta^{k+1}) = e^\theta\mathcal{D}(\delta^{k+1})$, it immediately follows that $\mathcal{R}(\delta_\theta^{k+1}) = e^\theta\mathcal{R}(\delta^{k+1})$. Since multiplication by $e^\theta$ defines a bounded invertible operator on $L^2\Lambda$, it preserves closed subspaces. Thus, $\mathcal{R}(\delta_\theta^{k+1})$ is closed in $L^2\Lambda^{k}$. By the equivalence of the $L^2$- and $L^2_\theta$-norms, these ranges remain closed in the weighted space $L^2_\theta\Lambda^k$.
\end{remark}

\begin{theorem}[weighted Hodge decomposition]\label{thm:hodge}
The weighted Hilbert space admits the weighted orthogonal decomposition
\begin{equation}\label{eq:weighted_hodge_decomposition}
  L^2_\theta\Lambda^k(M)
  =\mathcal R(d^{k-1})
   \oplus \mathcal R(\delta_\theta^{k+1})
   \oplus \HA.
\end{equation}
\end{theorem}

\begin{proof}
We write \(\perp_\theta\) for the orthogonal complement with respect to the
weighted inner product \((\cdot,\cdot)_\theta\). By \cref{lem:adjoint},
\(\delta_\theta^{k+1}\) is the \(L^2_\theta\)-adjoint of \(d^k\). Hence, by
the closed range theorem and the closedness of the ranges,
\[
  \mathcal R(\delta_\theta^{k+1})
  =
  \mathcal N(d^k)^{\perp_\theta}.
\]
Therefore,
\begin{equation}\label{eq:first_weighted_split}
  L^2_\theta\Lambda^k(M)
  =
  \mathcal N(d^k)
  \oplus
  \mathcal R(\delta_\theta^{k+1}).
\end{equation}

It remains to decompose the closed subspace \(\mathcal N(d^k)\). Since
\[
  \mathcal R(d^{k-1})\subset \mathcal N(d^k)
\]
and \(\mathcal R(d^{k-1})\) is closed, we have the orthogonal decomposition
inside \(\mathcal N(d^k)\)
\[
  \mathcal N(d^k)
  =
  \mathcal R(d^{k-1})
  \oplus
  \left(
    \mathcal N(d^k)\cap \mathcal R(d^{k-1})^{\perp_\theta}
  \right).
\]
Again by \cref{lem:adjoint}, we have
\(
  \mathcal R(d^{k-1})^{\perp_\theta}
  =
  \mathcal N(\delta_\theta^k)
\).
Consequently,
\[
  \mathcal N(d^k)\cap \mathcal R(d^{k-1})^{\perp_\theta}
  =
  \mathcal N(d^k)\cap \mathcal N(\delta_\theta^k)
  =
  \mathcal H_\theta^k .
\]
Substituting this decomposition of \(\mathcal N(d^k)\) into
\eqref{eq:first_weighted_split} gives
\[
  L^2_\theta\Lambda^k(M)
  =
  \mathcal R(d^{k-1})
  \oplus
  \mathcal R(\delta_\theta^{k+1})
  \oplus
  \mathcal H_\theta^k .
\]
This proves the weighted Hodge decomposition.
\end{proof}

\begin{corollary}[cohomology isomorphism]\label{cor:weighted_cohomology}
The natural map
\[
  \Phi:\HA\to H^k_{\mathrm{dR}}(M),
  \qquad \Phi(\omega)=[\omega],
\]
is an isomorphism. In particular, \(\dim\HA=b_k\).
\end{corollary}

\begin{proof}
Let \([\omega]\in H^k_{\mathrm{dR}}(M)\) and choose a closed representative
\(\omega\). By Theorem~\ref{thm:hodge},
\[
  \omega=d\alpha+\delta_\theta\gamma+h_\theta,
  \qquad h_\theta\in\HA.
\]
Since \(d\omega=0\), we have \(d\delta_\theta\gamma=0\). Therefore
\[
  \|\delta_\theta\gamma\|_\theta^2
  =(\delta_\theta\gamma,\delta_\theta\gamma)_\theta
  =(d\delta_\theta\gamma,\gamma)_\theta=0,
\]
so \(\delta_\theta\gamma=0\). Hence \(\omega-h_\theta=d\alpha\), and
\([h_\theta]=[\omega]\). This proves existence.

If \(h_\theta,h_\theta'\in\HA\) represent the same cohomology class, then
\(h_\theta-h_\theta'=d\alpha\). Since \(\HA\perp_\theta\mathcal R(d)\), this
means $d\alpha=0$. Thus \(h_\theta=h_\theta'\), proving
uniqueness and the isomorphism.
\end{proof}

\begin{remark}
This weighted Hodge isomorphism was previously established in
\cite[Eq.~(2.13)]{lott2003some} for closed manifolds with a smooth positive
measure, which allows one to identify weighted harmonic one-forms with parallel
one-forms under nonnegative Bakry--Emery Ricci curvature.
\end{remark}

\subsection{Long-time behavior and the antidynamo mechanism}\label{subsec:weightlong}

Beyond the cohomology isomorphism, the weighted formulation also yields the
long-time behavior. Since \(\Delta_\theta^k\) is similar to a self-adjoint
non-negative operator with compact resolvent, it also has compact resolvent, and
its spectrum is real and non-negative. Decomposing the solution into its weighted harmonic and orthogonal components, the weighted harmonic component is time-invariant, while the remaining components dissipate. This gives the following corollary.

\begin{corollary}[antidynamo for potential flows]
Assume that $\bm\beta^\flat=\eta d\theta$.
Then the evolution
\[
    \partial_t\omega
    +
    \bigl(\eta\Delta_{\rm HL}+L_{\bm\beta}\bigr)\omega
    =
    0
\]
does not exhibit exponential growth. Moreover, the solution converges to its
weighted harmonic component as \(t\to\infty\).
\end{corollary}

\begin{remark}[boundary conditions]
We refer to \cite{al2012witten} for related discussions on Witten Laplacians with boundary conditions. Consider the de~Rham complex with homogeneous boundary conditions on \(d\), $\mathcal D(d^k)=\mathring H\Lambda^k(M)$ and $\mathcal D(\delta^k)=H^*\Lambda^k(M)$.  The same argument yields the corresponding weighted Hodge decomposition. In this
case, we have $\dim \mathring{\mathcal H}_\theta^{n-k}=b_k$ with $\mathring{\mathcal H}_\theta^{k}=\{\omega\in \mathring H\Lambda^k\cap (e^\theta H^*\Lambda^k):d^k\omega=0,\delta_\theta^k\omega=0\}$. 
\end{remark}

\subsection{Discrete weighted Hodge theory}\label{sec:weighted_feec} On the compact manifold $M$, we show that the weighted Hodge decomposition admits an analogue on the discrete complex $(V_h^\bullet,d)$ for \emph{every} mesh size $h$ and introduce a weighted Galerkin discretization of the advection--diffusion equation \eqref{eq:ad}. In this case, we prove that this scheme satisfies the antidynamo mechanism without any restriction on the mesh.

\begin{definition}[discrete weighted codifferential and harmonic forms]\label{def:discrete_weighted}
The discrete weighted codifferential $\delta_{\theta,h}^k:V_h^k\to V_h^{k-1}$ is defined by
\begin{equation}\label{eq:discrete_weighted_codiff}
  (\delta_{\theta,h}^k u_h,\tau_h)_\theta
  =(u_h,d^{k-1}\tau_h)_\theta,
  \qquad
  \forall\,\tau_h\in V_h^{k-1},
\end{equation}
i.e., $\delta_{\theta,h}^k$ is the adjoint of $d^{k-1}|_{V_h^{k-1}}$ with respect to the weighted inner product \eqref{eq:weighted_inner_product}. The space of discrete weighted harmonic $k$-forms is
\begin{equation}\label{eq:discrete_weighted_harmonic}
  \mathfrak H_{\theta,h}^k
  :=\bigl\{v_h\in V_h^k:\ d^kv_h=0,\ \delta_{\theta,h}^kv_h=0\bigr\}.
\end{equation}
\end{definition}

By \eqref{eq:discrete_weighted_codiff}, the condition $\delta_{\theta,h}^kv_h=0$ is equivalent to $(v_h,d^{k-1}\tau_h)_\theta=0$ for all $\tau_h\in V_h^{k-1}$. Therefore
\begin{equation}\label{eq:harmonic_characterization}
  \mathfrak H_{\theta,h}^k
  =\mathfrak Z_h^k\cap(\mathfrak B_h^k)^{\perp_\theta}.
\end{equation}
Since \eqref{eq:weight_bounds} guarantees that $(\cdot,\cdot)_\theta$ is an inner product on each finite-dimensional space $V_h^k$, the operator $\delta_{\theta,h}^k$ is well defined on $V_h^k$.

\begin{theorem}[discrete weighted Hodge decomposition]\label{thm:discrete_weighted_hodge}
For every $h>0$ and $0\le k\le n$, the discrete space admits the weighted orthogonal decomposition
\begin{equation}\label{eq:discrete_weighted_decomposition}
  V_h^k
  =\mathfrak B_h^k
  \oplus
  \mathfrak H_{\theta,h}^k
  \oplus
  \delta_{\theta,h}^{k+1}V_h^{k+1},
\end{equation}
where the sums are orthogonal with respect to $(\cdot,\cdot)_\theta$. Moreover, the natural map
\[
  \Phi_h:\mathfrak H_{\theta,h}^k\longrightarrow H_h^k,
  \qquad
  \Phi_h(v_h)=[v_h],
\]
is an isomorphism. In particular,
\begin{equation}\label{eq:discrete_weighted_dim}
  \dim\mathfrak H_{\theta,h}^k=\dim H_h^k=b_k .
\end{equation}
\end{theorem}

\begin{proof}
Since $\delta_{\theta,h}^{k+1}$ is the $(\cdot,\cdot)_\theta$-adjoint of the linear map $d^k:V_h^k\to V_h^{k+1}$ between finite-dimensional spaces, the standard duality gives
\[
  \mathcal R(\delta_{\theta,h}^{k+1})
  =\bigl(\mathcal N(d^k)\cap V_h^k\bigr)^{\perp_\theta}
  =(\mathfrak Z_h^k)^{\perp_\theta},
\]
whence
\begin{equation}\label{eq:first_discrete_split}
  V_h^k=\mathfrak Z_h^k\oplus_\theta\mathcal R(\delta_{\theta,h}^{k+1}).
\end{equation}
For $\mathfrak Z_h^k$, we split $\theta$-orthogonally along the subspace $\mathfrak B_h^k\subset\mathfrak Z_h^k$ and use \eqref{eq:harmonic_characterization}:
\begin{equation}\label{eq:second_discrete_split}
  \mathfrak Z_h^k
  =\mathfrak B_h^k\oplus
  \bigl(\mathfrak Z_h^k\cap(\mathfrak B_h^k)^{\perp_\theta}\bigr)
  =\mathfrak B_h^k\oplus\mathfrak H_{\theta,h}^k .
\end{equation}
Substituting \eqref{eq:second_discrete_split} into \eqref{eq:first_discrete_split} proves \eqref{eq:discrete_weighted_decomposition}.

For the dimension count, \eqref{eq:second_discrete_split} yields
\[
  \dim\mathfrak H_{\theta,h}^k
  =\dim\mathfrak Z_h^k-\dim\mathfrak B_h^k
  =\dim H_h^k
  =b_k.
\]
The map $\Phi_h$ is well defined because every element of $\mathfrak H_{\theta,h}^k$ is closed, and it is injective because
\[
  \mathcal N(\Phi_h)
  =\mathfrak H_{\theta,h}^k\cap\mathfrak B_h^k
  =\{0\}
\]
by the $\theta$-orthogonality. Injectivity together with the equality of dimensions implies that $\Phi_h$ is an isomorphism.
\end{proof}

\begin{remark}[boundary conditions]\label{rmk:weighted_discrete_bc}
If $M$ has a boundary and the complex with boundary conditions is used, $\mathring {V}_h^k=V_h^k\cap\mathring H\Lambda^k(M)$, then the statement and proof of \cref{thm:discrete_weighted_hodge} are similar with the relative discrete cohomology in place of $H_h^k$. In particular, the discrete
weighted harmonic space is
\[
  \mathring{\mathfrak H}_{\theta,h}^{k}
  :=\bigl\{v_h\in\mathring V_h^k:\ d^kv_h=0,\
  \mathring\delta_{\theta,h}^kv_h=0\bigr\}
  =\mathring{\mathfrak Z}_h^k\cap
  \bigl(\mathring{\mathfrak B}_h^k\bigr)^{\perp_\theta}.
\] with $\mathring\delta_{\theta,h}^k:\mathring V_h^k\to
\mathring V_h^{k-1}$
\[
  (\mathring\delta_{\theta,h}^k u_h,\tau_h)_\theta
  =(u_h,d^{k-1}\tau_h)_\theta,
  \qquad\forall\,\tau_h\in\mathring V_h^{k-1}.
\] 
We can obtain $\dim\mathring{\mathfrak H}_{\theta,h}^{k}=\dim \mathring{\mathfrak Z}_h^k-\dim \mathring{\mathfrak B}_h^k=b_k(M,\partial M)=b_{n-k}(M)$.
\end{remark}

\subsection{A weighted Galerkin discretization and its long-time behavior}
By \cref{lem:adjoint}, the weighted codifferential $\delta_\theta$ is the
adjoint of $d$ with respect to the weighted inner product
$(\cdot,\cdot)_\theta$, and by \eqref{eq:potential_operator_identity} the
advection--diffusion operator for potential flows factorizes as
\[
  \eta\Delta_\theta^k
  =\eta\bigl(d^{k-1}\delta_\theta^k+\delta_\theta^{k+1}d^k\bigr)
  =\eta\Delta_{\mathrm{HL}}^k+L_{\bm\beta}^k .
\]
This suggests discretizing \eqref{eq:ad} by a Galerkin method posed in the weighted inner product, in which no contraction operator is required: find $(\sigma_h,u_h)\in V_h^{k-1}\times V_h^k$ such that, for $t>0$,
\begin{equation}\label{eq:wadd}
\begin{aligned}
  (\sigma_h,\tau_h)_\theta
  -\eta(u_h,d^{k-1}\tau_h)_\theta &=0,
    &&\forall\,\tau_h\in V_h^{k-1},\\
  (\partial_tu_h,v_h)_\theta
  +(d^{k-1}\sigma_h,v_h)_\theta
  +\eta(d^ku_h,d^kv_h)_\theta &=0,
    &&\forall\,v_h\in V_h^k.
\end{aligned}
\end{equation}
The first equation of \eqref{eq:wadd} states that $\sigma_h=\eta\,\delta_{\theta,h}^ku_h$. Eliminating $\sigma_h$ yields the finite-dimensional evolution
\begin{equation}\label{eq:wadds}
  \partial_tu_h+A_{\theta,h}^ku_h=0,
  \qquad
  A_{\theta,h}^k
  :=\eta\Delta_{\theta,h}^k
  :=\eta\bigl(d^{k-1}\delta_{\theta,h}^k+\delta_{\theta,h}^{k+1}d^k\bigr).
\end{equation}
Since $d^{k+1}d^k=0$, we have the following commutativity
\[
  d^kA_{\theta,h}^k
  =\eta\, d^k\delta_{\theta,h}^{k+1}d^k
  =\eta\bigl(d^k\delta_{\theta,h}^{k+1}+\delta_{\theta,h}^{k+2}d^{k+1}\bigr)d^k
  =A_{\theta,h}^{k+1}d^k .
\]

The following theorem is the discrete counterpart of the antidynamo mechanism for gradient advection fields established in Section \ref{subsec:weightlong} for any $h,\eta>0$.

\begin{theorem}[finite element antidynamo theorem for potential flows]\label{thm:weighted_antidynamo}
Let $M$ be a connected, compact, oriented manifold, and let
$\bm\beta/\eta=(d\theta)^\sharp$ be a potential advection field. Then, for any
$h,\eta>0$ and every $0\le k\le n$, the evolution \eqref{eq:wadds} has the
following properties.
\begin{enumerate}
  \item {(no discrete dynamo modes)} The operator $A_{\theta,h}^k$ is
  self-adjoint and positive semidefinite on
  $\bigl(V_h^k,(\cdot,\cdot)_\theta\bigr)$. Consequently
  $\sigma(A_{\theta,h}^k)\subset[0,\infty)$, the map
  $t\mapsto\|u_h(t)\|_\theta$ is non-increasing, and
  \[
    \|e^{-tA_{\theta,h}^k}\|_{\mathcal L(V_h^k,\|\cdot\|_\theta)}\le1,
    \qquad t\ge0 .
  \]
  In particular, no solution of \eqref{eq:wadds} exhibits exponential (or any) growth in the weighted norm.
  \item {(steady states)} Every stationary solution is closed, and
  \[
    \mathcal N(A_{\theta,h}^k)
    =\mathfrak H_{\theta,h}^k
    \cong H_h^k,
    \qquad
    \dim\mathcal N(A_{\theta,h}^k)=b_k .
  \]
  \item {(exponential convergence)} Let $P_{\theta,h}^k$ denote the
  $\theta$-orthogonal projection onto $\mathfrak H_{\theta,h}^k$. Then, for
  every initial state $u_h(0)$,
  \begin{equation}\label{eq:weighted_decay}
    \bigl\|u_h(t)-P_{\theta,h}^ku_h(0)\bigr\|_\theta
    \le
    e^{-\lambda_{1,h} t}\,
    \bigl\|u_h(0)-P_{\theta,h}^ku_h(0)\bigr\|_\theta,
    \qquad t\ge0,
  \end{equation}
  where $\lambda_{1,h}>0$ is the smallest nonzero eigenvalue of
  $A_{\theta,h}^k$. Therefore $u_h(t)$ converges to $P_{\theta,h}^ku_h(0)$
  exponentially. In particular, for $k=n$, the total mass $\int_M u_h(t)$ is conserved if the manifold is closed, and
  hence the limit is the unique element of the one-dimensional space
  $\mathfrak H_{\theta,h}^n$ with total mass $\int_Mu_h(0)$.
\end{enumerate}
\end{theorem}

\begin{proof}
By \eqref{eq:discrete_weighted_codiff}, for all $u_h,v_h\in V_h^k$,
\begin{equation}\label{eq:weighted_quadratic_form}
  (A_{\theta,h}^ku_h,v_h)_\theta
  =\eta\bigl[(\delta_{\theta,h}^ku_h,\delta_{\theta,h}^kv_h)_\theta
  +(d^ku_h,d^kv_h)_\theta\bigr],
\end{equation}
so $A_{\theta,h}^k$ is self-adjoint and positive semidefinite on
$\bigl(V_h^k,(\cdot,\cdot)_\theta\bigr)$. Taking $v_h=u_h$ in
\eqref{eq:weighted_quadratic_form} shows that $A_{\theta,h}^ku_h=0$ if and
only if $d^ku_h=0$ and $\delta_{\theta,h}^ku_h=0$, i.e.,
\[
  \mathcal N(A_{\theta,h}^k)=\mathfrak H_{\theta,h}^k
\]
by \eqref{eq:discrete_weighted_harmonic}. In particular, every stationary
solution is closed, so
$\mathcal N(A_{\theta,h}^k)\cap\mathfrak Z_h^k=\mathcal N(A_{\theta,h}^k)$.
The isomorphism $\mathfrak H_{\theta,h}^k\cong H_h^k$ and the dimension
count $\dim\mathfrak H_{\theta,h}^k=b_k$ are the content of
\cref{thm:discrete_weighted_hodge}. This proves item~2.

Since $A_{\theta,h}^k$ is self-adjoint, it admits a $\theta$-orthonormal eigenbasis
$\{\psi_j\}_{j=1}^{N}$ with eigenvalues $0\le\mu_1\le\mu_2\le\cdots\le\mu_N$,
where the zero eigenspace is $\mathfrak H_{\theta,h}^k$ and
$\lambda_{1,h}:=\min\{\mu_j:\mu_j>0\}$. The semigroup admits the spectral
representation
\begin{equation}\label{eq:spectral_representation}
  u_h(t)=e^{-tA_{\theta,h}^k}u_h(0)
  =\sum_{j=1}^{N} e^{-\mu_j t}\,(u_h(0),\psi_j)_\theta\,\psi_j .
\end{equation}
Since $e^{-\mu_j t}\le1$ for all $j$ and $t\ge0$, we obtain
$\|e^{-tA_{\theta,h}^k}\|_{\mathcal L(V_h^k,\|\cdot\|_\theta)}\le1$, and
consequently $t\mapsto\|u_h(t)\|_\theta$ is non-increasing. This proves
item~1. Moreover, $P_{\theta,h}^ku_h(0)$ collects exactly the terms in
\eqref{eq:spectral_representation} with $\mu_j=0$, which are invariant in
time. The remaining terms carry the factors $e^{-\mu_jt}\le
e^{-\lambda_{1,h}t}$, whence
\[
  \bigl\|u_h(t)-P_{\theta,h}^ku_h(0)\bigr\|_\theta^2
  =\sum_{\mu_j>0}e^{-2\mu_jt}\,\bigl|(u_h(0),\psi_j)_\theta\bigr|^2
  \le e^{-2\lambda_{1,h}t}
  \bigl\|u_h(0)-P_{\theta,h}^ku_h(0)\bigr\|_\theta^2 .
\]
This proves \eqref{eq:weighted_decay} and the exponential convergence in
item~3.

Since $d^nu_h\equiv0$, the second equation of \eqref{eq:wadd} reads
$(\partial_tu_h+d^{n-1}\sigma_h,v_h)_\theta=0$ for all $v_h\in V_h^n$.
As $\partial_tu_h+d^{n-1}\sigma_h\in V_h^n$, it vanishes identically, i.e.,
$\partial_t u_h=-d^{n-1}\sigma_h$. Hence, by Stokes' theorem on the closed
manifold,
\[
  \frac{d}{dt}\int_M u_h
  =
  -\int_M d^{n-1}\sigma_h
  =
  0 ,
\]
so the total mass is conserved.

It remains to identify the limit. Since $d^n=0$,
\cref{thm:discrete_weighted_hodge} reduces to
\[
  V_h^n=\mathfrak B_h^n\oplus_\theta\mathfrak H_{\theta,h}^n,
  \qquad
  \dim\mathfrak H_{\theta,h}^n=b_n=1,
\]
where $b_n=1$ because $M$ is closed and oriented. By Stokes' theorem and the natural isomorphism
\(H_h^n\cong H_{\mathrm{dR}}^n(M)\), integration induces an isomorphism
\[
\int_M:\mathfrak H_{\theta,h}^n\longrightarrow\mathbb R.
\]
Hence, by \eqref{eq:weighted_decay} and mass conservation, \(u_h(t)\)
converges to the unique element of \(\mathfrak H_{\theta,h}^n\) with
mass \(\int_Mu_h(0)\). This completes the proof.
\end{proof}

\begin{remark}
Since $A_{\theta,h}^k$ is self-adjoint in $(\cdot,\cdot)_\theta$, the
$\theta$-orthogonal projection $P_{\theta,h}^k$ coincides with the spectral
projection onto $\mathcal N(A_{\theta,h}^k)$. This is special to the
potential-flow setting: for general advection fields the discrete operator
is non-normal in $L^2$, and the long-time limit is given by a spectral projection instead; see Remark \ref{rmk:spectral_proj} and
\Cref{sec:anti}.
\end{remark}

\section{Discretization for general advection field}\label{sec:discretization}
In this section, we discretize the generalized Hodge Laplacian
\(\Delta_{\eta,\bm\beta}^k\) by a FEEC-compatible mixed formulation for general $\bm\beta$. Let
\(A^k:=\Delta_{\eta,\bm\beta}^k\), so that the continuous evolution is
\(\partial_t\omega+A^k\omega=0\). The only nonstandard term is the contraction
\(i_{\bm\beta}\), which does not generally map \(V_h^k\) into
\(V_h^{k-1}\). We therefore choose a family of bounded
operators
\[
  i_{\bm\beta,h}^k:V_h^k\longrightarrow V_h^{k-1},
  \qquad k=1,\ldots,n,
\]
whose concrete realization will be discussed in \Cref{sec:verify}.

Introducing the auxiliary \((k-1)\)-form
\[
  \sigma_h=\eta\delta_h^k u_h+i_{\bm\beta,h}^k u_h,
\]
we define the semidiscrete problem as follows: find
\((\sigma_h,u_h)\in V_h^{k-1}\times V_h^k\) such that, for \(t>0\),
\begin{equation}\label{eq:add}
\begin{aligned}
  (\sigma_h,\tau_h)
  -(i_{\bm\beta,h}^k u_h,\tau_h)
  -\eta(u_h,d^{k-1}\tau_h) &=0,
    &&\forall\tau_h\in V_h^{k-1},\\
  (\partial_tu_h,v_h)
  +(d^{k-1}\sigma_h,v_h)
  +\eta(d^ku_h,d^kv_h)
  +(i_{\bm\beta,h}^{k+1}d^ku_h,v_h) &=0,
    &&\forall v_h\in V_h^k.
\end{aligned}
\end{equation}
The initial value is \(u_h(0)=\omega_h(0)\), for example
\(u_h(0)=\pi_h^k\omega(0)\), where \(\pi_h^k\) is a cochain projection.

\begin{definition}[discrete stationary pair and closed stationary form]
A pair \((\sigma_h,u_h)\in V_h^{k-1}\times V_h^k\) is called a
\emph{discrete stationary pair} if it satisfies
\begin{equation}\label{eq:sadd}
\begin{aligned}
  (\sigma_h,\tau_h)
  -(i_{\bm\beta,h}^k u_h,\tau_h)
  -\eta(u_h,d^{k-1}\tau_h) &=0,
    &&\forall\tau_h\in V_h^{k-1},\\
  (d^{k-1}\sigma_h,v_h)
  +\eta(d^ku_h,d^kv_h)
  +(i_{\bm\beta,h}^{k+1}d^ku_h,v_h)&=0,
    &&\forall v_h\in V_h^k.
\end{aligned}
\end{equation}
The form \(u_h\) is called a \emph{closed discrete stationary \(k\)-form} if,
in addition, \(d^ku_h=0\).
\end{definition}

The first equation of \eqref{eq:add} is the weak identity
\[
  \sigma_h=i_{\bm\beta,h}^ku_h+\eta\delta_h^ku_h.
\]
Substitution into the second equation gives the operator form
\begin{equation}\label{eq:adds}
  \partial_tu_h+A_h^ku_h=0,
\end{equation}
where
\begin{equation}\label{eq:Ah_def}
\begin{aligned}
  A_h^ku_h
  &:=d^{k-1}\sigma_h+\bigl(\eta\delta_h^{k+1}+i_{\bm\beta,h}^{k+1}\bigr)d^ku_h \\
  &=d^{k-1}\bigl(i_{\bm\beta,h}^ku_h+\eta\delta_h^ku_h\bigr)
   +\bigl(\eta\delta_h^{k+1}+i_{\bm\beta,h}^{k+1}\bigr)d^ku_h.
\end{aligned}
\end{equation}
Thus
\[
  A_h^k=\eta\Delta_h^k+L_{\bm\beta,h}^k,
  \qquad
  \Delta_h^k=d^{k-1}\delta_h^k+\delta_h^{k+1}d^k,
\]
with the discrete Lie advection operator
\begin{equation}\label{eq:discrete_lie}
  L_{\bm\beta,h}^k
  :=d^{k-1}i_{\bm\beta,h}^k+i_{\bm\beta,h}^{k+1}d^k.
\end{equation}
We denote the closed stationary space by
\begin{equation}\label{eq:Kh_Sh}
  \mathcal S_h^k:=\mathcal N(A_h^k)\cap\mathfrak Z_h^k.
\end{equation}

Because \(d^{k+1}d^k=0\), \eqref{eq:discrete_lie} gives the cochain
commutation relation
\begin{equation}\label{eq:discrete_lie_commute}
\begin{aligned}
  d^kL_{\bm\beta,h}^k
  &=d^k\bigl(d^{k-1}i_{\bm\beta,h}^k+i_{\bm\beta,h}^{k+1}d^k\bigr) \\
  &=d^ki_{\bm\beta,h}^{k+1}d^k
   =\bigl(d^ki_{\bm\beta,h}^{k+1}+i_{\bm\beta,h}^{k+2}d^{k+1}\bigr)d^k
   =L_{\bm\beta,h}^{k+1}d^k.
\end{aligned}
\end{equation}
The discrete Hodge Laplacian also commutes with \(d\), hence
\begin{equation}\label{eq:Ah_commute}
  d^kA_h^k=A_h^{k+1}d^k.
\end{equation}

\begin{remark}
The commutation in \eqref{eq:discrete_lie_commute} requires one compatible
family \(\{i_{\bm\beta,h}^k\}_{k=1}^n\). In particular, the same algebraic
contraction must be used whenever the same form degree appears in adjacent
equations.
\end{remark}

In \cite{arnold1998topological}, a lower bound for the dimension of the space of closed stationary solutions was established at the continuous level using the index of $\Delta_{\eta,\bm\beta}$. We now prove the corresponding result for the FEEC discretization. Specifically, we show that the space of closed discrete stationary $k$-forms has dimension at least the $k$-th Betti number. This result is not limited to closed manifolds, but also applies to periodic domains and, more generally, to the setting of the finite element de Rham complex without boundary conditions.

\begin{theorem}\label{thm:betti}
The closed discrete stationary space satisfies
\[
  \dim\mathcal S_h^k\ge b_k.
\]
\end{theorem}

\begin{proof}
Restrict \(A_h^k\) to \(\mathfrak Z_h^k\). If \(\omega_h\in\mathfrak Z_h^k\),
then \(d^k\omega_h=0\), and \eqref{eq:Ah_def} gives
\[
  A_h^k\omega_h
  =d^{k-1}\bigl(i_{\bm\beta,h}^k\omega_h+\eta\delta_h^k\omega_h\bigr)
  \in\mathfrak B_h^k.
\]
Therefore \(A_h^k(\mathfrak Z_h^k)\subset\mathfrak B_h^k\). By the rank--nullity theorem, we obtain
\[
\begin{aligned}
  \dim\mathcal S_h^k
  &=\dim\mathcal N\bigl(A_h^k|_{\mathfrak Z_h^k}\bigr) \\
  &=\dim\mathfrak Z_h^k-
   \dim A_h^k\mathfrak Z_h^k \\
  &\ge \dim\mathfrak Z_h^k-\dim\mathfrak B_h^k
   =\dim(\mathfrak Z_h^k/\mathfrak B_h^k)=b_k.
\end{aligned}
\]
\end{proof}

This result shows that the dimension of \(\ker A_h^k\) is at least the dimension of the space of harmonic \(k\)-forms, suggesting that the kernel may contain additional elements.

\begin{remark}[Evolution with boundary conditions]
    If we consider the finite element de Rham complex with homogeneous boundary conditions, i.e., $\mathring{V}_h^k \subset \mathring{H}\Lambda^k$, we can similarly show that
    $$
    \dim \mathcal{S}_h^k = \dim \mathcal{N}(A_h^k|_{\mathring{\mathfrak{Z}}_h^k}) \ge \dim (\mathring{\mathfrak{Z}}_h^k/\mathring{\mathfrak{B}}_h^k) = b_{n-k},
    $$
    where $\mathring{\mathfrak{Z}}_h^k = \mathfrak{Z}_h^k \cap \mathring{H}\Lambda^k$ and $\mathring{\mathfrak{B}}_h^k = d\mathring{V}_h^{k-1}$. In either case, the key conclusion is that adding advection cannot reduce the dimension of kernel, namely,
\begin{equation}
    \dim\mathcal N(A_h^k) \ge   \dim\mathcal N(\eta\Delta_{\rm HL}).
\end{equation}
\end{remark}

We close this section with an elementary lemma on the long-time behavior
of finite-dimensional semigroups, which is an immediate consequence of the spectral decomposition of finite-dimensional semigroups (see \cite[Theorem~V.3.1]{engel2000one}) and will be used repeatedly in
\Cref{sec:anti} to pass from the spectral information to convergence statements. Its proof is
given in Appendix~\ref{apx:expA}.

\begin{lemma}\label{lm:expA}
Let $X$ be a finite-dimensional subspace of a Hilbert space $H$, and $A\colon X\to X$ a linear operator. Suppose $0$ is an algebraically simple eigenvalue of $A$, and all other eigenvalues $\lambda \in \sigma(A)\setminus\{0\}$ satisfy $\text{\rm Re}(\lambda)>0$. Let $P\colon X\to X$ denote the spectral projection onto the kernel of $A$. Then, for any exponential decay rate 
\[
0 < \gamma < \min \left\{ \text{\rm Re}(\lambda) \mid \lambda\in \sigma(A)\setminus\{0\} \right\},
\]
there exists a constant $C>0$ such that
\[
\|e^{-tA}-P\|_{\mathcal L(X)}\le Ce^{-\gamma t}, \qquad \forall t\ge 0.
\]
Consequently, $e^{-tA} \to P$ exponentially as $t \to \infty$.
\end{lemma}

\begin{remark}
The same argument also covers the case \(0\notin\sigma(A)\), upon
setting \(P=0\), and the semigroup $e^{-tA}$ decays exponentially to $0$.
\end{remark}

\begin{remark}\label{rmk:spectral_proj}
The projection $P$ in \cref{lm:expA} is the spectral projection onto
$\mathcal N(A)$ along the complementary $A$-invariant subspace $Y$
(which equals $\mathcal R(A)$ when $0$ is simple). In general, $P$ is
\emph{not} the orthogonal projection onto $\mathcal N(A)$: the two coincide
if and only if $\mathcal R(A)\perp\mathcal N(A)$, e.g., when $A$ is
self-adjoint. Concretely, if $0$ is an algebraically simple eigenvalue of $A$, let $r$ and $\ell$ denote the associated
right and left null vectors, i.e., $Ar=0$ and $A^*\ell=0$, where $A^*$ is
the adjoint of $A$ with respect to the given inner product. Then
\[
  Px=\frac{(x,\ell)}{(r,\ell)}\,r ,
\]
since $\ell\perp \mathcal R(A)$.
\end{remark}

\section{Discrete preservation of long-time asymptotics}
\label{sec:anti}

In this section, we investigate the long-time asymptotic behavior of the semidiscrete scheme \eqref{eq:add} and identify conditions under which the structural properties of the continuous dynamics are preserved at the discrete level. These properties include the characterization of steady states and long-time limits for differential forms of different degrees. As a principal consequence, we establish a finite element antidynamo theorem in two dimensions.

To streamline the presentation, we first consider closed, connected,
oriented manifolds and periodic domains, for which no boundary terms arise.
The modifications required for manifolds with boundary, including the
effects of natural and essential boundary conditions on the long-time
asymptotics, are discussed separately in
\Cref{sec:asymptoticwithboundary}. We begin with a concise summary of the
main results in \cref{tab:summary}, including the characterization of
steady states, the existence of long-time limits, the relevant diffusion
regimes, and the assumptions under which each result holds.

\begin{table}[htbp]
	\centering
	\scriptsize
	\renewcommand{\arraystretch}{1.25}
	\setlength{\tabcolsep}{3pt}
	\begin{tabular}{
			|>{\centering\arraybackslash}m{0.100\textwidth}
			|>{\centering\arraybackslash}m{0.085\textwidth}
			|>{\centering\arraybackslash}m{0.280\textwidth}
			|>{\centering\arraybackslash}m{0.400\textwidth}
			|>{\centering\arraybackslash}m{0.045\textwidth}|}
		\hline
		Thm.
		& Assump.
		& Steady-state characterization
		& Long-time limit as $t\rightarrow\infty$
		& $\eta$ \\
		\hline
		
		Thm.~\ref{thm:all}
        (all $k$)
		& \ref{asm:stab}
		&
		\(\mathcal N(A_h^k)=\mathcal S_h^k\cong H_h^k\), \(\dim\mathcal  S_h^k=b_k\).
		&
		Closed data $\xrightarrow{\exp}$ unique stationary representative in its cohomology class; arbitrary data $\xrightarrow{\exp}$ a closed steady state.
		&
		\(>\eta_0\) \\
		\hline
		Thm.~\ref{thm:k0}
		($k=0$)
		& \ref{asm:spectrum}
		&
		\(\mathcal N(A_h^0)=\mathcal S_h^0=\mathfrak H_h^0\), global constants.
		&
		\(\omega_h(t)\xrightarrow{\exp}\) constant
		&
		\(>0\) \\
		\hline
		
		Thm.~\ref{thm:kn}
		($k=n$)
		& \ref{asm:spectrum}
		&
		\(\dim \mathcal N(A_h^n)=\dim \mathcal S_h^n=1\).
		&
		\(\omega_h(t)\xrightarrow{\exp}\omega_\infty\in \mathcal S_h^n\) with same mass $\int_M\omega_h(0)$
		&
		\(>0\) \\
		\hline
		
		Thm.~\ref{thm:k1}
		($k=1$)
		& \ref{asm:spectrum}
		&
		\(\mathcal S_h^1\cong H_h^1\), \(\dim \mathcal S_h^1=b_1\), and \(\mathcal S_h^1\cap \mathfrak B_h^1=\{0\}\).
		&
		closed data \(\xrightarrow{\exp}\) unique stationary representative in its cohomology class
		&
		\(>0\) \\
		\hline
		
		Thm.~\ref{thm:kn-1}
		($k=n-1$)
		& \ref{asm:spectrum}
		&
		Stationary \((n-1)\)-forms are closed:
		\(\mathcal N(A_h^{n-1})\subset \mathfrak Z_h^{n-1}\).
		&
		\(d^{\,n-1}\omega_h(t)\xrightarrow{\exp}0\)
		&
		\(>0\) \\
		\hline
		Thm.~\ref{thm:antidynamo2d}
        
		\makecell{($n=2,$ \\ ~$k=1$)}
		& \ref{asm:spectrum}
		&
		The limit is a closed stationary \(1\)-form; if \(b_1=0\), the stationary limit is \(0\).
		&
		\(\omega_h(t)\xrightarrow{\exp}\omega_\infty\in \mathcal S_h^1\); no dynamo modes.
		&
		\(>0\) \\
		\hline
		
	\end{tabular}
	\caption{Summary of the long-time asymptotic results on closed manifolds. 
		Here \(\mathcal S_h^k=\mathcal N(A_h^k)\cap\mathfrak Z_h^k\).}
	\label{tab:summary}
\end{table}

The two assumptions referenced in \cref{tab:summary} are introduced at
the beginning of \Cref{sec:anti_all} and \Cref{sec:anti_endpoint},
respectively, and are verified for concrete discretizations in
\Cref{sec:verify}. The modifications required in the presence of
boundaries and boundary conditions are discussed in
\Cref{sec:asymptoticwithboundary}.

\subsection{Long-time behavior for sufficiently large diffusion coefficients}\label{sec:anti_all}
In this subsection we show the long-time behavior of all form degrees $0\le k\le n$ simultaneously, requiring the diffusion coefficient to dominate the
advection. The only structural condition is the following uniform
$L^2$-boundedness of the discrete contraction.

\begin{assumption}[\(L^2\)-boundedness of the discrete contraction]\label{asm:stab}
There exists a constant \(C_{\bm\beta}\), independent of \(h\) and \(\eta\), such
that
\begin{equation}\label{eq:stab_assumption}
  \|i_{\bm\beta,h}^k\omega_h\|
  \le C_{\bm\beta}\|\omega_h\|,
  \qquad
  \forall \omega_h\in V_h^k,\quad k=1,\ldots,n .
\end{equation}
\end{assumption}

Under this assumption, we prove that, for sufficiently large \(\eta\), the
stationary space is characterized by the discrete cohomology and the
corresponding evolution converges to a stationary solution as \(t\to\infty\). The analogous result at the continuous level was established in \cite[Chapter V, Theorem~3.6]{arnold1998topological}.

\begin{theorem}\label{thm:all}
Under Assumption~\ref{asm:stab}, there exists \(\eta_0>0\), independent of
\(h\), such that for every \(\eta>\eta_0\) the following statements hold for
all degrees \(k\):
\begin{enumerate}
  \item Every stationary solution is closed, i.e., \(\mathcal S_h^k=\mathcal N(A_h^k)\).
  \item \(\dim\mathcal S_h^k=b_k\), and the natural map
  \(\mathcal S_h^k\to H_h^k\) is an isomorphism.
  \item If the initial form is closed, then the solution exponentially converges to the unique stationary form in its discrete cohomology class.
  \item For an arbitrary initial form, the solution exponentially converges to a closed stationary form.
\end{enumerate}
\end{theorem}

\begin{proof}
The proof is structured into three parts based on the discrete Hodge decomposition.

{Part 1: Every stationary solution of \eqref{eq:sadd} is closed.} Let $\omega_h \in V_h^k$ satisfy $A_h \omega_h = 0$. Applying $d$ yields
\begin{equation}
 dA_h\omega_h= d(\eta \delta_h d\omega_h + i_{\bm\beta,h} d\omega_h) = 0.
\end{equation}
Setting $\zeta_h = d\omega_h \in \mathfrak{B}_h^{k+1}$ and taking the inner product with $\zeta_h$, we obtain
\begin{equation}
  (\zeta_h, d(\eta\delta_h\zeta_h + i_{\bm\beta,h}\zeta_h)) = (\delta_h\zeta_h, \eta\delta_h\zeta_h + i_{\bm\beta,h}\zeta_h) = 0.
\end{equation}
Under Assumption~\ref{asm:stab}, applying the Cauchy--Schwarz inequality together with the discrete Poincar\'e inequality for exact forms, $\|\zeta_h\| \le C_P \|\delta_h \zeta_h\|$ for some constant $C_P$, which follows from a standard duality argument; see also the generalized Gaffney inequality in \cite{he2018generalized}, we obtain
\begin{equation} \label{eq:proof_theta}
  \eta\|\delta_h\zeta_h\|^2 = -(\delta_h\zeta_h, i_{\bm\beta,h}\zeta_h) \le C_{\bm\beta} \|\delta_h\zeta_h\| \|\zeta_h\| \le C_{\bm\beta} {C_P} \|\delta_h\zeta_h\|^2.
\end{equation}
Let $\eta_0 := C_{\bm\beta} {C_P} $ (independent of $h$). For any $\eta > \eta_0$, inequality \eqref{eq:proof_theta} forces $\delta_h\zeta_h = 0$ and $\zeta_h=0$ due to discrete Poincar\'e inequality. Thus, all stationary solutions are closed forms, which proves item 1.

{Part 2: Cohomology isomorphism and closed-form dynamics.} Since $A_h \mathfrak{Z}_h^k \subset \mathfrak{B}_h^k \subset \mathfrak{Z}_h^k$, the space of discrete closed forms is invariant under the evolution \eqref{eq:add}. Let $\omega_h(t)$ be an evolving closed form. Using the discrete Hodge decomposition, we uniquely decompose $\omega_h(t) = p_h(t) + h_h(t)$. 
The evolution equation gives
\begin{equation}
  \partial_t \omega_h = -A_h \omega_h = -d(i_{\bm\beta,h} \omega_h + \eta \delta_h \omega_h).
\end{equation}
Since the right-hand side is an exact form, the harmonic component must remain invariant, i.e., $\partial_t h_h = 0$, which implies $h_h(t) \equiv h_h(0)$. This confirms that the evolution strictly preserves the discrete cohomology class.

The exact component $p_h \in \mathfrak{B}_h^k$ satisfies a linear dynamical system
\begin{equation} \label{eq:proof_p}
  \partial_t p_h = -d i_{\bm\beta,h} p_h - \eta d \delta_h p_h - d i_{\bm\beta,h} h_h.
\end{equation}
The linear part $A_1 p_h := -d(i_{\bm\beta,h} p_h + \eta \delta_h p_h)$ satisfies
\begin{equation}
    (A_1 p_h,p_h) = -\eta\|\delta_h p_h\|^2-(i_{\bm\beta,h}p_h,\delta_hp_h)\le -(\eta-C_\beta C_P)\|\delta_hp_h\|^2 \le -{\eta-\eta_0\over C_P^2}\|p_h\|^2,
\end{equation}
for $\eta > \eta_0=C_{\bm\beta} C_P$. This implies $A_1$ restricted to $\mathfrak{B}_h^k$ is strictly dissipative. Consequently, $A_1$ generates an exponentially stable semigroup and is invertible on $\mathfrak{B}_h^k$.

Setting $\partial_t p_h = 0$ in \eqref{eq:proof_p} demonstrates that each harmonic form $h_h \in \mathfrak{H}_h^k$ uniquely dictates a stationary exact component $p_\infty = A_1^{-1}(d i_{\bm\beta,h} h_h) \in \mathfrak{B}_h^k$. Conversely, by the discrete Hodge decomposition, any stationary closed form in $\mathcal{S}_h^k$ is uniquely decomposed into exact and harmonic components that satisfy this identical relationship. Consequently, the mapping $h_h \mapsto p_\infty + h_h$ establishes a linear isomorphism between the harmonic space $\mathfrak{H}_h^k$ and the space of stationary closed forms $\mathcal{S}_h^k$.

Therefore, a one-to-one correspondence is established between the discrete cohomology classes and the stationary forms, which yields $\dim \mathcal{S}_h^k = \dim \mathfrak{H}_h^k = b_k$ (proving item 2). Furthermore, the contracting dynamics of \eqref{eq:proof_p} guarantee that any evolving closed form asymptotically converges to this unique stationary form $\omega_\infty = p_\infty + h_h$ determined by its initial cohomology class $h_h$ (proving item 3).

{Part 3: Convergence to the closed form.} Consider an initial form, and decompose the solution as $\omega_h(t) = p_h(t) + h_h(t) + q_h(t)$, where $q_h \in \mathfrak{B}_h^{k,*}$.
Applying $d$ to the evolution equation and utilizing the commutativity $d A_h = A_h d$, we have
\begin{equation}
  \partial_t (d\omega_h) = -A_h (d\omega_h).
\end{equation}
Since $d\omega_h \in \mathfrak{B}_h^{k+1}$ is an exact form, a similar strict contraction argument from Part 2 gives that
\begin{equation}
 {1\over 2}  \frac{d}{dt}\|d\omega_h\|^2 \le -\frac{\eta-\eta_0}{C_P^2}\|d\omega_h\|^2.
\end{equation}
By Gronwall's inequality, we obtain
\[
  \|d\omega_h(t)\|^2
  \le
  e^{-\gamma t}\|d\omega_h(0)\|^2,
\]
where $\gamma=\frac{2(\eta-\eta_0)}{C_P^2}>0$. Because $dq_h = d\omega_h$, the Poincar\'e inequality $\|q_h\| \le C_P\|dq_h\|$ ensures that the co-exact component $q_h(t)$ vanishes exponentially as $t \to \infty$, i.e., $\|q_h(t)\|^2\le C_P^2 e^{-\gamma t}\|d q_h(0)\|^2$. Thus, the system asymptotically approaches the space of closed forms.

To prove the existence of the limit $\omega_h(\infty)$, let $P_{\mathfrak H_h^k}$ denote the $L^2$ orthogonal projection onto $\mathfrak H_h^k$, projecting the evolution equation yields
\begin{equation}
  \partial_t h_h = -P_{\mathfrak{H}_h^k} A_h \omega_h = -P_{\mathfrak{H}_h^k} (i_{\bm\beta,h} dq_h).
\end{equation}
Since \(P_{\mathfrak{H}_h^k} i_{\bm\beta,h}\) is a bounded linear operator in $L^2$, we get
\[
  \|\partial_t h_h(t)\|
  \le \|P_{\mathfrak{H}_h^k} i_{\bm\beta,h}\|\,\|d q_h(t)\|
  \le C e^{-{\gamma\over 2} t}.
\]
Hence
\[
  \int_0^\infty \|\partial_t h_h(t)\|\,dt
  \le C\int_0^\infty e^{-{\gamma\over 2} t}\,dt
  =\frac{C'}{\gamma}<\infty.
\]
Therefore \(\partial_t h_h\in L^1([0,\infty);\mathfrak H_h^k)\). For \(0\le s<t\),
\[
  h_h(t)-h_h(s)=\int_s^t \partial_\tau h_h(\tau)\,d\tau,
\]
and
\[
  \|h_h(t)-h_h(s)\|
  \le \int_s^t \|\partial_\tau h_h(\tau)\|\,d\tau
  \le C\int_s^t e^{-{\gamma\over2} \tau}\,d\tau
  \le \frac{C'}{\gamma}e^{-{\gamma\over2} s}.
\]
The right-hand side tends to \(0\) as \(s\to\infty\), uniformly for \(t>s\). Thus \(\{h_h(t)\}_{t\ge0}\) is a Cauchy family in the finite-dimensional Hilbert
space \(\mathfrak H_h^k\). Hence there exists
\(h_\infty\in\mathfrak H_h^k\) such that
\[
  h_h(t)\to h_\infty
  \qquad \text{as } t\to\infty.
\]
Equivalently, the limit can be written as the absolutely convergent integral
\[
  h_\infty
  =
  h_h(0)-\int_0^\infty P_{\mathfrak{H}_h^k} i_{\bm\beta,h}d q_h(\tau)\,d\tau.
\]
Moreover, the convergence is exponential
\[
  \|h_h(t)-h_\infty\|
  \le \int_t^\infty \|\partial_\tau h_h(\tau)\|\,d\tau
  \le \frac{C'}{\gamma}e^{-{\gamma\over2} t}.
\]
We next prove that the exact component \(p_h(t)\) also has a finite limit.
Recall that \(p_h(t)\in\mathfrak B_h^k\) satisfies
\[
  \partial_t p_h=A_1p_h+F(t),
\]
where \(A_1:\mathfrak B_h^k\to\mathfrak B_h^k\) is strictly dissipative, and
\[
  F(t)
  =
  -d i_{\bm\beta,h}h_h(t)
  -
  d i_{\bm\beta,h}q_h(t)
  -
  P_{\mathfrak B_h^k} i_{\bm\beta,h}d q_h(t).
\]
Since the mesh is fixed, all discrete operators involved are bounded. Moreover,
from the previous steps,
\[
  h_h(t)\to h_\infty,
  \qquad
  q_h(t)\to0,
  \qquad
  dq_h(t)\to0
\]
exponentially as $t\rightarrow\infty$. Hence
\[
  F(t)\to F_\infty:=-d i_{\bm\beta,h}h_\infty \in \mathfrak B_h^k.
\]
More precisely, there exist constants \(C_F>0\) and \(\kappa>0\) such that
\[
  \|F(t)-F_\infty\|\le C_F e^{-\kappa t}.
\]
Indeed,
\[
\begin{aligned}
  F(t)-F_\infty
  &=
  -d i_{\bm\beta,h}(h_h(t)-h_\infty)
  -
  d i_{\bm\beta,h}q_h(t)
  -
  P_{\mathfrak B_h^k} i_{\bm\beta,h}d q_h(t),
\end{aligned}
\]
and each term on the right-hand side decays exponentially by boundedness of the
operators. Since \(A_1\) is strictly contracting on the finite-dimensional space
\(\mathfrak B_h^k\), its semigroup satisfies
\[
  \|e^{tA_1}\|\le M e^{-\alpha t}
\]
for some \(M\ge1\) and \(\alpha>0\). In particular, \(A_1\) is invertible on $\mathfrak B_h^k$.

Define
\[
  p_\infty:=-A_1^{-1}F_\infty.
\]
Then \(p_\infty\in\mathfrak B_h^k\) and it solves the limiting steady equation
\[
  A_1p_\infty+F_\infty=0.
\]
Set
\[
  r_h(t):=p_h(t)-p_\infty.
\]
Using the equations for \(p_h(t)\) and \(p_\infty\), we obtain
\[
  \partial_t r_h
  =
  A_1r_h+\bigl(F(t)-F_\infty\bigr).
\]
By Duhamel's formula,
\[
  r_h(t)
  =
  e^{tA_1}r_h(0)
  +
  \int_0^t e^{(t-s)A_1}\bigl(F(s)-F_\infty\bigr)\,ds.
\]
Therefore
\[
\begin{aligned}
  \|r_h(t)\|
  &\le
  M e^{-\alpha t}\|r_h(0)\|
  +
  M\int_0^t e^{-\alpha(t-s)}
    \|F(s)-F_\infty\|\,ds \\
  &\le
  M e^{-\alpha t}\|r_h(0)\|
  +
  MC_F\int_0^t e^{-\alpha(t-s)}e^{-\kappa s}\,ds.
\end{aligned}
\]
More explicitly, if \(\alpha\neq\kappa\), then
\[
  \int_0^t e^{-\alpha(t-s)}e^{-\kappa s}\,ds
  =
  \frac{e^{-\kappa t}-e^{-\alpha t}}{\alpha-\kappa},
\]
and if \(\alpha=\kappa\), then it equals
\[
  t e^{-\alpha t}.
\]
Thus, for every \(0<\nu<\min\{\alpha,\kappa\}\), there exists
\(C_\nu>0\) such that
\[
  \|p_h(t)-p_\infty\|=\|r_h(t)\|\le C_\nu e^{-\nu t}.
\]
Thus
\[
  p_h(t)\to p_\infty
  =
  -A_1^{-1}F_\infty
  =
  A_1^{-1}d i_{\bm\beta,h}h_\infty.
\]
Combining this with \(h_h(t)\to h_\infty\) and \(q_h(t)\to0\), we conclude that
\[
  \omega_h(t)
  =
  p_h(t)+h_h(t)+q_h(t)
  \to
  p_\infty+h_\infty=:\omega_\infty.
\]
Finally, the limit exists and is closed and stationary:
$A_h\omega_\infty=A_h p_\infty+A_h h_\infty
=-A_1p_\infty+d\,i_{\bm\beta,h}h_\infty=F_\infty+d\,i_{\bm\beta,h}h_\infty=0$,
using $\delta_h h_\infty=0$, $dh_\infty=0$, and the definition of $p_\infty$. This proves item 4 and completes the proof.
\end{proof}

\subsection{Long-time behavior for arbitrary diffusion coefficients}\label{sec:anti_endpoint}
When the diffusion coefficient is arbitrary, the discretization requires
additional care. To this end, we assume the following discrete spectral
preservation property for the cases \(k=0\) and \(k=n\).

\begin{assumption}[discrete spectral preservation for \(0\)- and \(n\)-forms]\label{asm:spectrum}
For \(k=0\) and \(k=n\), assume that \(0\) is an algebraically simple eigenvalue of $A_h^k$ , and the nonzero eigenvalues lie in the open right half-plane, i.e.,
\[
  \operatorname{Re}\lambda>0,
  \qquad
  \forall \lambda\in\sigma(A_h^k)\setminus\{0\}.
\]
\end{assumption}

\begin{remark}
    The existence of the zero eigenvalue in fact holds unconditionally for every $h$ and $\eta$ (see the beginning of \Cref{sec:verify}). The substantive content of the assumption is the simplicity and the location of the remaining spectrum.
\end{remark}

We now study the long-time asymptotic behavior of the FEEC discretization of
the advection--diffusion equation and derive the corresponding finite element
antidynamo theorems. The results at the continuous level were established
in \cite[Chapter V, Theorems~3.7--3.12]{arnold1998topological}. At the discrete level, the
arguments follow a similar strategy, but the finite-dimensional setting provides
additional simplifications: the discrete operators are bounded, and finite element inverse inequalities imply, for instance, that asymptotic convergence of \(u_h\) also yields asymptotic convergence of \(d u_h\).

We first establish the long-time asymptotic behavior for the cases \(k=0\) and \(k=n\) under Assumption~\ref{asm:spectrum}.

\begin{theorem}[\(k=0\)]\label{thm:k0}
Under Assumption~\ref{asm:spectrum}, for every \(\eta>0\) the stationary
\(0\)-forms are exactly the constants,
\(\mathcal N(A_h^0)=\mathcal S_h^0=\mathfrak H_h^0\). Moreover, every solution of
\eqref{eq:adds} with initial data in \(V_h^0\) converges exponentially to a
constant.
\end{theorem}

\begin{proof}
For $k=0$, the operator reduces to $A_h = (i_{\bm\beta,h} + \eta \delta_h)d$. Since the exterior derivative $d$ annihilates constant functions, all constant discrete forms belong to $\mathcal S_h^0$. According to Assumption~\ref{asm:spectrum}, the eigenvalue $\lambda = 0$ is algebraically simple. Therefore, any stationary solution $\bar{\omega}_h \in \mathcal N(A_h)$ must be a constant, establishing $\mathcal N(A_h^0)=\mathcal S_h^0=\mathfrak H_h^0$. 

For an arbitrary initial condition $\omega_h(0)\in V_h^0$,
Lemma~\ref{lm:expA} yields the exponential convergence
\begin{equation}
  \lim_{t\to\infty}\omega_h(t)
  =\lim_{t\to\infty}e^{-tA_h^0}\omega_h(0)
  =P^{\mathrm{sp}}_0\,\omega_h(0),
\end{equation}
where $P^{\mathrm{sp}}_0$ is the spectral projection onto
$\mathcal N(A_h^0)=\mathfrak H_h^0$ along $\mathcal R(A_h^0)$
(cf.\ Remark~\ref{rmk:spectral_proj}). In
particular, the limit is a global constant. We emphasize that
$P^{\mathrm{sp}}_0$ differs in general from the $L^2$-orthogonal projection
onto constants. This completes the proof.
\end{proof}

\begin{theorem}[\(k=n\)]\label{thm:kn}
Under Assumption~\ref{asm:spectrum}, for every \(\eta>0\) the stationary
\(n\)-form space is one-dimensional:
\(\dim\mathcal N(A_h^n)=\dim\mathcal S_h^n=1\). Moreover:
\begin{enumerate}
  \item Every solution with initial data in \(V_h^n\) converges exponentially
  to the unique stationary \(n\)-form with the same total mass.
  \item If the initial condition is exact, \(\omega_h(0)=d^{n-1}B_h\) with
  \(B_h\in V_h^{n-1}\), then the solution converges exponentially to zero.
\end{enumerate}
\end{theorem}

\begin{proof}
For $k=n$ the operator reduces to
$A_h^n = d^{n-1}(i_{\bm\beta,h}^n + \eta\,\delta_h^n)$. Since every
$n$-form is closed, $\mathcal S_h^n=\mathcal N(A_h^n)\cap\mathfrak Z_h^n
=\mathcal N(A_h^n)$, and Assumption~\ref{asm:spectrum} gives
$\dim\mathcal S_h^n=1$. Lemma~\ref{lm:expA} then yields the exponential
convergence
\begin{equation}
  \lim_{t\to\infty}\omega_h(t)
  = P^{\mathrm{sp}}_n\,\omega_h(0),
\end{equation}
where $P^{\mathrm{sp}}_n$ is the spectral projection onto
$\mathcal N(A_h^n)$ along $\mathcal R(A_h^n)$
(cf.\ Remark~\ref{rmk:spectral_proj}).

First, Stokes' theorem implies strict mass conservation:
\begin{equation}\label{eq:massconserving}
  \frac{d}{dt}\int_M \omega_h = -\int_M d(i_{\bm\beta,h}\omega_h + \eta\delta_h\omega_h) = 0.
\end{equation}
Thus, $\int_M \omega_h(t) = \int_M \omega_h(0)$ for all $t \ge 0$. Because initial conditions can carry arbitrary non-zero mass, the basis of the one-dimensional stationary space $\mathcal{S}_h^n$ must inherently have a non-zero mass to preserve it.

Finally, consider an initially exact $n$-form $\omega_h(0) = d B_h$ with $B_h \in V_h^{n-1}$. By Stokes' theorem, its initial mass is zero ($\int_M d B_h = 0$). By mass conservation, it must exponentially converge to the unique stationary state carrying zero mass. Since any non-trivial stationary form carries non-zero mass, this limit must be strictly $0$. This completes the proof.
\end{proof}

The result for \(0\)-forms implies a corresponding statement for
closed \(1\)-forms.

\begin{theorem}[\(k=1\)]\label{thm:k1}
Under Assumption~\ref{asm:spectrum}, for every \(\eta>0\) the closed stationary
space satisfies
\[
  \dim\mathcal S_h^1=b_1(M),
  \qquad
  \mathcal S_h^1\cap\mathfrak B_h^1=\{0\}.
\]
Consequently, each class in \(H_h^1\) contains a unique closed stationary
\(1\)-form, and the evolution of any closed initial \(1\)-form converges
exponentially to that representative.
\end{theorem}

\begin{proof}
First, consider an initially exact $1$-form $\omega_h(0) = d u_h(0)$ with $u_h(0) \in V_h^0$. By the commutativity $A_h^1d=dA_h^0$, its evolution satisfies
\begin{equation}
  \omega_h(t) = e^{-tA_h^1} d u_h(0) = d \left( e^{-tA_h^0} u_h(0) \right).
\end{equation}
By Theorem~\ref{thm:k0}, the $0$-form evolution $e^{-tA_h^0}u_h(0)$ converges exponentially to a global constant $c$ as $t \to \infty$. Since $\omega_h(t)=d(e^{-tA_h^0}u_h(0))$ and the exact subspace is invariant, this implies $\lim_{t \to \infty}\omega_h(t)=dc=0$. This guarantees that any exact $1$-form dissipates to zero. Consequently, the only exact stationary $1$-form is the trivial zero field, i.e., $\mathcal{S}_h^1 \cap \mathfrak{B}_h^1 = \{0\}$.

By Theorem~\ref{thm:betti} and its proof, we have $\dim \mathcal{S}_h^1 \ge b_1(M) = \dim H_h^1$. Since $\mathcal{S}_h^1 \cap \mathfrak{B}_h^1 = \{0\}$, the natural projection from $\mathcal{S}_h^1$ onto the cohomology space $H_h^1 = \mathfrak{Z}_h^1 / \mathfrak{B}_h^1$ is injective. This forces $\dim \mathcal{S}_h^1 \le \dim H_h^1$. Hence, $\dim \mathcal{S}_h^1 =\dim H_h^1= b_1(M)$. This isomorphism dictates that each cohomology class uniquely corresponds to a closed stationary form, proving the dimension and uniqueness statements.

Finally, recall from the proof of \cref{thm:all} that the evolution of an initially closed $1$-form strictly preserves its discrete cohomology class. Since the natural projection from $\mathcal{S}_h^1$ to the cohomology space $H_h^1 = \mathfrak{Z}_h^1 / \mathfrak{B}_h^1$ is an isomorphism, the initial condition $\omega_h(0)$ shares its cohomology class with the unique stationary form $\bar{\omega}_h \in \mathcal{S}_h^1$. By the definition of cohomology classes, their difference must be an exact form, allowing us to uniquely decompose the initial condition as $\omega_h(0) = \bar{\omega}_h + \omega_h^\text{ex}(0)$, where $\omega_h^\text{ex}(0) \in \mathfrak{B}_h^1$. Since stationary forms are time-invariant and exact forms dissipate to zero, we obtain
\begin{equation}
  \lim_{t \to \infty} \omega_h(t) = \lim_{t \to \infty} e^{-tA_h^1}\bar{\omega}_h + \lim_{t \to \infty} e^{-tA_h^1}\omega_h^\text{ex}(0) = \bar{\omega}_h + 0 = \bar{\omega}_h.
\end{equation}
Thus, the solution asymptotically converges to the unique stationary form $\bar{\omega}_h$ within the cohomology class of the initial condition. This completes the proof.
\end{proof}

\begin{corollary}\label{cor:decomposition}
For \(k=1\), we have the decomposition
\[
  V_h^1=\mathcal S_h^1\oplus\mathfrak B_h^1\oplus\mathfrak B_h^{1,*}.
\]
\end{corollary}

\begin{proof}
The discrete Hodge decomposition implies
\[
  V_h^1=\mathfrak B_h^{1,*}\oplus\mathfrak B_h^1\oplus\mathfrak H_h^1.
\]
Theorem~\ref{thm:k1} shows that the natural projection
\(\mathcal S_h^1\to H_h^1=\mathfrak Z_h^1/\mathfrak B_h^1\) is an
isomorphism. Hence every closed form can be written uniquely as a sum of a
stationary closed form and an exact form, that is,
\(\mathfrak Z_h^1=\mathcal S_h^1\oplus\mathfrak B_h^1\). Since
\(\mathfrak B_h^{1,*}\) is the orthogonal complement of
\(\mathfrak Z_h^1\), the stated direct decomposition follows.
\end{proof}

The result for \(n\)-forms gives the complementary statement in degree
\(n-1\).

\begin{theorem}[\(k=n-1\)]\label{thm:kn-1}
Under Assumption~\ref{asm:spectrum}, for every \(\eta>0\) the exterior
derivative of every evolving \((n-1)\)-form decays exponentially:
\[
  d^{n-1}\omega_h(t)\longrightarrow0
  \qquad (t\to\infty).
\]
In particular, every stationary \((n-1)\)-form is closed, i.e., $\mathcal N(A_h^{n-1})\subset\mathfrak Z_h^{n-1}$.
\end{theorem}

\begin{proof}
Since $d$ and $A_h$ commute, the exact $n$-form $d\omega_h$ satisfies the $n$-form equation \eqref{eq:adds}. Theorem~\ref{thm:kn} implies that it dissipates to zero. Thus, we have
\begin{equation}
  \lim_{t \to \infty} d\omega_h(t) = 0,
\end{equation}
which proves that the exterior derivative of the evolving $(n-1)$-form asymptotically vanishes. 

If $\bar\omega_h$ is stationary for the $(n-1)$-form problem, then
\(A_h^n d\bar\omega_h=dA_h^{n-1}\bar\omega_h=0\). Moreover, \(d\bar\omega_h\) has zero total mass ($\int_M d\bar{\omega}_h=0$). By Theorem~\ref{thm:kn}, the only exact
stationary $n$-form with zero mass is zero. Therefore \(d\bar\omega_h=0\), so
every stationary $(n-1)$-form is closed.
\end{proof}

\subsection{The two-dimensional finite element antidynamo theorem}
\label{sec:anti_2d}
In dynamo theory, the central question is whether the evolution semigroup
\(e^{-tA}\) admits exponentially growing modes. Cowling's antidynamo
theorem \cite{cowling1933magnetic} rules out the maintenance of a
self-sustained axisymmetric magnetic field. It inspired a broad family of
antidynamo results \cite{zeldovich1957magnetic,roberts1992dynamo}, often
summarized by the principle that genuinely two-dimensional flows cannot
sustain fast kinematic dynamo action. Arnold placed this principle in a
geometric framework based on differential forms
\cite[Chapter~V]{arnold1998topological}. Our objective is to determine
whether FEEC discretizations \cite{arnold2006finite,arnold2010finite}
inherit this topological obstruction and thereby exclude spurious discrete
dynamo growth.

We now specialize the results of \cref{sec:anti_endpoint} to \(n=2\). In
this case, the degrees \(k=1\) and \(k=n-1\) coincide, allowing the two
corresponding conclusions to be combined. Consequently, for every
\(\eta>0\) and every fixed discretization satisfying
Assumption~\ref{asm:spectrum}, the discrete \(1\)-form solution converges
exponentially to a closed stationary \(1\)-form. This provides a discrete
FEEC counterpart of Arnold's two-dimensional antidynamo theorem: no
exponentially growing dynamo modes can be sustained by the discrete
evolution.

\begin{theorem}[discrete antidynamo] \label{thm:antidynamo2d}
Under Assumption~\ref{asm:spectrum} and $n=2$, for every
\(\eta>0\) and every initial discrete \(1\)-form, the solution of
\eqref{eq:adds} converges exponentially to a closed stationary \(1\)-form. Thus, there are no exponentially growing discrete dynamo modes. In particular, if \(M\) is simply
connected, the limiting stationary form is zero.
\end{theorem}

\begin{proof}
Using the discrete Hodge decomposition, we express the solution as
\begin{equation}
  \omega_h(t) = p_h(t) + h_h(t) + q_h(t),
\end{equation}
where $p_h\in \mathfrak B_h$, $h_h\in\mathfrak H_h$ and $q_h\in \mathfrak B_h^*$. We analyze these components successively. First, the coexact component
\(q_h\) is controlled through the identity \(d\omega_h=dq_h\). Its decay
then determines the evolution of the harmonic component \(h_h\), whose
time derivative is driven by \(dq_h\). Finally, the exact component
\(p_h\) satisfies a forced evolution equation whose forcing is determined
by \(h_h\), \(q_h\), and \(dq_h\).

For $n=2$, a $1$-form is also an $(n-1)$-form. By Theorem~\ref{thm:kn-1}, $d\omega_h(t) = dq_h(t) \to 0$ exponentially for any $\eta > 0$. The discrete Poincar\'e inequality then ensures the co-exact component $q_h(t) \in \mathfrak{B}_h^{1,*}$ also vanishes exponentially.

The argument proceeds similarly to the proof of \cref{thm:all}. Projecting the evolution onto $\mathfrak H_h^1$ gives
\[
  \partial_t h_h = -P_{\mathfrak H_h^1}\, i_{\bm\beta,h}^2\, d q_h.
\]
Since $dq_h\to0$ exponentially, the integral argument in the proof of
\cref{thm:all} (Part~3) yields $h_h(t)\to h_\infty\in\mathfrak H_h^1$
exponentially.

For the exact component \(p_h(t)\in\mathfrak B_h^1\), since the discrete Hodge decomposition is time-independent, we have
\[
  p_h(t)=P_{\mathfrak B_h^1}\omega_h(t),
  \qquad
  \partial_t p_h(t)
  =
  P_{\mathfrak B_h^1}\partial_t\omega_h(t).
\]
Projecting the evolution equation onto \(\mathfrak B_h^1\), and using
\(dp_h=0\), \(dh_h=0\), \(\delta_h h_h=0\), and
\(d\omega_h=dq_h\), gives
\begin{equation}\label{eq:ph_2d}
  \partial_t p_h=A_1p_h+F(t),
\end{equation}
where
\[
  A_1=-A_h^1|_{\mathfrak B_h^1}=-d(i_{\bm\beta,h}+\eta\delta_h),
\]
and
\[
  F(t)
  =
  -d i_{\bm\beta,h}h_h(t)
  -
  d i_{\bm\beta,h}q_h(t)
  -
  P_{\mathfrak B_h^1}i_{\bm\beta,h}dq_h(t).
\]
Since $N(A_h^0)=\mathbb R$ is invariant, $A_h^0$ induces an operator on
$V_h^0/\mathbb R$ whose spectrum is $\sigma(A_h^0)\setminus\{0\}$ ($0$ being
algebraically simple). The map $d^0:V_h^0/\mathbb R\to\mathfrak B_h^1$ is a
bijection and, by (4.8), the operator $A_h^1|_{\mathfrak B_h^1}$ is similar to $A_h^0|_{V_h^0/\mathbb R}$.
Hence $\sigma(A_h^1|_{\mathfrak B_h^1})=\sigma(A_h^0)\setminus\{0\}\subset\{\operatorname{Re}z>0\}$,
so $A_1=-A_h^1|_{\mathfrak B_h^1}$ generates an exponentially stable semigroup
and is invertible on $\mathfrak B_h^1$.

By Duhamel's
formula and the exponential stability of $A_1$, the same argument as in
the proof of \cref{thm:all} (Part~3) shows that
$p_h(t)\to p_\infty:=-A_1^{-1}F_\infty\in\mathfrak B_h^1$ exponentially for $F_\infty:=\lim_{t\to\infty}F(t)
=-d\,i_{\bm\beta,h}^1h_\infty$.

We conclude that
\[
  \omega_h(t)
  =
  p_h(t)+h_h(t)+q_h(t)
  \to
  p_\infty+h_\infty
  =:\omega_\infty.
\]
Moreover, $\omega_\infty$ is stationary. Since $q_h(t),\,dq_h(t)\to0$, the
forcing converges to $F_\infty$, and $p_\infty$ solves
$A_1p_\infty+F_\infty=0$. Using $dp_\infty=0$ and
$dh_\infty=0,\ \delta_hh_\infty=0$, we compute
\[
  A_h^1\omega_\infty
  =A_h^1p_\infty+A_h^1h_\infty
  =-A_1p_\infty+d^0 i_{\bm\beta,h}^1h_\infty
  =F_\infty - F_\infty=0 .
\]
Thus \(\omega_h(t)\) converges exponentially to the closed stationary
\(1\)-form \(\omega_\infty\).

If the manifold is simply connected, then $b_1(M)=0$ and $\mathfrak{H}_h^1$ is trivial, forcing $h_\infty = 0$ and $F_\infty = 0$. Consequently, $\omega_\infty = p_\infty=0$, completing the proof.
\end{proof}





\subsection{Long-time asymptotics under boundary conditions}\label{sec:asymptoticwithboundary}

The preceding results were stated for closed manifolds. We now describe how they are modified when \(M\) has a boundary. The arguments based only on the Hilbert-complex structure extend
to de~Rham complex with  absolute or relative cohomology after replacing the
corresponding cohomology spaces accordingly. The results
\cref{thm:k0,thm:kn}, however, require separate consideration, since the
boundary conditions may eliminate the zero eigenvalue in degree \(0\) or
\(n\).

Since both the continuous formulation and its FEEC discretization
are built upon the de Rham complex, it is natural to consider
boundary conditions that are compatible with this structure.
We therefore consider two standard types of boundary conditions:
absolute and relative boundary conditions.
In the mixed formulation, the absolute boundary conditions arise
naturally from the variational formulation, whereas the relative
boundary conditions are imposed strongly by restricting the
function spaces to the subcomplex with vanishing traces.
More precisely, we consider

\begin{itemize}
    \item[1.] {Absolute boundary conditions:}
    \[
        \operatorname{tr}(\star u)=0,
        \qquad
        \operatorname{tr}(\star du)=0,
    \]
    with the mixed variables satisfying
    \[
        (\sigma,u)\in
        H\Lambda^{k-1}(M)\times H\Lambda^k(M).
    \]

    \item[2.] {Relative boundary conditions:}
    \[
        \operatorname{tr}(u)=0,
        \qquad
        \operatorname{tr}(\sigma)=0,
    \]
    with the mixed variables satisfying
    \[
        (\sigma,u)\in
        \mathring H\Lambda^{k-1}(M)
        \times\mathring H\Lambda^k(M).
    \]
\end{itemize}

These choices are consistent with the classical absolute and
relative boundary conditions for the Hodge Laplacian and preserve
the corresponding de Rham complex structures.

\begin{remark}[Consistency with the cases $k=0$ and $k=n$]
	This is consistent with the classical scalar cases and the standard boundary conditions used to establish well-posedness. Assume that  $\partial M$ is smooth and
	$\bm\beta\in W^{1,\infty}(M;TM)$. In degree \(0\), the relative and absolute conditions reduce respectively to \( u|_{\partial M}=0, g( \nabla_g u,n)=0\) where \(n\) denotes the outward unit normal and $g$ is the Riemannian metric. Writing an \(n\)-form as \(v=\rho\,\operatorname{vol}_g\), where $\operatorname{vol}_g=\star 1$, the complementary absolute condition is \(\rho|_{\partial M}=0\).  The identity
\(\delta v=-i_{\nabla_g\rho}\operatorname{vol}_g\) gives $\sigma=\eta \delta v+i_{\bm\beta}v 
  =-i_{\eta\nabla_g\rho-\rho\bm\beta}\operatorname{vol}_g$. Since $\operatorname{tr}(i_X\operatorname{vol}_g)
  =g(X,n)\operatorname{vol}_{\partial M,g}$ (see \cite[Proposition~4.1.53]{nicolaescu2026lectures}), the relative conditions reduce to
\(g(\eta\nabla_g\rho-\rho\bm\beta,n)=0\), expressing vanishing
normal flux.
\end{remark}

\subsubsection{Long-time asymptotics under absolute boundary conditions}
We first consider the finite element de~Rham complex, where no essential
trace constraint is imposed:
\[
  V^k=H\Lambda^k(M),
  \qquad
  V_h^k\subset H\Lambda^k(M).
\]
The associated cohomology is the absolute de Rham cohomology with $\dim \mathfrak H_h^k=b_k(M)$. Accordingly, the general long-time result in \cref{thm:all} for sufficiently large diffusion coefficient $\eta$ remains
valid with \(b_k(M)\) interpreted as the absolute Betti number.

The main difference from the closed manifold case occurs in top
degree. If \(M\) is connected and \(\partial M\neq\emptyset\), then $b_n(M)=0$.
Thus the top-degree zero mode is absent. Assumption \ref{asm:spectrum} for $k=n$ should be $\operatorname{Re}(\lambda)>0,\forall \lambda\in\sigma(A_h^n)$. In
particular, the conclusion of \cref{thm:kn} should be replaced by $\mathcal S_h^n=\mathcal N(A_h^n)=\{0\}$,
and every solution of the \(n\)-form evolution decays exponentially to
zero. Equivalently, there is no nontrivial stationary \(n\)-form and total mass is not conserved in general under the natural boundary
condition. The remaining structural conclusions are unchanged after replacing the
cohomology spaces by the corresponding absolute cohomology spaces. 

\subsubsection{Long-time asymptotics under relative boundary conditions}
The essential boundary conditions require ${\rm tr}(u)={\rm tr }(\sigma)=0$, therefore we should consider the  finite element complex with homogeneous boundary conditions
\[
  \mathring{V}^k=\mathring H\Lambda^k(M),
  \qquad
  \mathring V_h^k\subset \mathring H\Lambda^k(M).
\]
The relevant cohomology is the relative cohomology
\(H^k(M,\partial M)\). We denote the corresponding discrete harmonic
space by \(\mathring{\mathfrak H}_h^k\) satisfying $ \dim \mathring{\mathfrak H}_h^k=b_k(M,\partial M)=b_{n-k}(M)$. Thus the dimension statement in \cref{thm:all} becomes
\[
  \dim \mathcal S_h^k
  =
  b_k(M,\partial M)
  =
  b_{n-k}(M),
\]
and the natural map $ \mathcal S_h^k\longrightarrow H_h^k(M,\partial M) := \mathring {\mathfrak Z}_h^k/ \mathring {\mathfrak B}_h^k$ is an isomorphism.

In degree \(0\), the essential boundary condition removes the constant
zero mode. Indeed, if \(M\) is connected and \(\partial M\neq\emptyset\),
then $ b_0(M,\partial M)=0$. The Assumption \ref{asm:spectrum} for $k=0$ should be $\operatorname{Re}(\lambda)>0,\forall \lambda\in\sigma(A_h^0)$.
Consequently, $\mathcal S_h^0=\mathcal N(A_h^0)=\{0\}$, and every solution of the \(0\)-form evolution decays exponentially to
zero. The conclusion of \cref{thm:k0} should be modified
accordingly under essential boundary conditions.

\cref{thm:kn} and \cref{thm:kn-1} remain valid. The corresponding result for closed stationary \(1\)-forms should also be interpreted in relative cohomology. More precisely,
\[
  \mathcal S_h^1\simeq \mathring{\mathfrak H}_h^1,
  \qquad
  \dim \mathcal S_h^1
  =
  b_1(M,\partial M)
  =
  b_{n-1}(M),
\]
and
\[
  \mathcal S_h^1\cap \mathring{\mathfrak B}_h^1=\{0\}.
\]
Hence each
relative cohomology class contains a unique closed stationary
\(1\)-form, and the evolution of any closed initial \(1\)-form converges
exponentially to that representative.

Finally, in two dimensions the degrees \(1\) and \(n-1\) coincide. Hence
the finite element antidynamo conclusion remains valid under either
natural or essential boundary conditions, with the corresponding
absolute or relative cohomology space: every discrete \(1\)-form
solution converges exponentially to a closed stationary \(1\)-form.

\section{Verification of the assumptions}\label{sec:verify}
In this section, we investigate the conditions under which the critical assumptions for the finite element antidynamo theorem are satisfied for closed manifolds. The first assumption is the $L^2$ boundedness of the discrete contraction operator. The second assumption involves the spectral preservation of the discretization. For sufficiently small $h$, this follows from spectral convergence under the usual approximation and regularity assumptions. We also explore whether an $M$-matrix discretization can be derived to guarantee the desired spectral distribution for arbitrary mesh size $h$. 
The results are summarized in \cref{tab:assumption_verification_summary}.

\begin{table}[htbp]
	\centering
	\scriptsize
	\renewcommand{\arraystretch}{1.25}
	\setlength{\tabcolsep}{4pt}
	\begin{tabular}{
			|>{\centering\arraybackslash}m{0.180\textwidth}
			|>{\centering\arraybackslash}m{0.080\textwidth}
			|>{\centering\arraybackslash}m{0.130\textwidth}
			|>{\centering\arraybackslash}m{0.35\textwidth}
			|>{\centering\arraybackslash}m{0.135\textwidth}|}
		\hline
		Verification choice
		& Verified assump.
		& Validity in \(h\)
		& Additional requirements
		& Used for \\
		\hline
		\(i_{\bm\beta,h}^k=Q_h^{k-1}i_{\bm\beta}^k\)
		&
		\ref{asm:stab}
		&
		arbitrary \(h\)
		&
		\(\bm\beta\in L^\infty\)
		&
		\cref{thm:all} \\
		\hline
		\(i_{\bm\beta,h}^k=Q_h^{k-1}i_{\bm\beta}^k\)
		&
		\ref{asm:spectrum}
		&
		\(h<h_0(\eta,\bm\beta)\)
		&
		\(\bm\beta\in L^\infty\);  $H^1(M)\hookrightarrow\hookrightarrow L^2(M)$; shape-regular mesh
		&
		\cref{thm:k0}--\ref{thm:antidynamo2d} \\
		\hline
		
		Exponential fitting \(M\)-matrix discretization
		&
		\ref{asm:spectrum}
		&
		arbitrary \(h\)
		&
		$\bm\beta\in W^{1,\infty}$; Euclidean metric; mesh condition in Assumption~\ref{ass:mesh0}
		&
		\cref{thm:k0}--\ref{thm:antidynamo2d} \\
		\hline
	\end{tabular}
	\caption{Choices used to verify the assumptions in Table~\ref{tab:summary}.}
	\label{tab:assumption_verification_summary}
\end{table}

\subsection{Existence of the zero eigenvalue}
First, we establish the existence of the zero eigenvalue at the discrete level for any $h>0$ and closed manifold $M$.

For $k=0$, the discrete space $V_h^0$ contains globally constant functions. Since the discrete exterior derivative satisfies $d^0 {1}=0$, it follows that $A_h^0{1}=0$. Thus, the zero eigenvalue is present for $0$-forms.

For $k=n$, let $m_h\in V_h^n$ be the Riesz representation of the integration functional, i.e.,
\begin{equation}
    (m_h,v_h)=\int_M v_h,\quad \forall v_h\in V_h^n.
\end{equation}
Since the integration functional is not zero, we have $m_h\neq 0$. The mass conservation property in \eqref{eq:massconserving} gives $\int_M A_h^n \omega_h = 0$ for all $\omega_h \in V_h^n$. Equivalently, 
\begin{equation}
   (m_h, A_h^n \omega_h)
   =
   \int_M A_h^n \omega_h
   =
   0,
   \qquad \forall\,\omega_h \in V_h^n .
\end{equation}
Thus
\(
  m_h\in \mathcal N\bigl((A_h^n)^*\bigr),
\)
where \((A_h^n)^*\) denotes the discrete adjoint of \(A_h^n\). Hence, \(A_h^n\) is
not surjective. Since \(A_h^n\) acts on the finite-dimensional space \(V_h^n\), it
is not invertible. Therefore \(\lambda=0\) is an eigenvalue of \(A_h^n\) for every mesh size \(h\).

Regarding Assumption~\ref{asm:spectrum}, while we have confirmed the existence of the zero eigenvalue for $k \in \{0, n\}$ for any $h>0$, it remains to be verified that this eigenvalue is algebraically simple and that the remaining spectrum resides strictly in the right half of the complex plane.

\subsection{Choice 1: spectrum convergence for sufficiently small $h$}

In this subsection, we verify Assumption~\ref{asm:stab} for every
\(h>0\) and Assumption~\ref{asm:spectrum} for sufficiently small \(h\).
We assume that \(\bm\beta\in L^\infty(M)\) and that $H^1(M)\hookrightarrow\hookrightarrow L^2(M)$. The Rellich compact embedding holds for closed manifolds and, more generally, for compact manifolds with Lipschitz boundary. For the discrete contraction, we define $i_{\bm\beta,h}^k = Q_h^{k-1} i_{\bm\beta}^k$, where $Q_h^{k-1}$ is the $L^2$-projection onto $V_h^{k-1}$. 

This is equivalent to applying the contraction operator $i_{\bm\beta}$ directly within the discrete variational formulation, since
\begin{equation}
  (i_{\bm\beta,h}^k u_h,v_h) = (Q_h^{k-1}i_{\bm\beta}^k u_h,v_h) = (i_{\bm\beta}^k u_h,v_h), \qquad \forall u_h\in V_h^k, \ v_h\in V^{k-1}_h.
\end{equation}
This choice naturally inherits $L^2$ boundedness, ensuring that Assumption~\ref{asm:stab} holds through the stability of the $L^2$-projection
\begin{equation}
  \|i_{\bm\beta,h}^ku_h\| = \|Q_h^{k-1}i_{\bm\beta}^ku_h\| \le \|i_{\bm\beta}^k u_h\| \le C_{\bm\beta}\|u_h\|,
\end{equation}
for $\bm\beta\in L^\infty$.

Regarding Assumption~\ref{asm:spectrum}, we have the following result.

\begin{lemma} \label{lm:spectrum_small_h}
  Let $\mathcal T_h$ be a shape regular family, and define $i_{\bm\beta,h}^k=Q_h^{k-1}i_{\bm\beta}^k$. Assume that $\bm\beta \in L^\infty(M)$. There exists a threshold $h_0 > 0$ depending on $\eta$ and $\bm\beta$ such that for all $h < h_0$, Assumption~\ref{asm:spectrum} is satisfied.
\end{lemma}

The proof is organized around a classical spectral convergence theorem for
families of operators, which we state first. Throughout,
$\mu=\mu_k>0$ denotes the spectral gap of the continuous operator
$A^k=\Delta^k_{\eta,\bm\beta}$ provided by \cref{thm:spec} (for $k=n$
through the duality with the $0$-form operator with reversed drift, whose
hypotheses are invariant under $\bm\beta\mapsto-\bm\beta$): the spectrum of
$A^k$ is contained in $\{0\}\cup\{\operatorname{Re}z\ge\mu\}$ and $0$ is an
algebraically simple eigenvalue. We consider a shift $c$ for $k\in\{0,n\}$ to be determined later
and write $\widetilde A^k:=A^k+cI$, $\widetilde A_h^k:=A_h^k+cI$.

The spectral-preservation result is based on the following classical
theorem for collectively compact operator approximations; see
\cite{anselone1968spectral}.

\begin{theorem}[spectral convergence under collective compactness]
\label{thm:abstract_cc}
Let $X$ be a complex Hilbert space, let $T\in\mathcal L(X)$ be compact, and
let $\{T_h\}_{0<h\le \bar h}\subset\mathcal L(X)$ satisfy
\begin{enumerate}
\item[\textup{(H1)}] \emph{pointwise convergence:} $T_hf\to Tf$ in $X$ for
every $f\in X$;
\item[\textup{(H2)}] \emph{collective compactness:} for every sequence
$h_j\to0$ and every bounded sequence $(f_j)\subset X$, the sequence
$(T_{h_j}f_j)$ contains a strongly convergent subsequence.
\end{enumerate}
The spectrum of $T_h$ converges to that of $T$ in the following sense. Let $\Omega\subset\mathbb C$ be open and bounded with rectifiable Jordan
boundary $\Gamma$ such that $\Gamma\cap\sigma(T)=\emptyset$ and
$0\notin\overline\Omega$. Then there exists $h_\Gamma>0$ such that for all
$h<h_\Gamma$:
\begin{enumerate}
\item $\Gamma\subset\rho(T_h)$ and
$\displaystyle\sup_{z\in\Gamma}\|(zI-T_h)^{-1}\|\le C_\Gamma<\infty$;
\item the Riesz spectral projections
$E:=\frac{1}{2\pi i}\oint_\Gamma(zI-T)^{-1}\,dz$ and
$E_h:=\frac{1}{2\pi i}\oint_\Gamma(zI-T_h)^{-1}\,dz$
satisfy $\|(E_h-E)|_{\mathcal R(E)}\|\to0$ and
$\operatorname{rank}E_h=\operatorname{rank}E<\infty$;
\item consequently, the total algebraic multiplicity of
$\sigma(T_h)\cap\Omega$ equals that of $\sigma(T)\cap\Omega$, and
$\sigma(T_h)\cap\Omega\to\sigma(T)\cap\Omega$ in the Hausdorff sense.
\end{enumerate}
\end{theorem}

We therefore apply \cref{thm:abstract_cc} to the shifted resolvent operators
\begin{equation}
T_h=(\widetilde A_h^k)^{-1},
\qquad
T=(\widetilde A^k)^{-1},
\qquad k\in\{0,n\}.
\end{equation}
Once \textup{(H1)}--\textup{(H2)} are established for \(T_h\), the theorem
yields spectral convergence of the shifted resolvents. Through the spectral
mapping
\begin{equation}
\lambda\in\sigma(A_h^k)
\quad\Longleftrightarrow\quad
(\lambda+c)^{-1}\in\sigma(T_h),
\end{equation}
this gives the corresponding convergence and preservation of the
eigenvalues of \(A_h^k\) for sufficiently small \(h\).

Firstly, for \emph{$k=0$}. Let $X:=L^2\Lambda^0(M)$ and consider the shifted forms
\begin{equation}\label{eq:form0}
  a_c(u,v):=\eta(du,dv)+(i_{\bm\beta}du,v)+c(u,v),
  \qquad u,v\in H\Lambda^0(M).
\end{equation}
Define $T\in\mathcal L(X)$ by $Tf:=u$, where $u\in H\Lambda^0(M)$ solves
$a_c(u,v)=(f,v)$ for all $v\in H\Lambda^0(M)$, and let $T_hf:=u_h\in V_h^0$ be its
Galerkin approximation. Since
\[
 |a_c(u,u)|
  \ \ge\ \eta\|du\|^2-C_{\bm\beta}\|du\|\|u\|+c\|u\|^2
  \ \ge\ \frac\eta2\|du\|^2
  +\Bigl(c-\frac{C_{\bm\beta}^2}{2\eta}\Bigr)\|u\|^2
  \ \ge\ \alpha\,\|u\|_{H\Lambda^0}^2
\]
with $\alpha:=\min\{\eta/2,\ c-C_{\bm\beta}^2/(2\eta)\}>0$ if $c\ge C_{\bm\beta}^2/\eta$, the form $a_c$ is uniformly coercive on $H\Lambda^0$
and on every subspace $V_h^0$. Then both problems are well posed for every $h$,
with
\begin{equation}\label{eq:stab0}
  \|Tf\|_{H\Lambda^0}+\|T_hf\|_{H\Lambda^0}\ \le\ \frac{2}{\alpha}\,\|f\| .
\end{equation}

\emph{\textup{(H1)}.} By C\'ea's lemma, we have
$\|Tf-T_hf\|_{H\Lambda^0}\le\frac{M_a}{\alpha}
\inf_{v_h\in V_h^0}\|Tf-v_h\|_{H\Lambda^0}$ with the continuity constant
$M_a$ of $a_c$. The right-hand side tends to $0$ for each fixed $f$ by the approximation property of
the family $\{V_h^0\}$ \cite{arnold2006finite}.

\emph{\textup{(H2)} and compactness of $T$.} By \cref{eq:stab0}, for any
$h_j\to0$ and $\|f_j\|\le1$ the sequence $T_{h_j}f_j$ is bounded in
$H^1(M)$. The compact embedding $H^1(M)\hookrightarrow\hookrightarrow
L^2(M)$ extracts a strongly convergent subsequence. The compact embedding also shows that $T$ is compact.

Then for \emph{$k=n$}. Let $X:=L^2\Lambda^n(M)$. The shifted problem in mixed form reads: given
$f\in X$, find $(\sigma,u)\in H\Lambda^{n-1}(M)\times L^2\Lambda^n(M)$ with
\begin{equation}\label{eq:mixed_primal}
\begin{aligned}
  (\sigma,\tau)-(i_{\bm\beta}u,\tau)-\eta(u,d\tau)&=0
  &&\forall\,\tau\in H\Lambda^{n-1}(M),\\
  (d\sigma,v)+c(u,v)&=(f,v)
  &&\forall\,v\in L^2\Lambda^n(M),
\end{aligned}
\end{equation}
and set $Tf:=u$; the discrete solution operator $T_hf:=u_h$ is defined by
the same equations on $V_h^{n-1}\times V_h^n$, with $i_{\bm\beta}$ replaced
by $i_{\bm\beta,h}=Q_h^{n-1}i_{\bm\beta}$ (the projection is invisible in
\cref{eq:mixed_primal} tested against $\tau_h\in V_h^{n-1}$). Writing
$b\bigl((\sigma,u),(\tau,v)\bigr)$ for the sum of the two left-hand sides,
the choice of test functions
\begin{equation}\label{eq:test_pair}
  (\tau,v):=\Bigl(\sigma,\ \eta u+\gamma\,d\sigma\Bigr),
  \qquad \gamma:=\frac{\eta}{2c},
\end{equation}
gives, after absorbing the advective and coupling terms by Young's
inequality and selecting $c\ge C_{\bm\beta}^2/\eta$,
\begin{equation}\label{eq:stab_kn}
  |b\bigl((\sigma,u),(\tau,v)\bigr)|
  \ \ge\ \frac12\|\sigma\|^2+\frac{\gamma}{2}\|d\sigma\|^2
  +\frac{\eta c}{4}\|u\|^2,
  \qquad
  \|(\tau,v)\|_{H\Lambda\times L^2}\le C\,\|(\sigma,u)\|_{H\Lambda\times L^2},
\end{equation}
i.e.\ a uniform inf-sup condition. Hence \cref{eq:stab_kn} holds on
$V_h^{n-1}\times V_h^n$, with constants independent of $h$. In particular
both problems are well posed for every $h$, with
\begin{equation}\label{eq:stab_kn2}
  \|\sigma\|_{H\Lambda}+\|u\|+\|\sigma_h\|_{H\Lambda}+\|u_h\|
  \ \le\ C_s\,\|f\| .
\end{equation}

\emph{\textup{(H1)}.} By the standard argument of mixed finite elements, by
\cref{eq:stab_kn}, we have
\[
  \|\sigma-\sigma_h\|_{H\Lambda}+\|u-u_h\|
  \ \le\ C\inf_{(\tau_h,v_h)\in V_h^{n-1}\times V_h^n}
  \bigl(\|\sigma-\tau_h\|_{H\Lambda}+\|u-v_h\|\bigr),
\]
with $C$ independent of $h$, and the right-hand side tends to $0$ for each
fixed $f$ by density of smooth forms and the approximation property of the
finite element complex \cite{arnold2006finite}.

\emph{$H^*\Lambda^n$ regularity.} The first equation in 
\cref{eq:mixed_primal} gives $(u,d\tau)=(g,\tau)$ for $g=(\sigma-i_{\bm\beta} u)/\eta \in L^2\Lambda^{n-1}$ and $\forall \tau\in H\Lambda^{n-1}$, thus $u\in H^*\Lambda^n$ and $\delta u = g$. This also implies that 
\begin{equation}\label{eq:primal_H1}
  \|u\|_{H^*\Lambda^n}\lesssim \|u\|+\|\sigma\| \lesssim\|f\| .
\end{equation}
Together with the compact embedding $H\Lambda^n\cap H^*\Lambda^n\hookrightarrow \hookrightarrow L^2\Lambda^n$, this shows in particular that
$T:X\to X$ is compact.

\emph{\textup{(H2)}.} We now prove collective compactness of the discrete solution operators.
Let \(h_j\to0\), let \(\{f_j\}\subset X\) be bounded, and let
\[
  (\sigma_j,u_j)
  :=
  (\sigma_{h_j},T_{h_j}f_j)
  \in
  V_{h_j}^{n-1}\times V_{h_j}^{n}
\]
denote the corresponding discrete solutions. By
\eqref{eq:stab_kn2},
\begin{equation}\label{eq:discrete_mixed_uniform_bound}
  \sup_j
  \left(
    \|\sigma_j\|_{H\Lambda}
    +
    \|u_j\|
  \right)
  <\infty.
\end{equation}
By the definition of the discrete coderivative, the first discrete
equation is equivalent to
\[
  \sigma_j
  =
  i_{\bm\beta,h_j}^nu_j
  +
  \eta\delta_{h_j}^nu_j.
\]
It follows from the uniform boundedness of
\(i_{\bm\beta,h_j}^n\) that
\begin{equation}\label{eq:discrete_top_graph_bound}
\begin{aligned}
  \eta\|\delta_{h_j}^nu_j\|
  &\le
  \|\sigma_j\|
  +
  \|i_{\bm\beta,h_j}^nu_j\|
\\
  &\le
  \|\sigma_j\|
  +
  C_{\bm\beta}\|u_j\|,
\end{aligned}
\end{equation}
and therefore
\begin{equation}\label{eq:discrete_top_uniform_graph_bound}
  \sup_j
  \left(
    \|u_j\|
    +
    \|\delta_{h_j}^nu_j\|
  \right)
  <\infty.
\end{equation}
For the standard FEEC families considered here on shape-regular mesh
families, the uniformly \(L^2\)-bounded commuting projections and the
compact embedding
\[
  H^*\Lambda^n(M)
  \simeq
  H^1\Lambda^0(M)
  \hookrightarrow\hookrightarrow
  L^2\Lambda^0(M)
  \simeq
  L^2\Lambda^n(M)
\]
yield the top-degree discrete Rellich compactness property in the
\(\delta_h^n\)-graph norm; see
\cite{christiansen2007stability,arnold2010finite}.
Since discrete compactness results are often formulated in
terms of \(d\), we briefly verify the corresponding
\(\delta_h^n\)-version needed here.

Let $\pi_h^\ell:L^2\Lambda^\ell(M)\longrightarrow V_h^\ell$ for $\ell=n-1,n$ be uniformly \(L^2\)-bounded commuting projections satisfying $ \pi_h^\ell|_{V_h^\ell}=I$ and
  $d^{n-1}\pi_h^{n-1}=
  \pi_h^nd^{n-1}$. We define
their \(L^2\)-adjoints
\begin{equation}\label{eq:adjoint_reconstruction}
  R_h^\ell:V_h^\ell\longrightarrow L^2\Lambda^\ell(M),
  \qquad
  R_h^\ell:=(\pi_h^\ell)^*,
\end{equation}
by
\begin{equation}\label{eq:adjoint_reconstruction_characterization}
  (w,R_h^\ell v_h)
  =
  (\pi_h^\ell w,v_h),
  \qquad
  \forall w\in L^2\Lambda^\ell(M),\quad
  \forall v_h\in V_h^\ell.
\end{equation}
By the Riesz representation theorem, \(R_h^\ell\) is well defined, and
the uniform \(L^2\) boundedness of \(\pi_h^\ell\) gives
\begin{equation}\label{eq:adjoint_reconstruction_bound}
  \|R_h^\ell v_h\|
  \le
  \|\pi_h^\ell\|_{\mathcal L(L^2,L^2)}\|v_h\|
  \le
  C_\pi\|v_h\|.
\end{equation}
Moreover, for every
\(\tau\in H\Lambda^{n-1}\) and \(v_h\in V_h^n\),
\begin{equation}\label{eq:adjoint_reconstruction_commutation}
\begin{aligned}
  (d^{n-1}\tau,R_h^nv_h)
  &=
  (\pi_h^nd^{n-1}\tau,v_h)
   = 
  (d^{n-1}\pi_h^{n-1}\tau,v_h) \\
  &=
  (\pi_h^{n-1}\tau,\delta_h^nv_h)
   =
  (\tau,R_h^{n-1}\delta_h^nv_h).
\end{aligned}
\end{equation}
Since \(\delta^n=(d^{n-1})^*\) and the right-hand side is bounded with respect to the \(L^2\)-norm of
	\(\tau\):
	$
	|(\tau,R_h^{n-1}\delta_h^nu_h)|
	\le
	C_\pi\|\delta_h^nu_h\|\,\|\tau\|
	$, this proves
\begin{equation}\label{eq:adjoint_reconstruction_delta}
  R_h^nv_h\in H^*\Lambda^n(M),
  \qquad
  \delta^nR_h^nv_h
  =
  R_h^{n-1}\delta_h^nv_h,
\end{equation}
and hence
\begin{equation}
\begin{aligned}
  \|R_h^nv_h\|_{H^*\Lambda^n}^2
  &=
		\|R_h^nv_h\|^2
		+
		\|\delta^nR_h^nv_h\|^2
		\notag\\
		&=
		\|R_h^nv_h\|^2
		+
		\|R_h^{n-1}\delta_h^nv_h\|^2
		\notag\\
		&\le
		C_\pi^2
		\left(
		\|v_h\|^2
		+
		\|\delta_h^nv_h\|^2
		\right).
        \end{aligned}
\end{equation}
Taking  \(v_h=u_j\) and using
\eqref{eq:discrete_top_uniform_graph_bound}, we have that
\(\{R_{h_j}^nu_j\}\) is bounded in \(H^*\Lambda^n(M)\). The compact
embedding above therefore gives, after passing to a subsequence,
\begin{equation}\label{eq:reconstructed_strong_convergence}
  R_{h_j}^nu_j
  \longrightarrow u
  \qquad
  \text{strongly in }L^2\Lambda^n(M) \text{ as } h_j\rightarrow0
\end{equation}
for some \(u\in L^2\Lambda^n(M)\).

Let \(Q_h^n\) be the \(L^2\)-orthogonal projection onto \(V_h^n\).
Because \(\pi_h^n\) is a projection onto \(V_h^n\), one has
\begin{equation}\label{eq:QhRh_identity}
  Q_h^nR_h^n=I
  \qquad\text{on }V_h^n.
\end{equation}
Indeed, for \(v_h,w_h\in V_h^n\),
\[
  (w_h,Q_h^nR_h^nv_h)
  =
  (\pi_h^nw_h,v_h)
  =
  (w_h,v_h).
\]
Furthermore, the standard \(L^2\)-approximation property gives
\(Q_h^nw\to w\) strongly for every \(w\in L^2\Lambda^n(M)\) as $h\rightarrow0$.
Consequently, using \(u_j=Q_{h_j}^nR_{h_j}^nu_j\),
\begin{equation}\label{eq:transfer_compactness_to_discrete}
\begin{aligned}
  \|u_j-u\|
  &\le
  \left\|
    Q_{h_j}^n
    \bigl(R_{h_j}^nu_j-u\bigr)
  \right\|
  +
  \|(Q_{h_j}^n-I)u\| \\
  &\le
  \|R_{h_j}^nu_j-u\|
  +
  \|(Q_{h_j}^n-I)u\|
  \longrightarrow0,  \quad\text{ as } h_j\rightarrow 0.
\end{aligned}
\end{equation}
Thus \(\{u_j\}=\{T_{h_j}f_j\}\) contains a strongly convergent
subsequence in \(X\), which proves \textup{(H2)}. Thus \(\{T_h\}\) is collectively compact.

Therefore, we obtain the spectral convergence for any eigenvalue of $T$ as $h\rightarrow 0$. However, it does not exclude the possibility that, for a fixed \(h\), discrete eigenvalues far from the origin lie in the left half-plane and hence violate Assumption~\ref{asm:spectrum}, even though such eigenvalues disappear in the limit \(h\to0\). To rule out this possibility for sufficiently small \(h\), the next lemma shows that all discrete eigenvalues that could violate Assumption~\ref{asm:spectrum} lie in a fixed compact subset of \(\mathbb C\), uniformly in \(h\). Once such a compact confinement is available, one may choose a fixed bounded region containing only the zero eigenvalue of the continuous operator. The preservation of algebraic multiplicity in \cref{thm:abstract_cc} then excludes any additional discrete eigenvalues in this region for all sufficiently small \(h\), thereby proving Assumption~\ref{asm:spectrum}. To provide a unified proof for this lemma, we set $D_hu_h:=d^0u_h$ for $k=0$ and $D_hu_h:=\delta_hu_h$ for $k=n$, the
discrete codifferential is defined by
$(\delta_hu_h,\tau_h)=(u_h,d^{n-1}\tau_h)$ for all
$\tau_h\in V_h^{n-1}$.

\begin{lemma}[distribution of eigenvalues]\label{lem:localization}
Define  $R_d:=\frac{C_{\bm\beta}+\sqrt{C_{\bm\beta}^2+2\eta\mu}}{2\eta}$. Let $k\in\{0,n\}$ and let
$(\widetilde\lambda_h,u_h)$ be an eigenpair of $\widetilde A_h^k$ with
$\operatorname{Re}\widetilde\lambda_h\le\frac{\mu}{2}+c$ for fixed $c\ge\eta R_d^2+C_{\bm\beta}R_d+1$. Then
\begin{equation}\label{eq:localization}
  \|D_hu_h\|\le R_d\,\|u_h\|,
  \qquad
  |\widetilde\lambda_h-c|\le \eta R_d^2+C_{\bm\beta}R_d .
\end{equation}
Consequently
\begin{equation}\label{eq:K_def}
  \sigma(\widetilde A_h^k)\subset
  K\cup\Bigl\{\operatorname{Re}z>\frac{\mu}{2}+c\Bigr\},
  \qquad
  K:=\Bigl\{z\in\mathbb C:\ |z-c|\le \eta R_d^2+C_{\bm\beta}R_d,\
  \operatorname{Re}z\le\frac{\mu}{2}+c\Bigr\},
\end{equation}
and $0\notin K$.
\end{lemma}

\begin{proof}
For $k=0$,
\begin{equation}\label{eq:k0id}
    \tilde{\lambda}_h\|u_h\|^2=(\widetilde A_h^0u_h,u_h)
=\eta\|du_h\|^2+(i_{\bm\beta}du_h,u_h)+c\|u_h\|^2.
\end{equation}
For
$k=n$, the mixed formulation gives
$\sigma_h=Q_h^{n-1}(i_{\bm\beta}u_h)+\eta\,\delta_hu_h$ and
$\widetilde A_h^nu_h=d\sigma_h+cu_h$, whence, using the definition of
$\delta_h$ with $\tau_h=\sigma_h$, we have
\begin{equation}\label{eq:knid}
\tilde{\lambda}_h\|u_h\|^2=  (\widetilde A_h^nu_h,u_h)
  =(\sigma_h,\delta_hu_h)+c\|u_h\|^2
  =\eta\|\delta_hu_h\|^2+(i_{\bm\beta}u_h,\delta_hu_h)+c\|u_h\|^2 .
\end{equation}
In both cases, taking the real part gives
\[
  \operatorname{Re}\widetilde\lambda_h\,\|u_h\|^2
  \ \ge\ \eta\|D_hu_h\|^2-C_{\bm\beta}\|D_hu_h\|\|u_h\|+c\|u_h\|^2 .
\]
If $\operatorname{Re}\widetilde\lambda_h\le\frac\mu2+c$, the quadratic
inequality
$\eta t^2-C_{\bm\beta}t-\frac\mu2\le0$ for $t:=\|D_hu_h\|/\|u_h\|$ yields
$t\le R_d$. Then both \eqref{eq:k0id} and \eqref{eq:knid} imply
\[
  |\widetilde\lambda_h-c|
  \le\frac{\eta\|D_hu_h\|^2+C_{\bm\beta}\|D_hu_h\|\|u_h\|}{\|u_h\|^2}
  \le\eta R_d^2+C_{\bm\beta}R_d .
\]
Finally, if $0\in K$, then $c=|0-c|\le \eta R_d^2+C_\beta R_d<c$ by the choice
of $c$, which is a contradiction.
\end{proof}

Finally, we present the proof of \cref{lm:spectrum_small_h}.

\begin{proof}[Proof of \cref{lm:spectrum_small_h}]
Take $c\ \ge\ c_*:=\max\Bigl\{\eta R_d^2+C_{\bm\beta}R_d+1,\ \frac{C_{\bm\beta}^2}{\eta}\Bigr\}$. Fix $k\in\{0,n\}$ and let $T$, $T_h$ be as above, so that
$T=(\widetilde A^k)^{-1}$ and $T_h=(\widetilde A_h^k)^{-1}Q_h^k$ on
$X=L^2\Lambda^k$. The inverses exist for every $h$ by
\cref{eq:stab0} and \cref{eq:stab_kn}. The nonzero eigenvalues of the two operators are reciprocal to each other and have the same multiplicities.

By \cref{thm:spec},
$\sigma(\widetilde A^k)\subset\{c\}\cup
\{\operatorname{Re}z\ge\mu+c\}$ with $c$ algebraically simple. The compact
set $K$ of \cref{eq:K_def} therefore satisfies
$\sigma(\widetilde A^k)\cap K=\{c\}$ and $\partial K\cap\sigma(\widetilde A^k)=\emptyset$.
Moreover, $c$ lies in the interior of $K$ and $0\notin K$. The inversion
$\iota(z):=1/z$ is a homeomorphism of a neighborhood of $K$. Set
$\Omega:=\iota(\operatorname{int}K)$ and $\Gamma:=\iota(\partial K)$, then
$\Omega$ is open and bounded with rectifiable Jordan boundary $\Gamma$,
$0\notin\overline\Omega$ (as $\infty\notin K$), and
\[
  \sigma(T)\cap\Omega=\{1/c\}
  \quad\text{with algebraic multiplicity }1,
  \qquad
  \Gamma\cap\sigma(T)=\emptyset .
\]

The hypotheses of
\cref{thm:abstract_cc} were verified above for both degrees. Hence there is
$h_\Gamma>0$ such that for $h<h_\Gamma$ the total algebraic multiplicity of
$\sigma(T_h)\cap\Omega$ equals $1$.

As shown in the beginning, the value $c$ is an eigenvalue of
$\widetilde A_h^k$ for \emph{every} $h$. Since $1/c\in\Omega$, it is the only eigenvalue in $\Omega$ and has multiplicity of 1, then we have
\[
  \sigma(\widetilde A_h^k)\cap K=\{c\}\ \text{algebraically simple}.
\]

Let $h_0:=h_\Gamma$ (minimized over the two
degrees). For $h<h_0$ and any eigenvalue
$\widetilde\lambda_h\ne c$ of $\widetilde A_h^k$, we have
$\widetilde\lambda_h\notin K$. Then \cref{lem:localization} forces
$\operatorname{Re}\widetilde\lambda_h>\frac\mu2+c$. Finally undoing the shift, we have that $0$
is an algebraically simple eigenvalue of $A_h^k$, and every other
eigenvalue satisfies $\operatorname{Re}\lambda_h>\mu/2>0$. This is
Assumption~\ref{asm:spectrum}, with the quantitative gap $\mu/2$.
\end{proof}

\subsection{Choice 2: $M$-matrix discretization for $\mathcal P_{1}^{-}\Lambda^k$ for arbitrary $h$}

An alternative way to guarantee Assumption \ref{asm:spectrum} for arbitrary mesh size is to construct, for the degrees $k\in\{0,n\}$, discretizations based on the lowest-order Whitney forms whose system matrices are singular irreducible $M$-matrices: by the Perron--Frobenius theorem, the zero eigenvalue of such a matrix is simple, its left and right eigenvectors are strictly positive, and all remaining eigenvalues have strictly positive real parts, which is exactly the required spectral distribution. The construction is restricted to the Euclidean setting because it relies on Delaunay or angle conditions. Throughout this subsection, \(M\) is a polyhedral domain with periodic boundary conditions (for instance a flat torus), so that every facet is interior and is endowed with the Euclidean metric. \(\mathcal T_h\) is an affine simplicial triangulation.

In this singular case the discrete maximum principle takes the form $A_h\omega_h\le0\implies A_h\omega_h=0$. If essential boundary conditions render the matrix nonsingular, it degenerates to the usual monotonicity $A_h\omega_h\le0\implies\omega_h\le0$. Throughout this subsection, $\bm\beta\in W^{1,\infty}(M)$. The construction combines the edge-averaged finite element (EAFE) method \cite{xu1999monotone} and its simplex extension (SAFE) \cite{wu2020simplex} with an exponential fitting of the flux.

\subsubsection{Exponential-fitting discretization}
Grouping the diffusive and convective terms, the mixed scheme \eqref{eq:add} can be written in flux form,
\begin{equation*}
  \begin{aligned}
    (\sigma_h,\tau_h) - \bigl(u_h,(\eta d+j_{\bm\beta})\tau_h\bigr) &= 0, \quad && \forall\tau_h\in V_h^{k-1},\\
   (\partial_t u_h,v_h) + (d\sigma_h,v_h) + \bigl(du_h,(\eta d+j_{\bm\beta})v_h\bigr) &= 0, \quad && \forall v_h\in V_h^k,
  \end{aligned}
\end{equation*}
where $J_{\eta,\bm\beta}:=\eta d+j_{\bm\beta}$ acts as a $(k+1)$-form flux operator. Our goal is to construct a discrete flux operator $J_h^k:V_h^k\to V_h^{k+1}$ leading to an $M$-matrix discretization. The discrete contraction is then recovered as the formal adjoint of the flux remainder,
\begin{equation}\label{eq:contraction_adjoint0}
  (i_{\bm\beta,h}^{k+1} u_h, v_h) = \bigl(u_h, (J_h^k - \eta d^k)v_h\bigr), \qquad u_h\in V_h^{k+1},\ v_h\in V_h^k,
\end{equation}
so that both contributions of diffusion and advection are encoded in the single operator $J_h^k$.

When $\bm\beta$ is a gradient field, the Witten deformation of \Cref{sec:weight} motivates an exponential fitting of the flux: if $\bm\beta=(d\psi)^\sharp$, then
\begin{equation}
    \eta e^{-\psi/\eta}d(e^{\psi/\eta}u)
    =\eta du+d\psi\wedge u
    =\eta du+j_{\bm\beta}u,
\end{equation}
so the flux is a conjugated exterior derivative. In general, $\bm\beta$ is not a gradient field. For $k=0$, we use an edgewise potential as in \cite{xu1999monotone}; for $k>0$, we use the elementwise construction of \cite{wu2020simplex}.

Denote by $\Sigma_T^k$ and $\Sigma_h^k$ the sets of oriented $k$-dimensional subsimplices of $T\in \mathcal T_h$ and $\mathcal T_h$, respectively, and let $\ell_S$ be the canonical degree of freedom associated with $S$. We write $\varphi_S$ for the corresponding Whitney basis form.

For $k=0$, we use the edgewise construction of EAFE
\cite{xu1999monotone}. For each oriented edge
$E=\overrightarrow{\bm a_i\bm a_j}$, define $\psi_E$,
up to a constant, by
\begin{equation}\label{eq:edge_potential}
  \frac{\partial\psi_E}{\partial\tau_E}
  =\bm\beta\cdot\tau_E
  \qquad\text{on }E,
\end{equation}
where $\tau_E$ is the unit tangent vector of $E$. Set
\[
  \mathcal E_E
  :=\Bigl(\fint_E e^{\psi_E/\eta}\Bigr)^{-1},
  \qquad
  \delta_E(w):=w(\bm a_j)-w(\bm a_i).
\]
The fitted flux $\tilde J_h^0v_h\in V_h^1$ is defined by
\begin{equation}\label{eq:edge_fitted_flux}
  \ell_E(\tilde J_h^0v_h)
  :=\eta\mathcal E_E\,\delta_E(e^{\psi_E/\eta}v_h),
  \qquad E\in\Sigma_h^1.
\end{equation}
These degrees of freedom are independent of the additive
constant in $\psi_E$ and are single-valued on shared edges.

For $k>0$, an analogous construction of the fitted flux on
$S\in\Sigma_h^{k+1}$ would require solving a PDE for a local
potential, whose compatibility conditions need not hold for a
general vector field $\bm\beta$. We use the elemental construction of SAFE
\cite{wu2020simplex}. On each element $T$, set
\begin{equation}\label{eq:betaT}
  \bm\beta_T:=\fint_T\bm\beta\,dx,
  \qquad
  \psi_T(\bm x):=\bm\beta_T\cdot\bm x.
\end{equation}
Let $\Pi_T^k$ be the canonical Whitney interpolation operator
and define
\begin{equation}\label{eq:WE}
  \mathcal W_T^k v:=\Pi_T^k(e^{\psi_T/\eta}v),
  \qquad
  \mathcal E_T^k:=(\mathcal W_T^k)^{-1}.
\end{equation}
The inverse is taken on $\mathcal P_1^-\Lambda^k(T)$.
Both operators are diagonal in the Whitney basis; explicitly,
\[
  \ell_S(\mathcal E_T^k w_h)
  =\mathcal E_S^T\,\ell_S(w_h),
  \qquad
  \mathcal E_S^T
  :=\Bigl(\fint_S e^{\psi_T/\eta}\Bigr)^{-1},
  \qquad S\in\Sigma_T^k.
\]
The elemental fitted flux is
\begin{equation}\label{eq:Jtilde}
  \tilde J_T^k
  :=\eta\mathcal E_T^{k+1}d^k\mathcal W_T^k,
  \qquad 1\le k\le n-1.
\end{equation}
We define $\tilde J_h^k$ elementwise by
\[
  (\tilde J_h^k v_h)|_T
  :=\tilde J_T^k(v_h|_T).
\]
Since $\psi_T$ depends on the element, the flux degrees of
freedom on shared subsimplices need not agree.
For $k=0$, $\tilde J_T^0$ denotes the restriction of the
edgewise construction \eqref{eq:edge_fitted_flux} to $T$.

The following result quantifies the exponential fitting by showing that the degrees of freedom of the fitted flux are obtained from the exponentially weighted continuous flux.

\begin{lemma}[DOF representation]\label{lem:dof_flux}
For a smooth $0$-form $v$ and $E\in\Sigma_h^1$,
\begin{equation}\label{eq:dof_flux_edge}
  \ell_E(\tilde J_h^0v)
  =\mathcal E_E\int_E
    e^{\psi_E/\eta}\operatorname{tr}_E
    \bigl(\eta\,dv+\bm\beta^\flat v\bigr).
\end{equation}
For a smooth $k$-form $v$ on $T$, $1\le k\le n-1$, and
$S\in\Sigma_T^{k+1}$,
\begin{equation}\label{eq:dof_flux}
  \ell_S(\tilde J_T^kv)
  =\mathcal E_S^T
    \int_S\operatorname{tr}_S
    \Bigl[e^{\psi_T/\eta}
      \bigl(\eta\,d^kv+\bm\beta_T^\flat\wedge v\bigr)\Bigr].
\end{equation}
For the elemental SAFE fluxes,
$\tilde J_T^{k+1}\tilde J_T^k=0$ for $1\le k\le n-2$.
\end{lemma}

\begin{proof}
By \eqref{eq:edge_potential},
\[
  \eta\,\frac{\partial}{\partial\tau_E}
    (e^{\psi_E/\eta}v)
  =e^{\psi_E/\eta}
    \left(\eta\frac{\partial v}{\partial\tau_E}
          +(\bm\beta\cdot\tau_E)v\right).
\]
Integrating along $E$ and using \eqref{eq:edge_fitted_flux}
gives \eqref{eq:dof_flux_edge}.

For $k>0$, the definition of $\Pi_T^k$ and Stokes' theorem give
$\ell_S(d\Pi_T^kw)=\int_{\partial S}\operatorname{tr}\Pi_T^kw
=\int_{\partial S}\operatorname{tr}w=\ell_S(dw)$ for smooth $w$.
Taking $w=e^{\psi_T/\eta}v$ and using
$\eta\,d(e^{\psi_T/\eta}v)
=e^{\psi_T/\eta}(\eta\,dv+d\psi_T\wedge v)$
with $d\psi_T=\bm\beta_T^\flat$ yields \eqref{eq:dof_flux}.
Finally, for $1\le k\le n-2$,
\[
  \tilde J_T^{k+1}\tilde J_T^k
  =\eta^2\mathcal E_T^{k+2}d
    (\mathcal W_T^{k+1}\mathcal E_T^{k+1})
    d\mathcal W_T^k
  =\eta^2\mathcal E_T^{k+2}dd\mathcal W_T^k=0,
\]
since $d\mathcal W_T^kv\in\mathcal P_1^-\Lambda^{k+1}(T)$ and
$\mathcal W_T^{k+1}\mathcal E_T^{k+1}=\mathrm{Id}$.
\end{proof}

For \(k=0\), we exploit the graph-Laplacian representation of the
\(\mathcal P_1\) stiffness matrix. The direct \(L^2\) pairing of
\(d^0u_h\) with the fitted flux \(\tilde J_h^0v_h\) does not
 give the desired $M$-matrix structure. We therefore insert a local projection onto constant \(1\)-forms. More
precisely, for \(u_h,v_h\in\mathcal P_1\Lambda^0(T)\),
\begin{equation}
  (d^{0}u_h,d^{0}v_h)_T
  =
  \sum_{E\in \Sigma_T^{1}}
  \omega_{E}^T\,
  \delta_E(u_h)\,\delta_E(v_h),
  \qquad
  \delta_E(u_h)=\ell_E(d^{0}u_h)=u_h(\bm a_j)-u_h(\bm a_i),
\end{equation}
where $\omega_E^T:=-(d^0\varphi_i,d^0\varphi_j)_T$ for the oriented edge $E=\overrightarrow{\bm a_i \bm a_j}$ and $\Sigma_T^k$ denotes the $k$-dimensional subsimplices of $T$. Accordingly, define $\bar{\Pi}_T^1:\mathcal P_1^-\Lambda^1(T)\to \mathcal P_0\Lambda^1(T)$ as the unique element satisfying
\begin{equation}
  (\bar{\Pi}_T^1u_h,q)_T
  =
  \sum_{E\in \Sigma_T^{1}}
  \omega_E^T\,
  \ell_E(u_h)\,\ell_E(q),
  \qquad
  \forall q\in \mathcal P_0\Lambda^1(T),
\end{equation}
and $\bar\Pi_h^1|_T:=\bar\Pi_T^1$. Then
\begin{equation}\label{eq:quad0}
(d^0 u_h,\bar{\Pi}_h^1\tilde{J}^0_h v_h) = \sum_{T\in\mathcal T_h}\sum_{E\in\Sigma_T^1} \omega_E^T\,\delta_E(u_h)\,\ell_E\bigl(\tilde{J}^0_h v_h\bigr),
\end{equation}
where, for the oriented edge $E = \overrightarrow{\bm a_i \bm a_j}$, the flux degree of freedom expands as
\begin{equation}\label{eq:exp0}
  \ell_E(\tilde J_h^0v_h)
  =\eta\mathcal E_E\,\delta_E(e^{\psi_E/\eta}v_h)
  =\eta\mathcal E_E
    \Bigl(e^{\psi_E(\bm a_j)/\eta}v_h(\bm a_j)
          -e^{\psi_E(\bm a_i)/\eta}v_h(\bm a_i)\Bigr).
\end{equation}
Furthermore, to transfer the \(M\)-matrix structure of the spatial
discretization to the semidiscrete evolution operator, we replace the
 mass matrix associated with the temporal term by a positive diagonal one. This avoids
the mixing of degrees of freedom caused by the inverse of a  mass
matrix and ensures that multiplication by the inverse lumped mass matrix
preserves the $M$-matrix structure. We therefore use the nodal mass-lumped pairing: 
\begin{equation}
    (\partial_t u_h,v_h)\approx (\partial_t u_h,v_h)_h := \sum_{z\in \Sigma_h^0} m_z\,\partial_t u_h(z)\,v_h(z),
    \qquad
    m_z=\int_M \varphi_z>0,
\end{equation}
and denote by $\tilde{M}_h^0 : V_h^{0} \to V_h^{0}$ the induced self-adjoint operator, 
$$(u_h,v_h)_h = (\tilde{M}_h^0u_h,v_h)=(u_h,\tilde{M}_h^0v_h).$$ 
The lumping contributes only a positive diagonal and will allow us to identify the eigenvalues of $A_h^0$ below.

For the volume-form case \(k=n\), the auxiliary variable
\(\sigma_h\in V_h^{n-1}\) is eliminated through the inverse of its facet
mass matrix. A consistent mass matrix would introduce nonlocal couplings. We therefore replace it by a positive diagonal facet mass
matrix, so that the elimination remains local and preserves the
off-diagonal sign pattern required for the \(M\)-matrix property. Following
\cite{sacco1997stabilization,baranger1994application,bossavit1999yee}, we
define
\[
(\sigma_h,\tau_h)\approx (\sigma_h,\tau_h)_h :=\sum_{F\in\Sigma_h^{n-1}}\mathbb M_F\,\ell_F(\sigma_h)\,\ell_F(\tau_h),
\qquad
\mathbb M_F:=\frac{1}{|F|}\sum_{T \supset F} d_{T,F},
\]
where $d_{T,F}$ is the {signed} distance from the circumcenter of $T$ to the hyperplane of $F$ (positive when the circumcenter lies on the same side as $T$). This quadrature diagonalizes the mass matrix with an $\mathcal O(h)$ consistency error. Its positivity is ensured under a strict Delaunay triangulation \cite{sacco1997stabilization,baranger1994application,bossavit1999yee}. Indeed, $\mathbb M_F$ is the diagonal entry of the circumcentric discrete Hodge star, and the resulting bilinear form $(\cdot,\cdot)_h$ is the associated Discrete Exterior Calculus (DEC) inner product on primal \((n-1)\)-cochains \cite{hirani2003discrete,hirani2013delaunay}. Denote by $\tilde{M}_h^{n-1}: V_h^{n-1} \to V_h^{n-1}$ the induced self-adjoint operator, $(\sigma_h,\tau_h)_h=(\tilde M_h^{n-1}\sigma_h,\tau_h)=(\sigma_h,\tilde M_h^{n-1}\tau_h)$. Introducing the change of variables \(\zeta_h=\tilde{M}^{n-1}_h\tau_h\), the first equation in \eqref{eq:sadd} can be rewritten as
\begin{equation}
  (\sigma_h,\tau_h)_h - (u_h,\tilde{J}_h^{n-1}\tau_h) = 0,
  \qquad\Leftrightarrow\qquad
  (\sigma_h,\zeta_h) - \bigl(u_h,\tilde{J}_h^{n-1}(\tilde{M}^{n-1}_h)^{-1}\zeta_h\bigr) = 0.
\end{equation}
No lumping is needed for the temporal term, since the mass matrix of $\mathcal P^-_1\Lambda^n$ is already diagonal and positive. Note that the flux operators for 0-forms and $(n-1)$-forms require additional modification, whereas for the other forms an $L^2$-projection is used to obtain a conforming flux. The unified flux operator $J_h^k : V_h^k \to V_h^{k+1}$ is thus
\begin{equation}\label{eq:Jh_unified}
  J_h^k
  :=
  \begin{cases}
    Q_h^1 \bar{\Pi}_h^1\tilde{J}_h^0(\tilde{M}_h^0)^{-1}, & k=0,\\[2mm]
    Q_h^{k+1}\tilde{J}_h^k, & 1 \le k \le n-2,\\[2mm]
    \tilde{J}_h^{n-1}(\tilde{M}_h^{n-1})^{-1}, & k= n-1,
  \end{cases}
\end{equation}
where $Q_h^{k+1}$ is the $L^2$-orthogonal projection onto $V_h^{k+1}$, and the discrete contraction is given by \eqref{eq:contraction_adjoint0}. First-order convergence results for the corresponding EAFE and SAFE model problems are established under the assumptions of \cite{xu1999monotone,wu2020simplex}.

We first establish a bound for the discrete contraction operator which is uniform in \emph{both} $h$ and $\eta$.

\begin{lemma}[uniform elemental bound]\label{lem:elemental}
There exists $C>0$, depending only on shape regularity and $n,k$, such that
\[
  \|(\tilde J_T^k-\eta d^k)v_h\|_{L^2(T)}
  \le C|\bm\beta_T|\,\|v_h\|_{L^2(T)},
  \quad
  v_h\in\mathcal P_1^-\Lambda^k(T),\quad
  1\le k\le n-1,\quad h,\eta>0.
\]
For $k=0$, the same estimate holds with $|\bm\beta_T|$ replaced by
$\|\bm\beta\|_{L^\infty(T)}$.
Consequently, for $1\le k\le n-2$, the operators
\eqref{eq:contraction_adjoint0} satisfy Assumption \ref{asm:stab} with
$C_{\bm\beta}=C\|\bm\beta\|_{L^\infty}$, uniformly in $h$ and $\eta$.
\end{lemma}

\begin{proof}
For $k>0$, the weight
$g_S:=\mathcal E_S^T e^{\psi_T/\eta}|_S$ is positive with unit mean,
$\fint_Sg_S=1$. By \eqref{eq:dof_flux},
\[
  \ell_S((\tilde J_T^k-\eta d^k)v_h)
  =\eta\int_S(g_S-1)\operatorname{tr}_S(d^kv_h)
   +\int_Sg_S\operatorname{tr}_S(\bm\beta_T^\flat\wedge v_h).
\]
For lowest-order elements,
$d^kv_h\in\mathcal P_0\Lambda^{k+1}(T)$, so the first term vanishes.
The second term satisfies
\[
  \left|\int_Sg_S
    \operatorname{tr}_S(\bm\beta_T^\flat\wedge v_h)\right|
  \le C|\bm\beta_T|\,\|v_h\|_{L^\infty(T)}\,|S|.
\]
Expanding $(\tilde J_T^k-\eta d^k)v_h$ in the Whitney basis
$\{\varphi_S\}_{S\in\Sigma_T^{k+1}}$ and combining
$|S|\le Ch_T^{k+1}$,
$\|\varphi_S\|_{L^2(T)}\le Ch_T^{n/2-k-1}$, and
$\|v_h\|_{L^\infty(T)}\le Ch_T^{-n/2}\|v_h\|_{L^2(T)}$
proves the first assertion.

For $k=0$, use
$g_E:=\mathcal E_Ee^{\psi_E/\eta}$ and
\eqref{eq:dof_flux_edge}. The tangential derivative of $v_h$ is
constant on $E$, so the diffusive term again vanishes after subtraction.
Moreover,
\[
  \left|\int_Eg_E\operatorname{tr}_E(\bm\beta^\flat v_h)\right|
  \le \|\bm\beta\|_{L^\infty(T)}
       \|v_h\|_{L^\infty(T)}\,|E|.
\]
The same scaling argument gives the estimate.
The last assertion follows from \eqref{eq:contraction_adjoint0},
the $L^2$ contractivity of $Q_h^{k+1}$, and
$|\bm\beta_T|\le\|\bm\beta\|_{L^\infty(T)}$.
\end{proof} 

\begin{remark}[regularity of $\bm\beta$]
In this section, we assume that $\bm\beta\in W^{1,\infty}$.
For $k>0$, the construction only requires a piecewise constant
approximation of $\bm\beta$. Element averages remain defined for
$\bm\beta\in L^\infty$ and converge in $L^p$, $1\le p<\infty$,
as $h\to0$, although no convergence rate need be available.
For $k=0$, the edgewise ODE additionally requires a well-defined
tangential restriction of $\bm\beta$ to each edge. The assumption $\bm\beta\in W^{1,\infty}$ suffices for both constructions
and is consistent with the convergence analyses in
\cite{xu1999monotone,wu2020simplex}.
\end{remark}

\subsubsection{The $M$-matrix property for spectral preservation}
We now show that the exponential-fitting discretization
\eqref{eq:Jh_unified} satisfies the spectral preservation property
(Assumption~\ref{asm:spectrum}) for arbitrary $h$ and $\eta$. The
strategy is algebraic: under suitable mesh conditions, the system
matrices at the degrees $k\in\{0,n\}$ are singular irreducible
$M$-matrices, so the Perron--Frobenius theorem delivers exactly the
required spectral distribution. The mesh conditions
are as follows.

\begin{assumption}[mesh condition]\label{ass:mesh0}
  The triangulation is conforming and connected. For every edge
\(E=[\bm a_i,\bm a_j]\), the assembled stiffness weight satisfies
\[
  \omega_E
  :=
  -\sum_{T\supset E}
  (d^0\varphi_i,d^0\varphi_j)_T
  \ge0.
\]
For every facet degree of freedom \(F\), the lumped mass
weight satisfies
\[
  \mathbb M_F
  :=
  \frac{1}{|F|}
  \sum_{T\supset F}d_{T,F}
  >0.
\]
\end{assumption}

\begin{remark}[geometric meaning of Assumption \ref{ass:mesh0}]
The relevant conditions are imposed on the assembled weights:
\[
  \omega_E\ge0
  \quad\text{for every edge }E,
  \qquad
  \mathbb M_F>0
  \quad\text{for every facet }F.
\]

\emph{Two dimensions.}
Let \(E=[\bm a_a,\bm a_b]\) be an interior edge shared by
\(T_1,T_2\), and let \(\alpha_i\) be the angle of \(T_i\) opposite
\(E\). Then
\[
  \omega_E^{T_i}
  =
  -(d^0\varphi_a,d^0\varphi_b)_{T_i}
  =
  \frac12\cot\alpha_i,\quad
  \omega_E
  =
  \frac12(\cot\alpha_1+\cot\alpha_2)
  =
  \frac{\sin(\alpha_1+\alpha_2)}
       {2\sin\alpha_1\sin\alpha_2}.
\]
Thus \(\omega_E\ge0\) is equivalent to the Delaunay
condition \(\alpha_1+\alpha_2\le\pi\), whereas
\(\omega_E>0\) is equivalent to strict Delaunay.
Moreover,
\[
  d_{T_i,E}=\frac{|E|}{2}\cot\alpha_i,
  \qquad
  \mathbb M_E
  =
  \frac{d_{T_1,E}+d_{T_2,E}}{|E|}
  =
  \omega_E.
\]
Consequently, on a periodic triangular mesh, strict Delaunay ensures both
requirements. 

\emph{Three dimensions.}
For an edge \(E=[\bm a_a,\bm a_b]\) of a tetrahedron \(T\), let
\(E_T'\) be the opposite edge and let \(\theta_{E_T'}^T\) be the
interior dihedral angle along \(E_T'\). The cotangent formula is
\[
  \omega_E
  =
  \frac16\sum_{T\supset E}
  |E_T'|\cot\theta_{E_T'}^T.
\]
Thus the  sign requirement is an assembled dihedral-angle
condition \cite{xu1999monotone}. Elementwise nonobtuse dihedral
angles are sufficient for \(\omega_E\ge0\), but this is stronger
than the assembled condition. 

For an interior facet \(F=T^+\cap T^-\), the precise positivity
condition is
\[
  \mathbb M_F
  =
  \frac{d_{T^+,F}+d_{T^-,F}}{|F|}
  >0.
\]
A strict Delaunay pair \(T^+,T^-\) satisfies
this condition in any dimension, by positivity of its signed circumcentric dual length
\cite[Theorem~3]{hirani2013delaunay}.
\end{remark}

The spectral conclusion is drawn from the following standard consequence of the Perron--Frobenius theorem \cite{perron1907theorie,frobenius1912matrizen}.

\begin{lemma}[singular irreducible $M$-matrices]\label{lem:pf}
Let $\mathbb A\in\mathbb R^{N\times N}$, $N\ge2$, have nonpositive
off-diagonal entries and zero row sums (or zero column sums). Assume further that \(\mathbb A\) has a symmetric sparsity pattern and is irreducible; equivalently, the graph on the vertex set $\{1,\dots,N\}$ whose edges are the pairs
$\{i,j\}$ with $\mathbb A_{ij}<0$ and $\mathbb A_{ji}<0$ is connected.
Then $0$ is an algebraically simple eigenvalue of $\mathbb A$ and every
other eigenvalue satisfies $\operatorname{Re}\lambda>0$. The hypotheses,
in the zero-row-sum case, are preserved under left multiplication by a
positive diagonal matrix.
\end{lemma}

\begin{proof}
Without loss of generality, assume $\mathbb A\mathbf 1=0$. Set $s:=\max_i\mathbb A_{ii}$ and $\mathbb B:=s I-\mathbb A$, then $\mathbb B\ge0$ entrywise, $\mathbb B$ is irreducible (its off-diagonal pattern is that of $\mathbb A$), and $\mathbb B\mathbf 1=s\mathbf 1$. By the Perron–Frobenius theorem, an irreducible nonnegative matrix with a strictly positive eigenvector has that eigenvalue equal to its spectral radius, and the latter is algebraically simple. Hence $\rho(\mathbb B)=s$ is algebraically simple and $0=s-\rho(\mathbb B)$ is an algebraically simple eigenvalue of $\mathbb A$. For $\xi\in\sigma(\mathbb B)\setminus\{s\}$ we have $|\xi|\le s$ and $\xi\ne s$, which forces $\operatorname{Re}\xi<s$. Thus every nonzero $\lambda=s-\xi\in\sigma(\mathbb A)$ has $\operatorname{Re}\lambda>0$. Meanwhile, left multiplication by a positive diagonal matrix does not alter the required pattern and therefore preserves the same property.
\end{proof}

We next verify that the proposed exponential-fitting discretization satisfies the required spectral preservation assumption.
\begin{lemma}\label{lem:mmatrix}
  Under Assumption \ref{ass:mesh0}, the discretization \eqref{eq:add} with $i_{\bm\beta,h}$ defined by \eqref{eq:contraction_adjoint0} and \eqref{eq:Jh_unified} satisfies Assumption \ref{asm:spectrum} for every $h>0$ and every $\eta>0$.
\end{lemma}

\begin{proof}
{Case $k=0$.}
An eigenpair $(\lambda_h,u_h)$ of $A_h^0$ satisfies
\begin{equation}
 \lambda_h(u_h,v_h)=  (A_h^0 u_h, v_h) = (d^0 u_h, J_h^0 v_h)
 = \bigl(d^0u_h,\bar{\Pi}^1_h \tilde{J}_h^0(\tilde{M}_h^0)^{-1} v_h\bigr),
 \qquad\forall v_h\in V_h^0 .
\end{equation}
Substituting $v_h\mapsto\tilde M_h^0v_h$ and using \eqref{eq:quad0}--\eqref{eq:exp0} yields the generalized eigenvalue problem
\begin{equation}\label{eq:k0expeigen}
 \lambda_h(u_h,v_h)_h = (d^0u_h,\bar{\Pi}^1_h \tilde{J}_h^0 v_h)
  =\eta \sum_{E\in\Sigma_h^1}\omega_E\, \mathcal E_E\, \delta_E(u_h)\,\delta_E(e^{\psi_E/\eta}v_h),
\end{equation}
where the elemental weights were summed to $\omega_E=\sum_{T\supset E}\omega_E^T$ because $\psi_E$ and $\mathcal E_E$ are defined edgewise. Let $\mathsf u\in\mathbb R^{N_0}$ denote the coefficient vector of $u_h$ in
the nodal basis. In matrix form, this is $\mathbb A^0\mathsf u=\lambda_h\,\mathbb M_L\,\mathsf u$ with $\mathbb M_L=\operatorname{diag}(m_z)$. In particular, for the oriented edge $E=\overrightarrow{\bm a_i\bm a_j}$, the off-diagonal entry is
\begin{equation}\label{eq:A0}
  \mathbb A^0_{ij} = \eta\, \omega_E\,\mathcal E_E\, \delta_E(\varphi_j)\,\delta_E(e^{\psi_E/\eta}\varphi_i)=-\eta\,\omega_E\,\mathcal E_E\, e^{\psi_E(\bm a_i)/\eta}\le 0,
\end{equation}
and likewise $\mathbb A^0_{ji}=-\eta\,\omega_E\,\mathcal E_E\, e^{\psi_E(\bm a_j)/\eta}\le 0$, while the entry vanishes for non-adjacent vertices. The off-diagonal entries are nonpositive, and the sparsity pattern is symmetric, since $ \operatorname{sign}(A_{ij}) = \operatorname{sign}(A_{ji}) = \operatorname{sign}(-\omega_E)$ and $\omega_E\ge 0$. Moreover, since \(d^0{1}=0\), the row sums vanish: $ \sum_j \mathbb A^0_{ij} = \bigl(d^0 {1},\,\bar\Pi_h^1\tilde J_h^0\varphi_i\bigr) = 0$.

Let \(G_+\) be the undirected graph whose vertices are the mesh
vertices and whose edges are those with \(\omega_E>0\).
We first show that \(G_+\) is connected. By the definition of $\omega_E$, every \(v_h\in V_h^0\) satisfies
\[
  \|d^0v_h\|_{L^2(M)}^2
  =
  \sum_{E=[\bm a_i,\bm a_j]\in\Sigma_h^1}
  \omega_E
  \bigl|v_h(\bm a_i)-v_h(\bm a_j)\bigr|^2.
\]
If \(G_+\) were disconnected, choose \(v_h\) with nodal value \(1\)
on one connected component and \(0\) on all other components.
Every positive-weight edge then joins vertices with equal nodal
values, while zero-weight edges contribute nothing. Hence
\(d^0v_h=0\). Since \(M\) is connected and \(v_h\) is continuous,
\(v_h\) must be constant, contradicting its construction. Thus
\(G_+\) is connected.

By \eqref{eq:A0}, an edge \(E=[\bm a_i,\bm a_j]\) belongs to
\(G_+\) if and only if $\mathbb A^0_{ij}<0$ and $\mathbb A^0_{ji}<0$, i.e., $\omega_E>0$. Therefore the directed
adjacency graph of \(\mathbb A^0\) is obtained from \(G_+\) by
replacing each edge with two oppositely directed edges. It is
strongly connected, and hence \(\mathbb A^0\) is irreducible.

Left multiplication by the positive diagonal matrix $\mathbb M_L^{-1}$ preserves the sign pattern, the zero row sums, and the irreducibility. Hence $\mathbb M_L^{-1}\mathbb A^0$ is a singular irreducible $M$-matrix, and by \cref{lem:pf}, its zero eigenvalue is algebraically simple while all other eigenvalues have positive real parts. This proves Assumption \ref{asm:spectrum} for $k=0$.

\emph{Case $k=n$.}
The stationary problem reads
\[
\begin{aligned}
(\sigma_h,\tau_h)_h
  -(u_h,\tilde J_h^{n-1}\tau_h)
  &=0,
  &&\forall\tau_h\in V_h^{n-1},\\
(d^{n-1}\sigma_h,v_h)
  &=0,
  &&\forall v_h\in V_h^n.
\end{aligned}
\]
For each cell $T$, let $\varphi_T$ be its Whitney basis
function, normalized by $\int_T\varphi_T=1$, and define
the corresponding auxiliary flux $\sigma_T\in V_h^{n-1}$ by
\begin{equation}\label{eq:auxvar}
  (\sigma_T,\tau_h)_h
  =(\varphi_T,\tilde J_h^{n-1}\tau_h),
  \qquad \forall\tau_h\in V_h^{n-1}.
\end{equation}
Then $A_h^n\varphi_T=d^{n-1}\sigma_T$.

Writing $\mathcal T_h=\{T_i\}$, define
\[
  \mathbb A^n_{ij}
  :=(d^{n-1}\sigma_{T_j},\varphi_{T_i}),
  \qquad
  \mathbb M^n
  :=\operatorname{diag}(1/|T_i|).
\]
The eigenvalue problem for $A_h^n$ is therefore
$\mathbb A^n\mathsf u=\lambda_h\mathbb M^n\mathsf u$.
Hence its eigenvalues are those of
$\hat{\mathbb A}:=(\mathbb M^n)^{-1}\mathbb A^n$.
Since $\varphi_{T_i}|_{T_i}=\operatorname{vol}_{T_i}/|T_i|$,
\begin{equation}\label{eq:top_matrix_entries}
  \hat{\mathbb A}_{ij}
  =|T_i|(d^{n-1}\sigma_{T_j},\varphi_{T_i})
  =\int_{T_i}d^{n-1}\sigma_{T_j}.
\end{equation}

We determine the signs of these entries from the facet fluxes.
Since the facet mass is diagonal and $\varphi_T$ is supported
on $T$, \eqref{eq:auxvar} implies that the facet degrees of
freedom of $\sigma_T$ vanish outside $\partial T$.
For a facet $F\subset\partial T$, orient $F$ outward from $T$
and let $\mathbb M_F=(\varphi_F,\varphi_F)_h>0$.
Testing \eqref{eq:auxvar} with $\varphi_F$, applying
\cref{lem:dof_flux} with $k=n-1$ and $S=T$, and using
$d\psi_T=\bm\beta_T^\flat$ and Stokes' theorem, we obtain
\[
\begin{aligned}
\mathbb M_F\int_F\sigma_T
&=(\varphi_T,\tilde J_h^{n-1}\varphi_F)\\
&=\frac{\mathcal E_T}{|T|}
  \int_T e^{\psi_T/\eta}
  \bigl(\eta\,d\varphi_F+\bm\beta_T^\flat\wedge\varphi_F\bigr)\\
&=\frac{\eta\mathcal E_T}{|T|}
  \int_{\partial T}
  e^{\psi_T/\eta}\operatorname{tr}\varphi_F\\
&=\frac{\eta\mathcal E_T}{|T|}
  \fint_F e^{\psi_T/\eta}>0.
\end{aligned}
\]
Here we used that
$\operatorname{tr}_F\varphi_F=\operatorname{vol}_F/|F|$
in the chosen orientation, while its trace vanishes on
the other facets of $T$.

Now suppose that $T_i\ne T_j$ share a facet $F$, oriented
outward from $T_j$. The induced orientation from $T_i$
is opposite, so \eqref{eq:top_matrix_entries} gives
\[
  \hat{\mathbb A}_{ij}
  =\int_{\partial T_i}\operatorname{tr}\sigma_{T_j}
  =-\int_F\sigma_{T_j}
  =-\frac{\eta\mathcal E_{T_j}}{\mathbb M_F|T_j|}
    \fint_F e^{\psi_{T_j}/\eta}<0.
\]
If $T_i$ and $T_j$ do not share a facet, then
$\hat{\mathbb A}_{ij}=0$.
Thus $\hat{\mathbb A}$ has nonpositive off-diagonal entries
and a symmetric sparsity pattern.

Moreover, since $M$ is periodic, Stokes' theorem
yields, for every column $j$,
\[
  \sum_i\hat{\mathbb A}_{ij}
  =\sum_i\int_{T_i}d^{n-1}\sigma_{T_j}
  =\int_M d^{n-1}\sigma_{T_j}
  =0.
\]
Since $M$ is connected and the triangulation is conforming,
the cell adjacency graph is connected. The strictly negative
entries in both directions between neighboring cells
therefore imply that $\hat{\mathbb A}$ is irreducible.

By \cref{lem:pf}, zero is algebraically simple and all
remaining eigenvalues have positive real parts.
This proves Assumption~\ref{asm:spectrum} for $k=n$.

\end{proof}

\begin{remark}
The main advantage of the exponential-fitting discretization is that it establishes the spectral property (Assumption \ref{asm:spectrum}) for every $\eta>0$ and every $h>0$. The uniform elemental bound in \cref{lem:elemental} does not directly establish Assumption \ref{asm:stab} for the cases $k=0,n-1$ of \eqref{eq:Jh_unified} due to mass-lumping operation. For these cases, we do not explicitly verify the $L^2$-boundedness uniformly in $h$ and $\eta$, since this condition is mainly needed for \cref{thm:all}. Nevertheless, for each fixed $h$ and $\eta$, the resulting operator is a bounded linear operator on a finite-dimensional space. Hence this possible lack of uniform boundedness does not affect the finite element antidynamo result in \cref{thm:antidynamo2d}.
\end{remark}

\begin{remark}
If an elementwise constant approximation is also used for $k=0$ to construct the
flux operator \cite{wu2020simplex}, the exponential-fitting factors may differ between
elements sharing an edge. In that case, the double sum over elements
and edges cannot generally be reduced to \eqref{eq:k0expeigen} using
only the assembled weights $\omega_E=\sum_{T\supset E}\omega_E^T$. Instead, nonpositive off-diagonal entries require the corresponding
element-dependent weighted sums of \(\omega_E^T\) to be nonnegative.
In this case, a sufficient condition is
\[
  \omega_E^T\ge0
  \qquad
  \text{for every element \(T\) and edge \(E\subset T\)}.
\]
This is stronger than the assembled condition \(\omega_E\ge0\)
used here, although it is not necessary: negative elemental
contributions may still be compensated by positive ones after
weighted assembly.
\end{remark}

\section{Numerical experiments}\label{sec:experiments}

We present simple numerical experiments to illustrate the importance of an
appropriate spatial discretization. Although a complete analysis of the
fully discrete problem requires additional work, backward Euler is
compatible with the $0$- and $n$-form problems whenever the spatial
discretization satisfies the spectral preservation requirement. Motivated by this, we use backward
Euler in the following two-dimensional $1$-form problems to explore the
influence of the spatial discretization. We compute
\begin{equation}\label{eq:time_exp_be}
 u_h^{m+1}=(I+\Delta t A_h^1)^{-1}u_h^m,
 \qquad t_m=m\Delta t,
\end{equation}
where $\Delta t$ denotes the time-step size.

We consider the flat torus $M=[0,2\pi)^2$ with periodic boundary
conditions. We use the lowest-order Whitney finite elements on
acute simplicial triangulations satisfying the mesh
conditions, and denote the mesh resolution by $h=2\pi/N$.
All methods on a given mesh start from the interpolant of
\begin{equation}\label{eq:time_exp_initial}
 \begin{aligned}
 \omega_0={}&(0.6+0.25\cos x+0.15\sin(x+y))\,dx\\
 &+(0.6+0.2\sin y+0.1\cos(x-y))\,dy.
 \end{aligned}
\end{equation}
We plot the normalized $L^2$ norm
\begin{equation}\label{eq:time_exp_norm}
 R_0(t_m)=\frac{\|u_h^m\|}{\|u_h^0\|}.
\end{equation}

\subsection{Weighted Galerkin evolution}
\label{subsec:time_exp_weighted}

We first compare the weighted-inner-product Galerkin discretization with
the standard discretization, which uses the $L^2$-projected contraction.
We take $\eta=0.5$, $N\in\{8,16,32\}$, and $\Delta t=0.1$ on
$0\le t\le100$. The advection field is a potential flow
\[
 \bm\beta=(2\sin2x,2\sin2y)=\nabla\Phi,\qquad
 \Phi=-\cos2x-\cos2y,\qquad
 \theta=\Phi/\eta,\qquad \rho=e^{-\theta}.
\]
The weighted method is the Galerkin discretization \eqref{eq:wadd} with
the inner product weighted by $\rho$. In addition to $\|\cdot\|$, we
measure both solutions using
\begin{equation}\label{eq:time_exp_weighted_metrics}
 \|v\|_\theta^2=\int_M\rho v\wedge\star v,
 \qquad
 R_\theta(t_m)=\frac{\|u_h^m\|_\theta}{\|u_h^0\|_\theta}.
\end{equation}

\begin{figure}[tbp]
 \centering
 \includegraphics[width=\linewidth]{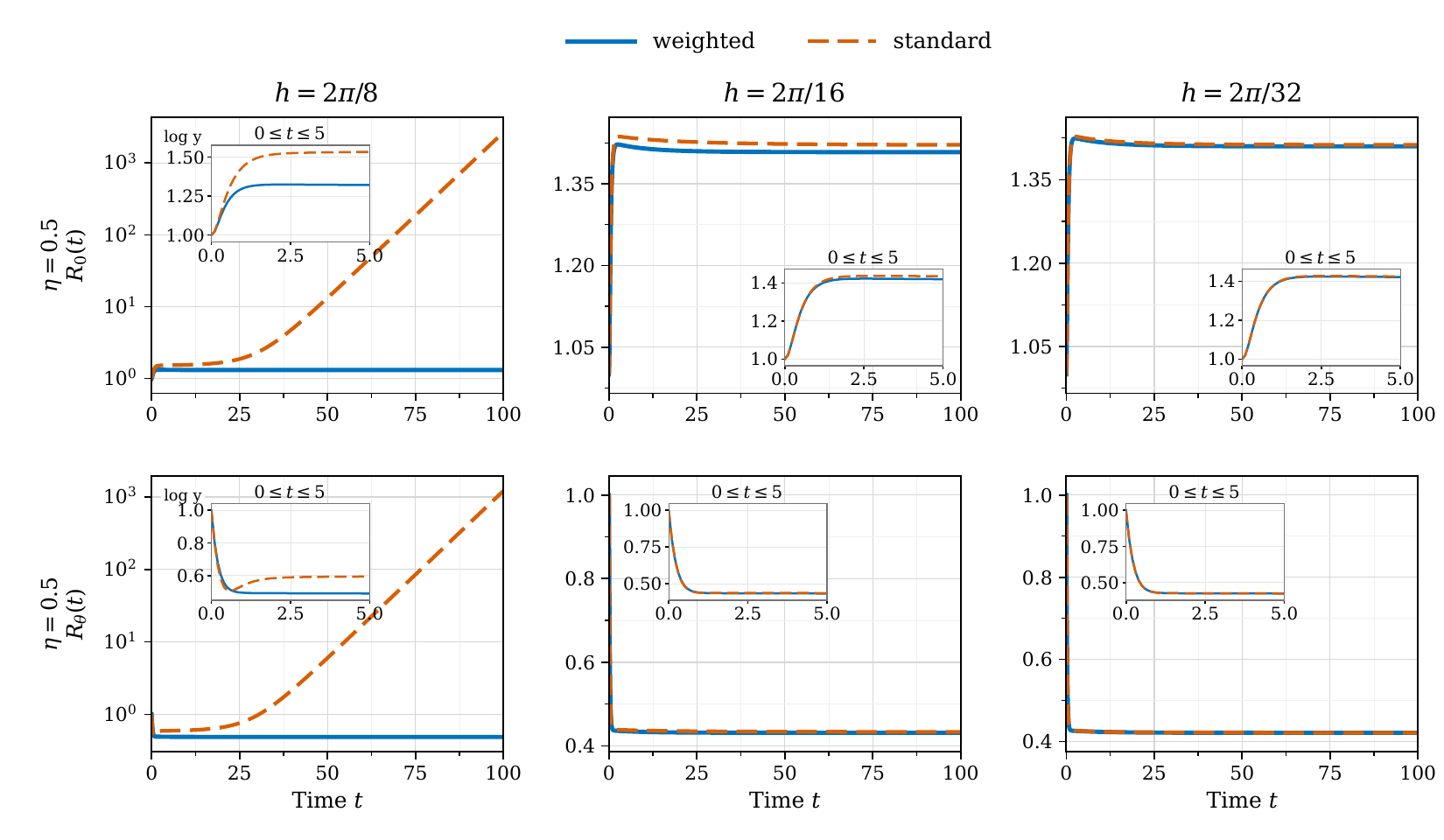}
 \caption{Weighted Galerkin and standard evolution for the potential
 field with $\eta=0.5$ and $\Delta t=0.1$. Columns correspond to
 $N=8,16,32$. The top row shows the normalized $L^2$ norm $R_0$,
 and the bottom row shows the normalized weighted norm $R_\theta$.
 The main panels cover $0\le t\le100$; the insets show $0\le t\le5$.}
 \label{fig:time_exp_weighted}
\end{figure}

\Cref{fig:time_exp_weighted} shows the evolution of the two norms.
On the coarsest mesh, the standard discretization exhibits growth in
the weighted norm after an initial decrease. This effect is
reduced under mesh refinement, and the two discretizations agree more
closely on the finest mesh. The mesh dependence is consistent with the
need for sufficiently small $h$ in \cref{lm:spectrum_small_h}.

In contrast, the weighted discretization gives a non-increasing weighted
norm on all three meshes. This agrees with the weighted dissipativity
in \cref{thm:weighted_antidynamo}, which is also preserved by backward
Euler. The $L^2$ norm need not decrease monotonically and may
exhibit transient growth, while the discrete solution still converges
exponentially to a stationary state.

\subsection{Comparison between exponential-fitting and standard discretizations}
\label{subsec:time_exp_mesh}

We next compare the {exponential-fitting} discretization with the
standard method, using the non-potential field
\[
 \bm\beta=(2\sin2x+0.5\sin y,\;2\sin2y-0.5\sin x).
\]
We consider $\eta\in\{0.5,0.1\}$ on the meshes
$N\in\{8,16,32\}$, with $\Delta t=0.1$ and $0\le t\le100$.

\begin{figure}[tbp]
 \centering
 \includegraphics[width=\linewidth]{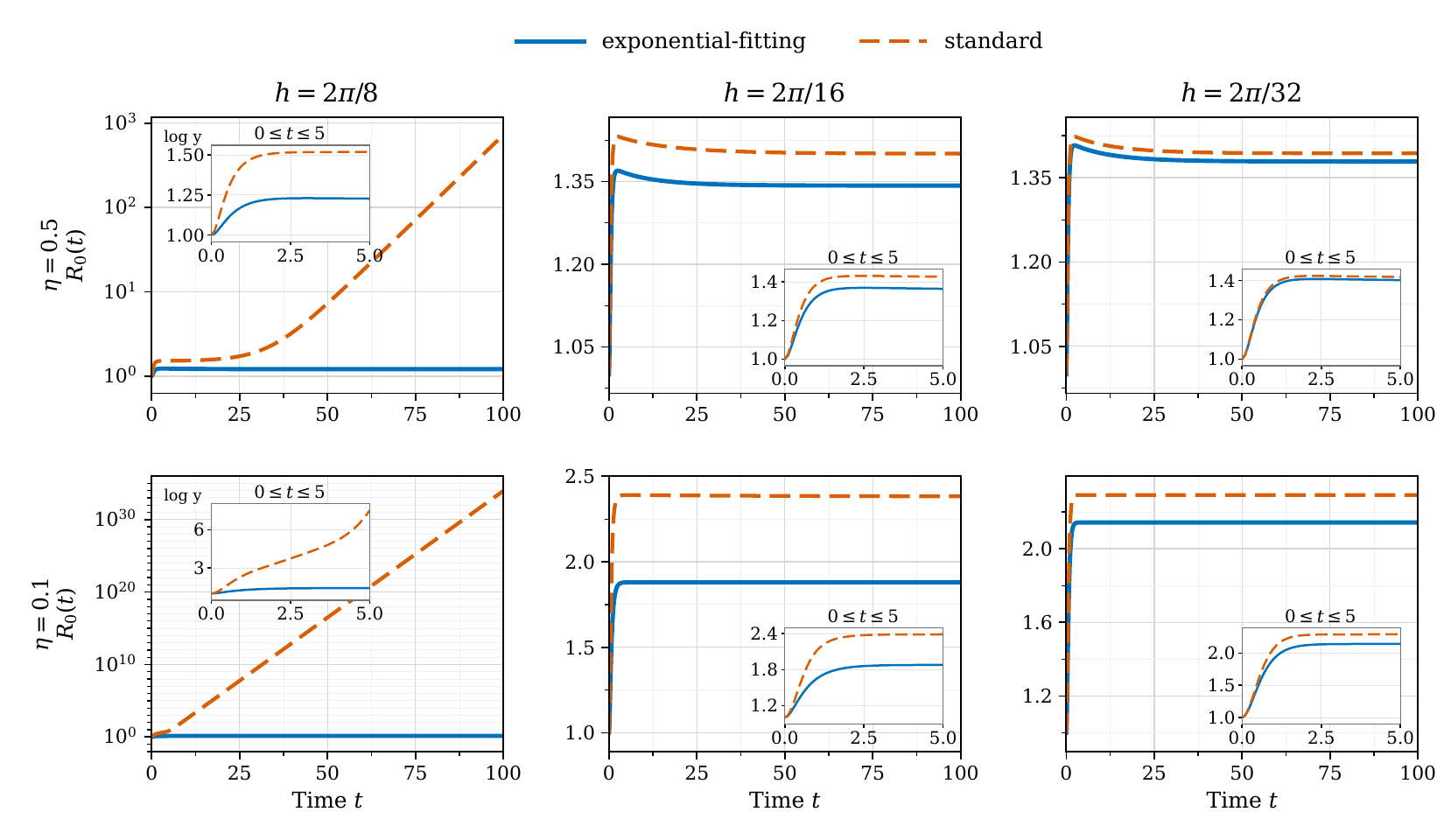}
 \caption{Exponential-fitting and standard evolution of the normalized
 $L^2$ norm $R_0$ with $\Delta t=0.1$. Rows correspond to
 $\eta=0.5$ and $\eta=0.1$; columns correspond to $N=8,16,32$.
 The main panels cover $0\le t\le100$, and the insets show
 $0\le t\le5$.}
 \label{fig:time_exp_mesh}
\end{figure}

\Cref{fig:time_exp_mesh} shows a separation between the two
trajectories on the coarsest mesh, where the standard discretization
exhibits substantial growth that becomes more pronounced as the diffusion
coefficient decreases.  This large growth is not observed on the two finer meshes over the
displayed interval. For both diffusion coefficients, refinement brings
the two trajectories closer. In particular, the exponential-fitting discretization produces lower $L^2$ norms than the standard discretization on all three meshes. This behavior
may reflect additional numerical diffusion associated with
exponential fitting. 

It is worth emphasizing that exponential fitting was originally
developed and has been widely used to stabilize advection--diffusion
problems. Here our motivation is its ability to preserve the correct
spectral distribution for the $0$- and $n$-form operators without a
mesh size restriction, under some mesh conditions. This spectral preservation ensures that the
discretization reproduces the correct qualitative long-time
asymptotic behavior predicted by the continuous theory, even on coarse meshes and
with a small diffusion coefficient.

\subsection{Influence of mass lumping}
\label{subsec:time_exp_lumping}

Finally, we examine whether exponential fitting alone is sufficient
when the mass lumpings are not used. Here \emph{nodal lumping}
refers to the mass on $0$-forms, whereas \emph{facet lumping} refers to
the auxiliary flux mass in the $2$-form discretization.

We use two fixed smooth periodic velocity fields, each at
$\eta\in\{0.005,0.001\}$. The first is the trigonometric polynomial
\begin{equation}\label{eq:time_exp_velocity_i}
 \bm\beta_{\mathrm I}(x,y)=\bm c+
 \sum_{\bm k\in\mathcal K}
 \bigl[\bm a_{\bm k}\cos(k_1x+k_2y)
       +\bm b_{\bm k}\sin(k_1x+k_2y)\bigr],
\end{equation}
where
\[
 \mathcal K=\{(0,j):1\le j\le8\}
 \cup\{(i,j):1\le i\le4,\ -8\le j\le8\}.
\]

\begin{figure}[tbp]
 \centering
 \includegraphics[width=\linewidth]{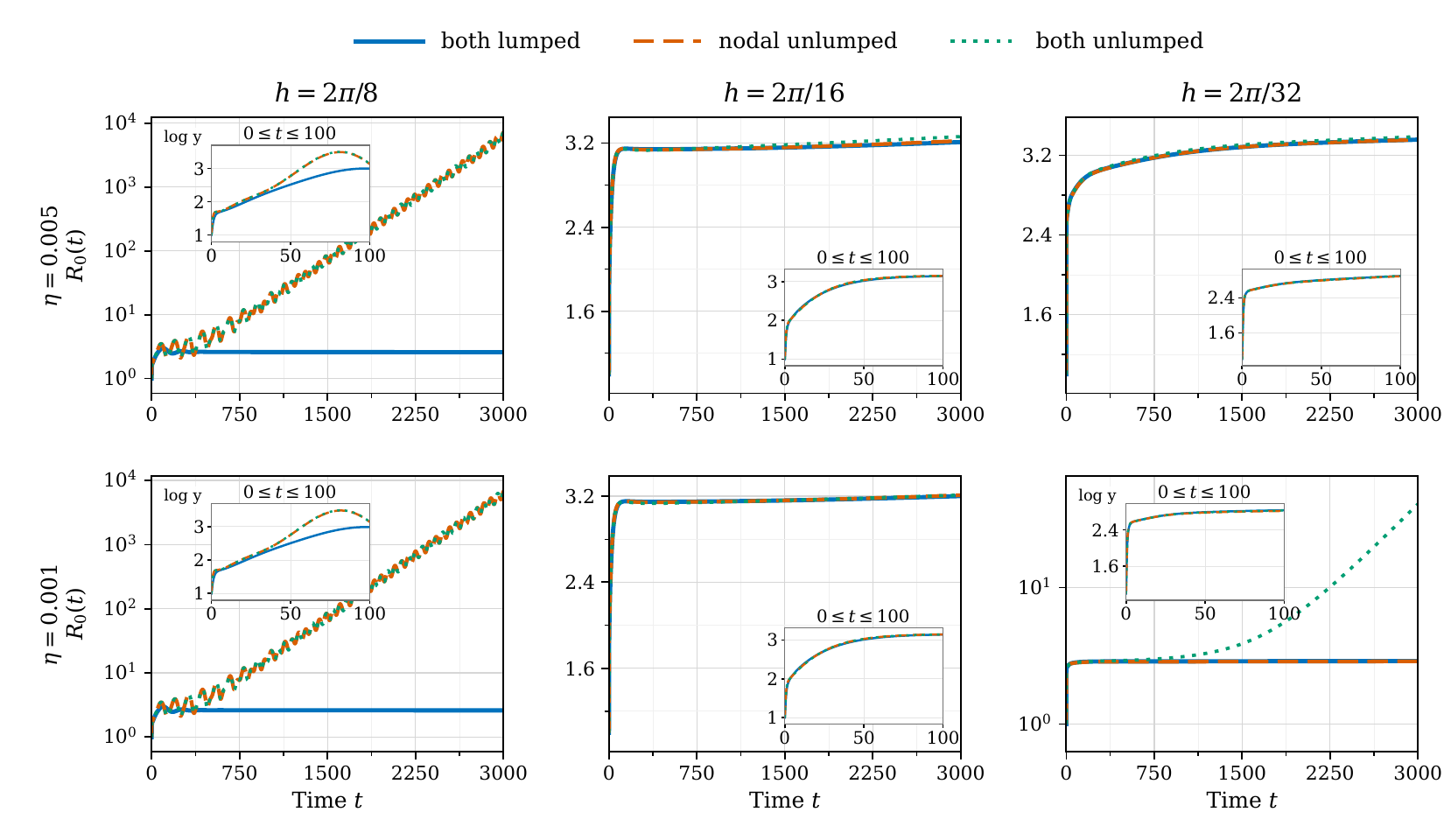}
 \caption{Effect of nodal lumping for velocity
 $\bm\beta_{\mathrm I}$.
 Rows correspond to $\eta=0.005$ and $\eta=0.001$; columns correspond
 to $N=8,16,32$. The main panels cover $0\le t\le3000$ with
 $\Delta t=0.1$; the insets show $0\le t\le100$.}
 \label{fig:time_exp_lumping_i}
\end{figure}

The second field has the simpler form
\begin{equation}\label{eq:time_exp_velocity_ii}
 \bm\beta_{\mathrm{II}}(x,y)=3\mathcal R(0.78)
 \begin{pmatrix}
  \sin(5x+0.85)\cos(5y+2.95)\\
  \cos(5x+0.85)\sin(5y+2.95)
 \end{pmatrix},
 \quad
 \mathcal R(\alpha)=
 \begin{pmatrix}\cos\alpha&-\sin\alpha\\
                 \sin\alpha& \cos\alpha\end{pmatrix}.
\end{equation}
For this field, we compare the results obtained with and without mass lumping.

\begin{figure}[tbp]
 \centering
 \includegraphics[width=\linewidth]{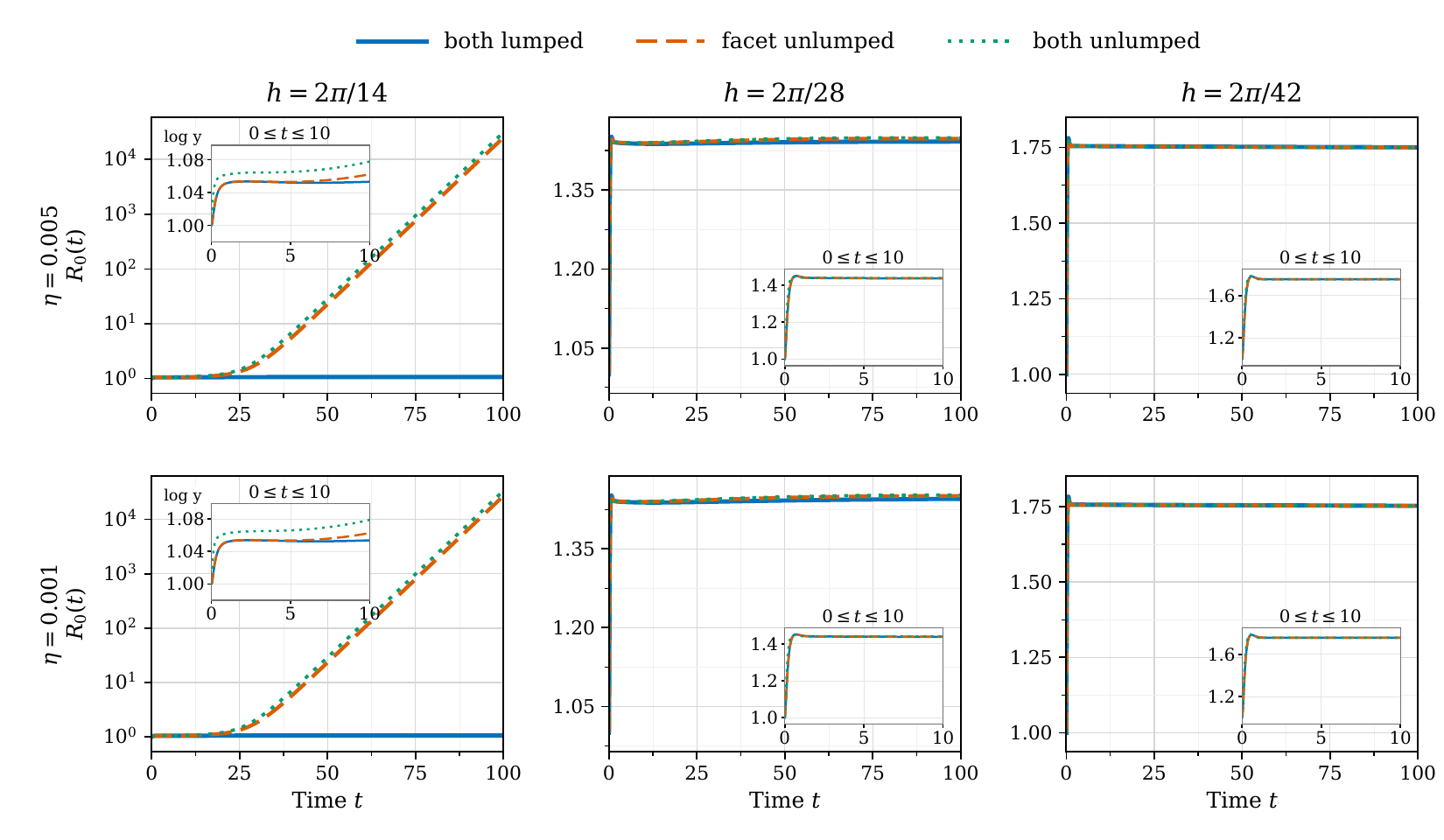}
 \caption{Effect of facet lumping for velocity
 $\bm\beta_{\mathrm{II}}$. Rows correspond to $\eta=0.005$ and
 $\eta=0.001$; columns correspond to $N=14,28,42$.
 The main panels show $R_0$ on $0\le t\le100$ with
 $\Delta t=0.025$; the insets show $0\le t\le10$.}
 \label{fig:time_exp_lumping_ii}
\end{figure}

\Cref{fig:time_exp_lumping_i,fig:time_exp_lumping_ii} show that either
partial removal of lumping can lead to sustained spurious growth,
depending on the velocity field. For $\bm\beta_{\mathrm I}$, removing
nodal lumping produces large oscillatory growth on the coarsest mesh
at both diffusion coefficients. For $\bm\beta_{\mathrm{II}}$, removing
facet lumping produces substantial growth on the coarsest mesh within
a much shorter time interval. Removing both lumpings also
produces growth in these coarse-mesh tests. In contrast, the fully
lumped trajectories remain of order one throughout the displayed
intervals. The early-time insets show that an initial increase in the
physical norm can occur while this transient
effect is distinct from the sustained growth of the unlumped variants.

We compute the spectrum of the discrete operator, which supports this distinction. On the
coarsest mesh for $\bm\beta_{\mathrm I}$, the nodal-unlumped
$0$-form operator has a conjugate pair with real part approximately
$-3.10\times10^{-3}$ at both diffusion coefficients. For
$\bm\beta_{\mathrm{II}}$, the facet-unlumped $2$-form operator has a
real eigenvalue approximately $-0.1400$ at $\eta=0.005$ and $-0.1408$
at $\eta=0.001$. The corresponding lumped operators have no
eigenvalues with negative real part.

These comparisons highlight the importance of mass lumping in the
exponential-fitting discretization: replacing the lumped mass matrices
by consistent ones need not preserve the same qualitative behavior.
Although mesh refinement may prevent such growth, its effect need not
be monotone in $h$. For example, for $\bm\beta_{\mathrm I}$ at
$\eta=0.001$, the variant without any lumping shows growth at $N=32$,
while none of the three variants shows strong growth at $N=16$
over the displayed time interval.

\subsection{Stationary kernels and dynamo action}
\label{subsec:dynamo3d}

The lower bound in \cref{thm:betti} allows the dimension of the
stationary space to exceed the Betti number. We investigate this
possibility numerically and examine its relation to dynamo action.
For the induction problem \eqref{eq:intro_induction_vector}, magnetic
fields are represented by closed $1$-forms in two dimensions and
closed $2$-forms in three dimensions. For a mode $u_h(t)=e^{s_ht}b_h$, the discrete evolution gives
\[
 A_h^k b_h=-s_hb_h, \qquad s_h\in\mathbb C.
\]
Thus, $\operatorname{Re}s_h>0$ indicates exponential growth,
$\operatorname{Re}s_h<0$ indicates decay, and $s_h=0$ corresponds
to a stationary mode.
In all examples, we compare
\[
 \dim\ker A_h^k,\qquad
 \dim\mathcal S_h^k=\dim(\ker A_h^k\cap\mathfrak Z_h^k),\qquad
 \dim(\ker A_h^k\cap\mathfrak B_h^k).
\]
Numerical zeros are counted by singular values in the
$L^2$ metric, using the common threshold
$10^{-12}\|A_h^k\|_{L^2\to L^2}$.

\subsubsection{Three-dimensional comparison}
\label{subsec:antidynamo2d_test}
We consider the Roberts flow \cite{roberts1972dynamo}
$\bm\beta=(\cos y,\sin x,\cos x+\sin y)$ on
$\mathbb T^3=[0,2\pi)^3$, with $b_2=3$, using the standard
lowest-order FEEC discretization with $h=2\pi/N$.
For each mesh, we seek positive values of the diffusion parameter
$\eta$ for which the stationary problem
\[
 A_h^2(\eta)b_h=0,
 \qquad 0\ne b_h\in\mathfrak B_h^2,
\]
admits a solution. We call these values stationary parameter roots
and select the largest one (denoted as $\eta_h^*$). We then track the corresponding eigenspace at
$\eta=0.99\eta_h^*$ and $1.01\eta_h^*$.

As shown in \cref{tab:dynamo3d}, the additional stationary modes are exact and satisfy the
magnetic closedness constraint.
At $\eta=0.99\eta_h^*$, these modes have positive growth rates,
demonstrating discrete dynamo action even though the stationary
dimension has returned to $b_2$.
At $\eta=1.01\eta_h^*$, the same modes decay.
The stationary crossing therefore marks a change from decay to
growth along this branch as $R_m=\eta^{-1}$ increases. These observations connect strictness of the topological lower bound with a transition between growth and decay of the selected magnetic modes. A similar mechanism is illustrated by Arnold's explicit example in \cite{arnold2013some}.

\begin{table}[htbp]
 \centering
 \caption{Three-dimensional Roberts flow numerical tests.}
 \label{tab:dynamo3d}
 \setlength{\tabcolsep}{5pt}
 \begin{tabular}{@{}crccccr@{}}
  \toprule
  $N$ & $\eta_h^*$ & $\eta/\eta_h^*$
      & $\dim\ker A_h^2$ & $\dim\mathcal S_h^2$
      & $\dim(\ker A_h^2\cap\mathfrak B_h^2)$ & $\operatorname{Re}s_h$ \\
  \midrule
  $4$ & $0.17070714$ & $0.99$ & $3$ & $3$ & $0$ & $ 0.0037064$ \\
      &              & $1$    & $5$ & $5$ & $2$ & $ 0$ \\
      &              & $1.01$ & $3$ & $3$ & $0$ & $-0.0036906$ \\
  \addlinespace[2pt]
  $6$ & $0.16816320$ & $0.99$ & $3$ & $3$ & $0$ & $ 0.0027894$ \\
      &              & $1$    & $5$ & $5$ & $2$ & $ 0$ \\
      &              & $1.01$ & $3$ & $3$ & $0$ & $-0.0027873$ \\
  \addlinespace[2pt]
  $8$ & $0.16832090$ & $0.99$ & $3$ & $3$ & $0$ & $ 0.0026720$ \\
      &              & $1$    & $5$ & $5$ & $2$ & $ 0$ \\
      &              & $1.01$ & $3$ & $3$ & $0$ & $-0.0026723$ \\
  \bottomrule
 \end{tabular}
\end{table}

\subsubsection{Two-dimensional comparison}
\label{subsec:antidynamo2d_test}

As a planar counterpart, we take $\bm\beta=(\cos y,\sin x)$ on
$\mathbb T^2=[0,2\pi)^2$, with $b_1=2$.
We compare the standard and exponential-fitting
methods on acute periodic meshes, with $h=2\pi/N$.
For the discretization on the meshes $N=8,16,32$,
no positive diffusivity $\eta$ was found for which a nonzero
stationary mode exists in $\mathfrak B_h^1$.

The standard discretization nevertheless admits additional
stationary modes.
\Cref{tab:antidynamo2d_test} reports one representative critical
point selected independently on each of the meshes $N=16$ and $32$. The reported $s_h$ is the growth rate of the additional mode
that becomes stationary at $\eta_h^*$.
This mode grows at $\eta=0.99\eta_h^*$ and decays at
$\eta=1.01\eta_h^*$. The closed stationary space retains dimension $b_1=2$
at all three parameter values.

For the exponential-fitting method, the dimensions remain $(2,2,0)$
throughout the parameter search.
The table reports the largest real part among the nonzero
full-space growth rates, so a negative value means that every
nonstationary eigenmode decays.
The reported values are negative on both meshes, consistent
with the discrete antidynamo result, which guarantees this
behavior for every $\eta>0$.

\begin{table}[htbp]
 \centering
 \small
 \caption{Two-dimensional comparison.}
 \label{tab:antidynamo2d_test}
 \setlength{\tabcolsep}{3pt}
 \begin{tabular}{@{}crlccccr@{}}
  \toprule
  $N$ & $\eta_h^*$ & Method & $\eta/\eta_h^*$
      & $\dim\ker A_h^1$ & $\dim\mathcal S_h^1$
      & $\dim(\ker A_h^1\cap\mathfrak B_h^1)$
      & $\operatorname{Re}s_h$ \\
  \midrule
  $16$ & $0.00278940$ & Standard
       & $0.99$ & $2$ & $2$ & $0$ & $0.00234894$ \\
       & & & $1$ & $3$ & $2$ & $0$ & $0$ \\
       & & & $1.01$ & $2$ & $2$ & $0$ & $-0.00234419$ \\
  \addlinespace[2pt]
  \cmidrule(lr){3-8}
       & & Exp-Fitted
       & $0.99$ & $2$ & $2$ & $0$ & $-0.05912090$ \\
       & & & $1$ & $2$ & $2$ & $0$ & $-0.05966233$ \\
       & & & $1.01$ & $2$ & $2$ & $0$ & $-0.06020200$ \\
  \midrule
  \midrule
  $32$ & $0.00127238$ & Standard
       & $0.99$ & $2$ & $2$ & $0$ & $0.00538226$ \\
       & & & $1$ & $3$ & $2$ & $0$ & $0$ \\
       & & & $1.01$ & $2$ & $2$ & $0$ & $-0.00539572$ \\
  \addlinespace[2pt]
  \cmidrule(lr){3-8}
       & & Exp-Fitted
       & $0.99$ & $2$ & $2$ & $0$ & $-0.03250329$ \\
       & & & $1$ & $2$ & $2$ & $0$ & $-0.03250333$ \\
       & & & $1.01$ & $2$ & $2$ & $0$ & $-0.03250339$ \\
  \bottomrule
 \end{tabular}
\end{table}

For the standard method, spectral preservation is guaranteed for
sufficiently small $h$ when $\eta$ and $\bm\beta$ are fixed.
The admissible set of mesh size may depend on these parameters, so the
result does not provide a uniform guarantee as they vary.
In particular, even on a fixed fine mesh, spurious spectrum and artificial growth may occur as the diffusivity is reduced. Under the stated mesh assumptions, the exponential-fitting
discretization satisfies the spectral preservation assumption
for every $h>0$ and $\eta>0$, making it more robust in
advection-dominated regimes.

\section{Conclusion}\label{sec:conclude}
In this paper, we developed a finite element framework for the
long-time analysis of Lie advection--diffusion equations of differential
forms. The main objective was to go beyond finite-time convergence as
\(h\to0\), and to understand whether a fixed finite element discretization
preserves the structural asymptotic behavior of the continuous system. Under a spectral preservation assumption for \(0\)- and
\(n\)-forms, we established long-time asymptotic results for arbitrary positive diffusion coefficients of advection--diffusion equations. In
particular, on two-dimensional compact manifolds, the discrete \(1\)-form
evolution converges to a stationary closed form for every \(\eta>0\).

We also verified the structural assumptions in representative cases. For
sufficiently fine meshes, the spectral preservation follows from
spectral approximation theory. For lowest-order Whitney forms, we adopt an exponential-fitting discretization motivated by the Witten deformation and the edge-averaged finite element method. Under suitable mesh conditions, the resulting operators possess the structure of singular irreducible $M$-matrices for arbitrary mesh sizes.

Several questions remain open. The first concerns boundary conditions.
As noted by Arnold and Khesin~\cite{arnold1998topological}, \textit{``the
preceding theory has not been settled for manifolds with boundary, though
it certainly deserves to be.''} Although \Cref{sec:asymptoticwithboundary} extends our discrete
results to absolute and relative boundary conditions through the
underlying de Rham complex structure, this extension is primarily
structural. A complete analysis connecting these discrete results
with the well-posedness, spectral properties, and long-time
asymptotic behavior of the corresponding continuous boundary-value
problems remains to be developed. In particular, extending the
theory to more general boundary conditions, especially for degrees $0<k<n$, is an important direction
for future research. A second direction is to identify broader classes of spatial discretizations that satisfy the structural assumptions introduced in this paper. In particular, semi-Lagrangian methods represent a promising candidate: since the Lie derivative is generated by pullback along the flow, a semi-Lagrangian FEEC method based on discrete pullbacks could naturally preserve the underlying cochain structure. 
These questions point toward a broader global numerical analysis of Lie advection--diffusion and its applications in topological hydrodynamics and magnetohydrodynamics.

\section*{Declaration of AI Use}
During the preparation of this work, the authors used ChatGPT 6 Pro to set up the numerical experiments for determining some critical velocity fields. In addition, ChatGPT 5.6 Sol and Claude Fable 5 were used to refine the exposition, check mathematical arguments and proofs for potential errors or typographical issues, and suggest improvements to their presentation. The authors take full responsibility for the mathematical content, numerical results, and final manuscript.

\appendix

\section{Auxiliary proofs}\label{apx:proofs}

\subsection{Proof of the scalar spectral theorem}
\label{apx:spectrum}

Throughout, $M$ is a smooth, compact, connected Riemannian
manifold, either closed or with smooth boundary, and the periodic
case is identified with $M=\mathbb T^n$. We prove \cref{thm:spec} for $\bm\beta\in C^{0,1}(M;TM)$ in the closed, periodic, and Neumann cases, where the classical strong
maximum principle and the Hopf lemma apply directly. The main
simplification is that the principal eigenfunction is explicit,
$e_0=\mathbf 1$ with $\lambda_0=0$, so no Krein--Rutman argument is needed. Under Dirichlet
conditions, constants are not admissible and the principal eigenpair is
instead constructed using Krein--Rutman theory
\cite{krein1948linear}, see also
\cite{berestycki1994principal}. For $\bm\beta\in L^\infty$ the conclusions also persist. The pointwise arguments are replaced by De~Giorgi--Nash regularity, the weak Harnack inequality, and the strong maximum principle for $H^1$ sub- and supersolutions \cite[Thms.~8.18--8.19]{gilbarg1977elliptic}. In the presence of a
boundary, where the classical Hopf argument is not directly available,
one may instead apply the Krein--Rutman theory on \(L^2(M)\)
to compact irreducible positive resolvents and semigroups
\cite{depagter1986irreducible}
\cite[Ch.~V, \S5]{schaefer1974banach}. For brevity, we write $L:=\Delta_{\eta,\bm\beta}^0
  =-\eta\Delta_g+\bm\beta\cdot\nabla$, where \(\Delta_g\) denotes the Laplace--Beltrami operator.

\begin{proof}[Proof of \cref{thm:spec} in the closed and Neumann cases]

On \(V=H^1(M)\), define
\[
  a(u,v)
  :=
  \eta(\nabla u,\nabla v)
  +(\bm\beta\cdot\nabla u,v).
\]
Writing \(B_\infty:=\|\bm\beta\|_{L^\infty(M)}\), for every
\(\varepsilon>0\),
\[
  |(\bm\beta\cdot\nabla u,u)|
  \le
  \varepsilon\eta\|\nabla u\|^2
  +
  \frac{B_\infty^2}{4\varepsilon\eta}\|u\|^2.
\]
Hence \(a\) is a closed sectorial form, and Kato's first representation
theorem \cite[Ch.~VI, Thm.~2.1]{kato1995perturbation} yields the sectorial
operator \(L\). Standard elliptic regularity gives
\[
  \mathcal D(L)
  =
  \begin{cases}
    H^2(M), & \partial M=\varnothing,\\[1mm]
    \{u\in H^2(M):\partial_\nu u=0\}, & \partial M\neq\varnothing.
  \end{cases}
\]
Since \(H^1(M)\hookrightarrow\hookrightarrow L^2(M)\), \(L\) has compact
resolvent. Hence its spectrum consists of isolated eigenvalues of finite
algebraic multiplicity.

Let \(Lu=\lambda u\), \(\|u\|_{L^2}=1\), and
\(X:=\|\nabla u\|_{L^2}\). Then
\[
  \operatorname{Re}\lambda
  \ge
  \eta X^2-B_\infty X
  \ge
  \frac{\eta}{2}X^2-\frac{B_\infty^2}{2\eta},
  \qquad
  |\operatorname{Im}\lambda|
  \le B_\infty X.
\]
Thus
\begin{equation}
\label{eq:scalar_spectral_parabola}
  |\operatorname{Im}\lambda|^2
  \le
  c_1(\operatorname{Re}\lambda+c_2),
  \qquad
  c_1:=\frac{2B_\infty^2}{\eta},
  \quad
  c_2:=\frac{B_\infty^2}{2\eta}.
\end{equation}

Since \(L\mathbf1=0\), \(0\in\sigma(L)\). Let
\[
  Lu=\lambda u,
  \qquad
  u\neq0,
\]
and set \(w:=|u|^2\). Since $\bm\beta\in C^{0,1}\subset C^{0,\alpha}$, elliptic regularity
upgrades the eigenfunction to a classical solution, so
that $w\in C^{2}(\overline M)$ and the maximum principle applies
pointwise.
\[
  -\eta\Delta_g w+\bm\beta\cdot\nabla w
  =
  2\operatorname{Re}\lambda\,w
  -2\eta|\nabla u|^2.
\]
In the Neumann case,
\[
  \partial_\nu w=0
  \qquad\text{on }\partial M.
\]
If \(\operatorname{Re}\lambda\le0\), then
\[
  \eta\Delta_g w-\bm\beta\cdot\nabla w\ge0.
\]
The strong maximum principle, together with the Hopf lemma in the Neumann
case, implies that \(w\) is constant. Substituting back gives
\[
  \operatorname{Re}\lambda=0,
  \qquad
  \nabla u=0,
\]
so \(u\) is constant and therefore \(\lambda=0\). Hence
\begin{equation}
\label{eq:strict_positive_real_part}
  \mathcal N(L)=\operatorname{span}\{\mathbf1\},
  \qquad
  \operatorname{Re}\lambda>0
  \quad
  \forall\lambda\in\sigma(L)\setminus\{0\}.
\end{equation}

Suppose that \(0\) is not algebraically simple. Since
\(\mathcal N(L)=\operatorname{span}\{\mathbf1\}\), after rescaling there is
\(v\in\mathcal D(L)\) such that
\[
  Lv=\mathbf1.
\]
Because the coefficients are real, \(q:=\operatorname{Re}v\) satisfies
\[
  Lq=1,
\]
with \(\partial_\nu q=0\) in the Neumann case. At an interior minimum of
\(q\),
\[
  1=Lq=-\eta\Delta_gq\le0,
\]
which is a contradiction. In the Neumann case, if the minimum occurs on the
boundary, the Hopf lemma applied to \(-q\) gives
\(\partial_\nu q<0\), again contradicting the boundary condition.
Thus \(0\) is algebraically simple. Consequently,
\[
  \lambda_0=0,
  \qquad
  e_0=\mathbf1
\]
is the real, algebraically simple principal eigenpair.

Set
\[
  K
  :=
  \sigma(L)\cap
  \{z\in\mathbb C:0\le\operatorname{Re}z\le1\}.
\]
By \eqref{eq:scalar_spectral_parabola}, \(K\) is bounded, and compact
resolvent implies that \(K\) is finite. Therefore
\[
  \mu
  :=
  \min\left\{
    1,\,
    \min_{\lambda\in K\setminus\{0\}}
    \operatorname{Re}\lambda
  \right\}
  >0.
\]
Hence
\[
  \operatorname{Re}\lambda\ge\mu,
  \qquad
  \lambda\in\sigma(L)\setminus\{0\}.
\]
Combining this with \eqref{eq:scalar_spectral_parabola},
\[
  \frac{|\operatorname{Im}\lambda|}
       {\operatorname{Re}\lambda}
  \le
  \frac{\sqrt{c_1(r+c_2)}}{r},
  \qquad
  r:=\operatorname{Re}\lambda\ge\mu.
\]
Since the right-hand side is bounded on \([\mu,\infty)\), there exists
\(\vartheta<\pi/2\) such that
\[
  |\arg\lambda|\le\vartheta
  \qquad
  \forall\lambda\in\sigma(L)\setminus\{0\}.
\]
Thus
\[
  \sigma(L)
  \subset
  \{0\}\cup
  \{z\in\mathbb C\setminus\{0\}:|\arg z|\le\vartheta\},
\]
which proves the claimed spectral gap and sectorial distribution.
\end{proof}

\subsection{Proof of \cref{lm:expA}}\label{apx:expA}
\begin{proof}[Proof of \cref{lm:expA}]
Since \(X\) is finite-dimensional and \(0\) has algebraic multiplicity \(1\), we have $\dim \mathcal N(A)=1$ and
$\mathcal N(A)\cap\mathcal R(A)=\{0\}$. Hence \(P\) is a rank-one projection. Set $Y:=\mathcal R(A)$. By the rank--nullity theorem,
\[
X=\mathcal N (A)\oplus Y,
\qquad
\sigma(A|_Y)=\sigma(A)\setminus\{0\}.
\]
Now let \(x\in X\). Write uniquely
\[
x=x_0+y,
\qquad x_0\in \mathcal N (A),\ y\in Y.
\]
Since \(Ax_0=0\), we have
\[
e^{-tA}x_0=x_0,
\qquad t\ge 0.
\]
Also, since \(Y\) is \(A\)-invariant,
\[
e^{-tA}x=x_0+e^{-tA|_Y}y.
\]
Set
\[
B:=A|_Y.
\]
Then
\[
\sigma(B)=\sigma(A)\setminus\{0\},
\]
so every eigenvalue \(\mu\in \sigma(B)\) satisfies \(\text{Re}(\mu)>0\). Since \(Y\) is finite-dimensional, the standard estimate for matrix exponentials implies that for every
\[
0<\gamma<\min\{\text{Re}(\mu):\mu\in \sigma(B)\},
\]
there exists \(C>0\) such that
\[
\|e^{-tB}\|_{\mathcal L(Y)}\le Ce^{-\gamma t},
\qquad t\ge 0.
\]
Therefore,
\[
e^{-tA}x=x_0+e^{-tB}y \to x_0
\qquad\text{as }t\to\infty.
\]
Since \(Px=x_0\), this proves that
\[
e^{-tA}x\to Px
\qquad\text{for every }x\in X.
\]
Finally, for \(x=x_0+y\) we have
\[
(e^{-tA}-P)x=e^{-tB}y.
\]
Thus \(e^{-tA}-P\) vanishes on \(\mathcal N (A) \) and coincides with \(e^{-tB}\) on \(Y\). Since \(X=\mathcal N (A)\oplus Y\) and all norms are equivalent on the finite-dimensional space \(X\), it follows that there exists \(C>0\) such that
\[
\|e^{-tA}-P\|_{\mathcal L(X)}\le Ce^{-\gamma t},
\qquad t\ge 0.
\]
This completes the proof.
\end{proof}

\bibliography{ref}
\bibliographystyle{abbrv}
\end{document}